\documentclass[12pt]{article}
 \usepackage{amsmath}
 \usepackage{amsfonts}

\usepackage{graphicx}
\usepackage{color}

\usepackage[left=2.5cm,right=2.5cm,top=2.5cm,bottom=2.5cm]{geometry}
\usepackage[utf8]{inputenc}
\usepackage[T1]{fontenc}
\usepackage{xcolor,graphicx}
\usepackage{hyperref}
\usepackage{subcaption}
\usepackage{tikz}

\hypersetup{%
colorlinks=true,
breaklinks=true,
linkcolor=red,
anchorcolor=black,
citecolor=brown,
filecolor=blue,
menucolor=red,
urlcolor=red,
}
\usepackage[round,longnamesfirst]{natbib}
\usepackage{mathtools, amsfonts, amssymb,amsthm,stmaryrd}
\mathtoolsset{showonlyrefs}

\usepackage{enumerate}
\usepackage{comment}
\usepackage{ifthen}
\usepackage{mathrsfs}

\newtheorem{definition}{Definition}

\newcommand{\EE}{\mathbb{E}}

\newcommand{\PP}{\mathbb{P}}

\newcommand{\Xtemp}[1]{%
\ifthenelse{\equal{#1}{0}}
{X^{(\n0)}}
{
\ifthenelse{\equal{#1}{1}}
{X^{[\n1]}}
{X^{(#1)}}
}}
\newcommand{\hXtemp}[1]{%
\ifthenelse{\equal{#1}{0}}
{\widehat X^{(\n0)}}
{
\ifthenelse{\equal{#1}{1}}
{\widehat X^{[\n1]}}
{\widehat X^{(#1)}}
}}

\newcommand{\Tt}{\boldsymbol t}

\newcommand{\TT}{\boldsymbol T}
\newcommand{\Ss}{\boldsymbol s}
\newcommand{\Uu}{\boldsymbol u}
\newcommand{\Vv}{\boldsymbol v}

\newcommand{\cT}{\mathcal T}
\newcommand{\cU}{\mathcal U}
\newcommand{\Ynm}{{Y^{(j)}_m}}
\newcommand{\Tnm}{{\boldsymbol t^{(j)}_m}}

\newcommand{\enm}{{\varepsilon^{(j)}_m}}
\newcommand{\unm}{{e^{(j)}_m}}
\newcommand{\Xp}[1]{X^{(#1)}}
\newcommand{\Xtp}[1]{\widetilde{X}^{(#1)}}

\newcommand{\Xn}{{\Xtemp{j}}}

\newcommand{\D}{\rm{d}}

\newcommand{\Eeta}{\boldsymbol \eta}

\newcommand{\Rplus}{\mathbb R_+}

\newtheorem{corollary}{Corollary}
\newtheorem{proposition}{Proposition}

\newcounter{assumptionHt}
\newcounter{assumptionH}
\newcounter{assumptionE}
\newcounter{assumptionLP}

{%
\begin{enumerate}[({G}1)]%
    \setcounter{enumi}{\value{assumptionHt}}%
    }{%
    \setcounter{assumptionHt}{\value{enumi}}%
\end{enumerate}
}
\newenvironment{assumptionH}%
{%
\medskip
\noindent\textbf{Assumptions.}
  \begin{enumerate}[({H}1)]%
  \setcounter{enumi}{\value{assumptionH}}%
}{%
  \setcounter{assumptionH}{\value{enumi}}%
  \end{enumerate}
}

{%
\medskip
\noindent \textbf{Assumptions.}
#1
\begin{enumerate}[({H}1)]%
  \setcounter{enumi}{\value{assumptionH}}%
  }{%
  \setcounter{assumptionH}{\value{enumi}}%
\end{enumerate}
}

{%
\medskip
\noindent \textbf{Assumptions.}
\begin{enumerate}[({E}1)]%
  \setcounter{enumi}{\value{assumptionE}}%
  }{%
  \setcounter{assumptionE}{\value{enumi}}%
\end{enumerate}
}

\newenvironment{assumptionLP}%
{%
\medskip
\noindent \textbf{Assumptions.}
\begin{enumerate}[({LP}1)]%
  \setcounter{enumi}{\value{assumptionLP}}%
  }{%
  \setcounter{assumptionLP}{\value{enumi}}%
\end{enumerate}
}

\newcommand{\assrefHt}[1]{(\hyperref[#1]{G\ref{#1}})}
\newcommand{\assrefH}[1]{(\hyperref[#1]{H\ref{#1}})}
\newcommand{\assrefE}[1]{(\hyperref[#1]{E\ref{#1}})}
\newcommand{\assrefLP}[1]{(\hyperref[#1]{LP\ref{#1}})}

\newcommand{\RR}{\mathbb{R}}

\newcommand{\Mmu}{\mathfrak{m}}

\title{Learning the regularity of multivariate functional data}
\author{Omar Kassi\footnote{Ensai, CREST - UMR 9194, France; omar.kassi@ensai.fr} \qquad	
	Nicolas Klutchnikoff\footnote{Univ Rennes, IRMAR - UMR 6625, F-35000 Rennes, France; Nicolas.klutchnikoff@univ-rennes2.fr}	\qquad
	Valentin Patilea\footnote{Ensai, CREST - UMR 9194, France; valentin.patilea@ensai.fr}
}
\date{\today}

\mathtoolsset{showonlyrefs}
\usepackage{xr}
\begin{document}

\maketitle

\begin{abstract}
Combining information both within and between sample realizations, we propose a simple estimator for the local regularity of surfaces in the functional data framework. The independently generated surfaces are measured with errors at possibly random discrete times. Non-asymptotic exponential bounds for the concentration of the regularity estimators are derived. An indicator for anisotropy is proposed and an exponential bound of its risk is derived. Two applications are proposed. We first consider the class of multi-fractional, bi-dimensional, Brownian sheets with domain deformation, and study the nonparametric estimation of the deformation. As a second application, we build minimax optimal, bivariate kernel estimators for the reconstruction of the surfaces.

	
\medskip	
	
	
	
\end{abstract}


\bigskip


\section{Introduction}\label{sec:introd}

Functional data analysis (FDA) provides methods for dealing with complex data such that collected by modern sensing devices. See, for instance, the textbooks \cite{ramsay_functional_2005}, \cite{horvath2012inference}, \cite{koko}. The paradigm consists of considering  that data are generated by a sample of functions, realizations of a stochastic process or random field defined over a continuous domain. However, the realizations are practically never observed over a continuous domain, and rarely without error. All the data points generated by such a realization then represent a single observation unit.  The remarkable advantage of functional data analysis is that it can combine information both within and between realizations. Restrictive assumptions, such as stationarity, stationary increments or Gaussianity on the data generating process or random field, can therefore be avoided.

We focus here on the case where the realizations are surfaces, \emph{i.e.}, the realizations are generated by a random scalar  field defined over a multi-dimensional continuous domain. We call this framework \emph{multivariate functional data}, and  focus on the case of a domain in the plane. 
Thus, in a different wat, we use existing FDA terminology, which  usually refers to a vector-valued  processes defined over an interval. In recent years, a wide panel of applications from different areas, including Astrophysics, Climate Sciences, Geophysics and Material Sciences, deal with data which can be considered as generated by random surfaces. 
For instance, it is now well admitted that the world ocean plays a key role in regulating Earth’s climate. Modern tools, such as floating sensors, provide  ocean heat transport measurements , which are made freely available by international programs, such as the Argo Program (http://www.argo.ucsd.edu, http://argo.jcommops.org). If one studies a specific area of the ocean, an observation unit is represented by the measurements collected at random points, sparsely distributed over the area, at some date in the year. See, for instance  \cite{kuusela_2018}, \cite{park2023}, and the references therein, for studies on Argo data. We aim at providing a new perspective for refined and effective analysis of multivariate functional data. 

Our main contribution is a new approach for studying the local regularity of random fields in the context of multivariate functional data. A main example of random field we have in mind is the  multifractional Brownian sheet, see \cite{herbin_06}. In the case of curves, that means for random fields defined over an interval, the local regularity can be defined naturally using the expectation of the squared increments to which one can impose a Hölder-like condition. The local regularity is then determined by the Hölder exponent and the Hölder constant. See \cite{GKP}. See also \cite{fracdim} where the local regularity exponent is linked to the fractal dimension for self-similar Gaussian processes. We here extend the ideas of \cite{GKP} to random fields defined over a domain in the plane. We thus introduce a general notion of  local regularity satisfied by a large class of random fields,  propose simple estimation procedures for the regularity parameters, and prove non-asymptotic results. \cite{fastandexact} consider a related estimation idea for a particular class of random fields, and use it to efficiently and exactly simulate fractional Brownian surfaces. \cite{hsing2020} study the estimation of the local regularity of a multifractional Brownian sheet from one realization of the process observed on a  regular grid, and provide asymptotic theory.

Knowing the regularity of the data generating random field has important consequences for FDA. For instance, there has been increasing interest in the nonparametric  estimation of the characteristics, such as the mean and the covariance structure, of the random field generating the functional data. See \cite{caponera2022} for a valuable review and an interesting approach. It is well-known that the optimal accuracy of nonparametric estimates depends on the regularity of the realizations. See, \cite{cai2010}, \cite{cai2011}, \cite{CMsphere}, \cite{GKP}. Inference methods for functional data should thus adapt to the regularity of the underlying process when aiming at optimality.  \cite{golovkine2023adaptive} and \cite{wei2023adaptive} used  regularity estimators to derive new, easy to implement, adaptive procedures for mean, covariance and functional principal components analysis for data generated by random curves. Adaptation to regularity for multivariate functional 
data seems yet unexplored.

The paper is organized as follows. In Section \ref{sec2}, the general observation scheme we consider is presented. It encompasses the  scenarios of \textit{common design} 
(the domain points where the random field realizations are  observed, possibly with noise,  are the same for all realizations)  and \textit{random design} (the observation domain points are randomly generated for each realization). Moreover, a general class of bivariate stochastic processes, for which the local regularity is well defined, is introduced. After discussion of some identification issues, in Section \ref{sec_*}, we present the estimation approach for the local regularity exponents, as well as the corresponding Hölder constants. Our estimators adapt to both isotropic and anisotropic settings. In Section \ref{sec4}, we provide  concentration bounds for the estimators of the local regularity, as well as a risk bound for the anisotropy detection. The new results are  of the non-asymptotic type, in the sense that they hold for any  number of random field realizations and observation domain points, provided these numbers are sufficiently large. Our estimation approach to local regularity for multivariate functional data opens the door to a large  array of adaptive procedures. Two applications are proposed. In Section \ref{BfMs}, we  consider the class of multi-fractional Brownian sheets with domain deformation, an example of a random field that belongs to the class defined in Section \ref{sec_*}. Deformed random fields have been studied in the literature, see for instance  \cite{Clerc2003},
\cite{anderes2008},
\cite{anderes2009consistent}, but yet seem little explored in the context of the functional data paradigm. 
As a second application, in Section \ref{sec6} we consider the problem of nonparametric reconstruction of the realizations of a random field from noisy measurements over a discrete set in the domain. A related problem was addressed by \cite{dun_aniso} in the context of Gaussian processes. With  the regularity estimates in hand, we build an adaptive Nadaraya-Watson pointwise estimates, and provide a sharp non-asymptotic bound for the average risk which achieves the optimal minimax rate expected in nonparametric statistics. The proofs of our results are given in the Appendix. Additional proofs and technical lemmas are provided in the Supplement.


\section{The framework}\label{sec2}
In this section we  present a formal mathematical setup for the local regularity  for bivariate   stochastic processes  (also called scalar random fields, or simply random fields) and the data observed for such processes.

\subsection{Data}\label{sec:data}
Consider $N$ independent realizations, also called sheets, $X^{(1)},\ldots,X^{(j)},\ldots X^{(N)}$   of a stochastic process $X $ defined over a continuous domain $\cT\in\mathbb R^2$. For simplicity, we focus here on domains $\mathcal T$ in the plane, though the extension to higher dimensions would not raise different challenges. For the purpose of describing our methodology, we distinguish three observational scenarios. 
First, the ideal, infeasible situation where the sheets $X^{(j)}$ are \emph{completely observed}, i.e. without error over the entire domain $\cT$.  
The second case is the one where the $X^{(j)}$ are observed (measured) at some \emph{discrete points} in the domain $\cT$, \emph{without noise}. The observation points can be fixed to be the same for all the $X^{(i)}$'s (common design), or can be randomly drawn for each sheets separately (independent design). Finally, the most realistic scenario is the one where in the second case, we admit that the realizations of $X$  are observed at discrete domain points \emph{with noise}.

 To formally describe the second and third scenarios, let  $M_1, \dotsc,M_N$ be an independent sample of an integer-valued random variable $M$ with expectation $\EE[M]=\Mmu$. 
In the independent design case, for each $1\leq j \leq N$, and given  $M_j$, let $\Tnm\in\cT$,  $1\leq m \leq  M_j$, be a random sample of a random vector $\TT\in\cT$. The $\Tnm$'s represent the observation points for the realization $\Xp{j}$. We assume that the realizations of $X$, $M$ and $\TT$ are mutually independent. 
In the common design case, $M\equiv \mathfrak m$ and the $\Tnm$'s are the same for all $j$.  Let $\mathcal T_{obs}^{(j)} $ denote the set of observation points $ \Tnm $, $1\leq m \leq M_j$, on the sheet $\Xp{j}$. With common design,  $\mathcal T_{obs}^{(j)} $ does not depend on $j$, while with independent design, the expected cardinal of $\mathcal T_{obs}^{(j)} $ can be random with mean $\mathfrak m$.   The following presentation  includes both independent design and common design cases. \color{black}  Finally, the data  consist of  the pairs  $(\Ynm , \Tnm ) \in\mathbb R \times \cT $ where $\Ynm$ is defined as
\begin{equation}\label{model_eq}
	\Ynm = \Xn (\Tnm) + \enm, \quad\text{with}  \quad  \enm = \sigma(\Tnm,\Xn(\Tnm)) \unm, 
	\quad 1\leq i \leq N,  \; 1\leq m \leq M_j.
\end{equation}
Here, the $\unm  \in\mathbb R $ are independent copies of a centered variable $e$ with unit variance, and $\sigma^2(\cdot,\cdot)\geq 0$ is some unknown, bounded conditional variance  function which account for possibly heteroscedastic measurement errors. The case  $\sigma(t,x)\equiv 0$ corresponds to our second scenario, while in the third scenario we have positive conditional variance. 

For each $1\leq j \leq N$,  let $\widetilde X^{(j)}$ denote an observable approximation of $X^{(j)}$. If   the sheets $X^{(j)}$ were  completely observed, as in our infeasible first scenario, $\widetilde X^{(j)} = X^{(j)}$. 
When $X^{(j)}$ are observed  only at  some discrete points $\Tnm$,  arbitrary  $\widetilde X^{(j)}(\Tt)$ can be obtained by simple interpolation or defined equal to the value of $\widetilde X^{(j)}$ at the nearest neighbor of $\Tt$.  Finally, with noisy, discretely observed sheets,   $\widetilde X^{(j)}$  is a pilot nonparametric estimator  of $X^{(j)}$, such as kernel smoothing, splines \emph{etc}.

Let us next introduce a general class of stochastic processes (random fields)  $X$ with  irregular realizations $X^{(j)}$, for which the regularity can vary over the domain $\mathcal T$.

\subsection{A class of multivariate processes}
Let $\mathcal{T}$ be an open, bounded rectangle with the closure included in $(0,\infty)^2$. In the following, $H_1,H_2 : \mathcal T \to (0,1)$ are two continuously differentiable functions such that 
\begin{equation}\label{low_thres}
\underline \beta=\min_{i=1,2} \inf_{t\in\mathcal T} H_i(\Tt) >0.
\end{equation} 	
Let $\mathbf{L}=(L_1^{(1)},L_2^{(1)},L_1^{(2)},L_2^{(2)}),$ be a  vector-valued function with  non-negative, Lipschitz  continuous, components  functions defined on $\mathcal{T}$ such that 
\begin{equation}\label{id_L}
	L_i^{(1)}(\Tt) +L_i^{(2)}(\Tt) >0,\qquad \forall \Tt\in\mathcal T, \; i=1,2.
\end{equation}

Let $X$ be a real-valued, second order stochastic process defined over $(0,\infty)^2$. Let $(e_1,e_2)$ be the canonical basis of $\mathbb R ^2,$ and, for sufficiently small scalars $\Delta$, let   
$$
\theta_{\Tt}^{(i)}(\Delta)=\EE\left[\left\{X\left(\Tt-\frac{\Delta}{2}e_i\right)-X\left(\Tt+\frac{\Delta}{2}e_i\right)\right\}^2\right],\quad i=1,2.
$$

\begin{definition}\label{def}
Let $H_1$, $H_2$ satisfy \eqref{low_thres}.	The class $\mathcal {H}^{H_1,H_2}(\mathbf{L},\mathcal{T})$ is the set of stochastic processes $X$ satisfying the following condition:  constants $  \Delta_0, C,\beta>0$ exist such that 
	for any $\Tt\in \mathcal T$ and $ 0<\Delta\leq\Delta_0$,  
	\begin{equation}\label{as_repr}
		\left|\theta_{\Tt}^{(i)}(\Delta)-L_1^{(i)}(\Tt)\Delta^{2H_1(\Tt)} -L_2^{(i)}(\Tt)\Delta^{2H_2(\Tt)}\right|\leq C\Delta^{2 \max\{H_1(\Tt),H_2(\Tt)\}+\beta}, \quad i=1,2.
	\end{equation}
	Let  
	$$
	\mathcal{H}^{H_1,H_2} 
	=\bigcup_{\mathbf{L}}\mathcal {H}^{H_1,H_2}(\mathbf{L},\mathcal{T}) ,
	$$
	where the union is taken over the set of four-dimensional functions $\mathbf{L}$ with non negative positive, Lipschitz  continuous components satisfying \eqref{id_L} 
	The functions $H_1,H_2$ define the local regularity of the process, while $\mathbf{L}$ represent  the local Hölder constants.
\end{definition}

Definition \ref{def}  is general, and extends the local regularity notion considered by \cite{GKP} to processes defined over a compact interval on the real line.  A main example we have in mind is the  multi-fractional Brownian  sheet (MfBs) with a domain deformation. MfBs  is a generalization of the standard fractional Brownian sheet, where the Hurst parameter is allowed to vary along the  domain. The definition of this general class of  processes and some of their properties are provided in Section \ref{BfMs}.


\section{Local regularity estimation approach}\label{sec_*}
The idea is to relate the functional parameters $H_1,H_2$ and $\mathbf L$ to quantities which are estimable from the data. In other words, we build estimating equations for each of the parameters we want to estimate. The parameter estimate is then obtained from the sample version of the estimating equation.  

Before proceeding with the local regularity estimation, let us discuss  the identification aspect. Definition \ref{def} is too general and does not allow  all the unknown parameters to be identifiedwithout further restrictions. To be more clear, let $H_1,H_2,\widetilde H_1$ and $\widetilde H_2$  be  continuously differentiable functions  taking values in $(0,1).$ Assume that  $X\in \mathcal {H}^{H_1,H_2}(\mathbf{L},\mathcal{T})$ and  $X\in \mathcal {H}^{\widetilde H_1,\widetilde H_2}(\widetilde{\mathbf L},\mathcal{T})$, for some $\mathbf{L}$ and $\widetilde{\mathbf L}$. Then necessarily  
$$
\min\{H_1(\Tt),H_2(\Tt)\} \! = \! \min\{\widetilde H_1(\Tt),\widetilde H_2(\Tt)\}\hspace{0.1cm}\text{ and } \hspace{0.1cm}\max\{H_1(\Tt),H_2(\Tt)\} \! = \! \max\{\widetilde H_1(\Tt),\widetilde H_2(\Tt)\},
$$ 
 for any $\Tt\in \cT$, and, modulo a permutation of the components,   $\mathbf{L}\equiv \widetilde{\mathbf L}$. In general, the permutation depends on the domain point $\Tt$. We deduce from these facts that, for instance, only  
$$\underline H (\Tt)=\min\{H_1(\Tt),H_2(\Tt)\}\quad \text{  and  } \quad \overline{H}(\Tt)=\max\{H_1(\Tt),H_2(\Tt)\},
$$ 
are expected to be identifiable in the general framework we consider. Concerning the components of  $\mathbf{L}$,  the identifiable quantities are provided below.   

\subsection{ Estimating equations for $\underline{H}(\Tt)$ and $\overline{H}(\Tt)$}
Let  $X\in \mathcal H^{H_1,H_2}$. 
Since $\Delta^{b}$ is negligible compared to $\Delta^{a}$ if $0< a <b$ and $\Delta$ is small, in view of Definition \ref{def} we first define  the estimation equation for $\underline{H}(\Tt)$, for some fixed $\Tt\in \mathcal T$. 

For $i=1,2$ and $\Delta$ sufficiently small, we have
$$\theta_{\Tt}^{(i)}(\Delta)= K_1^{(i)}(\Tt)\Delta^{2\underline H (\Tt)}+O(\Delta^{\widetilde{\beta}})=K_1^{(i)}(\Tt)\Delta^{2\underline H (\Tt)}+K_2^{(i)}(\Tt)\Delta^{2\overline H (\Tt)}+O(\Delta^{\overline H (\Tt)+\beta}),$$
where 
$$K_1^{(i)}(\Tt) =\left\lbrace\begin{array}{lll}
	\!\!  L_1^{(i)} (\Tt)\quad &\!\!\! \text{if }  H_1(\Tt)< H_2(\Tt)\\
	\!\! L_2^{(i)}(\Tt)\quad &\!\!\!\text{if }   H_2(\Tt)< H_1(\Tt)\\
	\!\! L_1^{(i)}(\Tt)+L_2^{(i)}(\Tt) \quad &\!\!\!\text{if }   H_1(\Tt)=H_2(\Tt) 
\end{array}\right. \!\!, \quad K_2^{(i)}(\Tt) =\left\lbrace\begin{array}{lll}
	\!\!  L_1^{(i)} (\Tt)\quad &\!\!\! \text{if }  H_1(\Tt)> H_2(\Tt)\\
	\!\! L_2^{(i)}(\Tt)\quad &\!\!\!\text{if }   H_2(\Tt)>H_1(\Tt)\\
	\!\! 0\quad &\!\!\!\text{if }   H_1(\Tt)=H_2(\Tt) 
\end{array}\right. \!\!,
$$
and 
$$
\widetilde \beta =\left\lbrace\begin{array}{lll}
	\!\! 2\overline{H}(\Tt)& \text{if}&  \underline H(\Tt)< \overline{H}(\Tt) \\
	\!\! 2\underline H(\Tt)+\beta & \text{if}&  \underline H(\Tt)=\overline{H}(\Tt) 
\end{array}\right. ,\;\; \text{  } \;\; 
$$
Related to the previous discussion on the identifiability, 
similarly to the role of $\underline H (\Tt)$ and $\overline H (\Tt)$ for $H_1(\Tt)$ and $H_2(\Tt)$, the functions $K_1^{(i)}(\Tt)$ and $K_2^{(i)}(\Tt)$, $i=1,2$
are the identifiable functionals of $\mathbf L$. More precisely, given  the order choice in the case $H_1(\Tt)\neq H_2(\Tt)$, the functions $K_1^{(i)}(\Tt)$ and $K_2^{(i)}(\Tt)$ represent  the identifiable components of $\mathbf L$. When $H_1(\Tt)=H_2(\Tt) $, only the $ L_1^{(i)}(\Tt)+L_2^{(i)}(\Tt) $ are identifiable. See also the discussion following Proposition \ref{conc_Lest}.

Next, we define $$\gamma_{\Tt}(\Delta)=\theta^{(1)}_{\Tt}(\Delta)+\theta^{(2)}_{\Tt}(\Delta).$$
The reason for considering this quantity, instead of considering $\theta^{(1)}_{\Tt}(\Delta)$ and $\theta^{(2)}_{\Tt}(\Delta)$ separately,
is that the Hölder constant associated with $\underline{H} (\Tt)$ can vanish, and this would prevent us from estimating the lower regularity exponent. On contrary,  $\gamma_{\Tt}(\Delta)$ can be written as 
\begin{multline}\label{eq:K1K2}
\gamma_{\Tt}(\Delta)=\left(K_1^{(1)}(\Tt)+K_1^{(2)}(\Tt)\right)\Delta^{2\underline H (\Tt)} +
\left(K_2^{(1)}(\Tt)+K_2^{(2)}(\Tt)\right)
\Delta^{2\overline H (\Tt)}+ O(\Delta^{ 2\overline H (\Tt) +\beta})\\=: K_1(\Tt)\Delta^{2\underline H (\Tt)}
+ K_2(\Tt)\Delta^{2\overline H (\Tt)}+O(\Delta^{2\overline H (\Tt) +\beta}),
\end{multline}
and  condition \eqref{id_L} guarantees  $K_1(\Tt), K_2(\Tt)>0$, and thus allows us to estimate $\underline H (\Tt)$ consistently. We also consider 
\begin{equation}\label{eq:alpha}
	\alpha_{\Tt}(\Delta)=\left|\frac{\gamma_{\Tt}(2\Delta)}{(2\Delta)^{2\underline{H}(\Tt)}}-\frac{\gamma_{\Tt}(\Delta)}{\Delta^{2\underline{H}(\Tt)}}\right|.
\end{equation}

\quad 

\begin{proposition}\label{proprox}
Let $X$ belong to the class $ \mathcal {H}^{H_1,H_2}(\mathbf{L},\mathcal{T})$, introduced by Definition \ref{def}. Then, for any $\Tt\in\cT$, 
	\begin{equation}\label{proxy_H_low}
		\underline{H}(\Tt) = \frac{\log(\gamma_{\Tt}(2\Delta))-\log(\gamma_{\Tt}(\Delta))}{2\log(2)} + O(\Delta^{\widetilde \beta -2\underline H(\Tt)}),
	\end{equation}
	and
	\begin{equation*}
		\overline{H}(\Tt)-\underline{H}(\Tt) = \frac{\log(\alpha_{\Tt}(2\Delta))-\log(\alpha_{\Tt}(\Delta))}{2\log(2)} + O(\Delta^{ \beta }).
	\end{equation*}
\end{proposition}

\smallskip

To estimate $\underline H (\Tt)$, we thus use the dominating term on the right-hand side of the  representation \eqref{proxy_H_low} as a proxy, for which we compute an estimate. 
To build a sample counterpart of the proxy quantity, we can estimate 
$\theta_{\Tt}^{(i)}(\Delta)$ by 
\begin{equation} \label{eq_theta_hat}
\widehat{\theta}_{\Tt}^{(i)}(\Delta)= \frac{1}{N}\sum_{j=1}^{N}\left\{\widetilde{X}^{(j)}(\Tt-(\Delta/2) e_i)-\widetilde{X}^{(j)}(\Tt+(\Delta/2) e_i)\right\}^2, \quad i=1,2,
\end{equation}
where $\widetilde X^{(j)}$ is the observable approximation of $\Xn$, In the ideal, infeasible scenario where    the sheets $X^{(j)}$ are completely observed, 
$\widetilde X^{(j)} = X^{(j)}$. When $X^{(j)}$ are observed  only at  some discrete points $\Tnm$,  the  $\widetilde X^{(j)}(\Tt)$'s can be obtained by simple interpolation   or using nearest neighbors. 
Finally, with noisy, discretely observed sheets,   $\widetilde X^{(j)}$  can be a pilot nonparametric estimator  of $X^{(j)}$, such as bivariate kernel smoothing, splines \emph{etc}.

An estimator of $\gamma_{\Tt}(\Delta)$ is then given by $\widehat{\gamma}_{\Tt}(\Delta)=\widehat{\theta}_{\Tt}^{(1)}(\Delta)+\widehat{\theta}_{\Tt}^{(2)}(\Delta).$ By plugging this estimator of $\gamma_{\Tt}(\Delta)$ into \eqref{proxy_H_low}, we obtain an estimator of $\underline{H}(\Tt)$~:
\begin{equation}\label{est_under}
	\widehat{\underline{H}}(\Tt)=\left \lbrace 
	\begin{array}{cl}
		\frac{\log(\widehat{\gamma}_{\Tt}(2\Delta))-\log(\widehat{\gamma}_{\Tt}(\Delta))}{2\log(2)}\quad &\text{if}\quad \widehat{\gamma}_{\Tt}(2\Delta),\widehat{\gamma}_{\Tt}(\Delta)>0\\
		1&\text{otherwise}
	\end{array}
	\right..
\end{equation}
Moreover, replacing $\gamma_{\Tt}$ by  $\widehat{\gamma}_{\Tt}$  and $\underline{H}(\Tt)$ by  $\widehat{\underline{H}}(\Tt)$ in \eqref{eq:alpha}, we get an estimator of $\alpha_{\Tt}$~: 
\begin{equation}\label{eq:alpha_hat}
	\widehat{\alpha}_{\Tt}(\Delta)=\left \lbrace 
	\begin{array}{cl}
		\left|\frac{\widehat{\gamma}_{\Tt}(2\Delta)}{(2\Delta)^{2\widehat{\underline{H}}(\Tt)}}-\frac{\widehat{\gamma}_{\Tt}(\Delta)}{\Delta^{2\widehat{\underline{H}}(\Tt)}}\right|\quad &\text{if}\quad \frac{\widehat{\gamma}_{\Tt}(2\Delta)}{(2\Delta)^{2\widehat{\underline{H}}(\Tt)}}\ne\frac{\widehat{\gamma}_{\Tt}(\Delta)}{\Delta^{2\widehat{\underline{H}}(\Tt)}}\\
		1&\text{otherwise}.
	\end{array}
	\right..
\end{equation}
Finally,  using the second part of Proposition \ref{proprox}, we estimator an estimator of $\overline{H}(\Tt)-\underline{H}(\Tt)$ under the form 
$$
\widehat{(\overline{H}-\underline{H})}(\Tt)=\frac{\log(\widehat{\alpha}_{\Tt}(2\Delta))-\log(\widehat{\alpha}_{\Tt}(\Delta))}{2\log(2)}.
$$


It will be shown below that, for the pointwise estimation of $\overline{H}$ and $\mathbf L$,  we have to distinguish between the isotropic and anisotropic cases. Here, the \emph{isotropic} and \emph{anisotropic} cases are defined locally, and correspond to  $\underline{H}(\Tt)=\overline{H}(\Tt)$ and $\underline{H}(\Tt)< \overline{H}(\Tt)$, respectively. We propose here an estimator of $\overline{H}(\Tt)$ which  adapts to isotropy. Let us consider  the  event 
\begin{equation}\label{def_A_N}
	 A_N(\tau)=A_N(\tau; \Tt)=  \left\{ \widehat{(\overline{H}-\underline{H})}(\Tt)\geq \tau\right\},
\end{equation}
for some appropriate, small threshold $\tau>0$. 
We then define  the following estimator for $\overline{H}(\Tt)$~:
\begin{equation}\label{est_over}
	\widehat{\overline{H}}(\Tt)=\widehat{\underline{H}}(\Tt)+\widehat{(\overline{H}-\underline{H})}(\Tt)\mathbf{1}_{A_N(\tau)}.
\end{equation}
 Here, for a set $S$, $\mathbf 1_S$ denotes the indicator of $S$.   In Section \ref{info_tau}, we provide an exponential bound  for the probability that the anisotropy detection rule defined by $\mathbf{1}_{A_N(\tau)}$ fails. In particular, that indicates how small $\tau $ is allowed to be such that $\mathbf{1}_{A_N(\tau)}$ detects anisotropy with high probability.

\subsection{Estimating equations for $\mathbf{L}(\Tt)$}

Assume for the moment that $\underline{H}(\Tt)< \overline{H}(\Tt)$. A sample-based diagnosis   procedure for detecting this situation  can be built using the results in Section \ref{info_tau} below. Without 
loss of generality, we consider  $\underline{H}(\Tt)=H_1(\Tt)$.
Let us recall that, for $i=1,2$,
$$
\theta_{\Tt}^{(i)}(\Delta)=L_1^{(i)}(\Tt)\Delta^{2H_1(\Tt)} +L_2^{(i)}(\Tt)\Delta^{2H_2(\Tt)}+O(\Delta^{2{H}_2+\beta}).
$$

\begin{proposition}\label{prop_Lest} Let $ X\in \mathcal {H}^{H_1,H_2}$.
	Denote $D(\Tt)= {H}_2(\Tt)-{H}_1(\Tt) >0$.
	For $i=1,2$,
	\begin{equation*}
		L_1^{(i)}(\Tt)= \frac{\theta_{\Tt}^{(i)}(\Delta)}{\Delta^{2 H_1 (\Tt)}}+O(\Delta^{2{D}(\Tt)}),
	\end{equation*}
and
	\begin{equation*}
		L_2^{(i)}(\Tt)= \frac{1}{(2^{2D(\Tt)}-1)\Delta^{2D(\Tt)}}\left|\frac{\theta_{\Tt}^{(i)}(2\Delta)}{(2\Delta)^{ 2H_1(\Tt)}}-\frac{\theta_{\Tt}^{(i)}(\Delta)}{\Delta^{ 2H_1(\Tt)}}\right|+O(\Delta^\beta).
	\end{equation*}
\end{proposition}
We denote the estimators of the local Hölder constants by
\begin{equation}\label{est_Lcomp}
\widehat{L_1^{(i)}}(\Tt), \quad \widehat{L_2^{(i)}}(\Tt), \qquad i=1,2.
\end{equation}
The estimators of $L_1^{(i)}(\Tt)$, $i=1,2$, are obtained by plugging into its dominating term derived in Proposition \ref{prop_Lest}, the estimators  in \eqref{eq_theta_hat}, \eqref{est_under}. For the estimators of $L_2^{(i)}(\Tt)$, we first consider  
$$
\widehat D(\Tt) =  \widehat{(\overline{H}-\underline{H})}(\Tt) \;\;\; \text{ if } \; \widehat{(\overline{H}-\underline{H})}(\Tt)\neq 0,\quad \text{ and } \; \widehat D(\Tt) = 0 \;\text{ otherwise}.
$$
If $\widehat D(\Tt) \neq  0$,  the estimators of $L_2^{(i)}(\Tt)$, $i=1,2$, are obtained by plugging into its dominating term, the estimated quantities, otherwise they are set equal to zero.


\section{Non-asymptotic results }\label{sec4}
We now derive concentration inequalities for the pointwise estimators of the  parameters $(H_1,H_2)$ and $\mathbf{L}= (L_1^{(1)},L_2^{(1)},L_1^{(2)},L_2^{(2)})$. For this purpose, we need to measure the error between each realizations $X^{(j)}$ of $X$ and its observable approximation of $\Xtp{j}$, as  defined in Section \ref{sec2}. We consider the following $\mathbb{L}^p$-risk~: 
$$
R_p (\mathfrak m) =\sup_{\Tt\in\cT }\EE[|\xi^{(j)}(\Tt)|^p] ,\qquad \xi^{(j)}(\Tt)=\Xtp{j}(\Tt)-\Xp{j}(\Tt).
$$
In general, the $\mathbb{L}^p$-risk depends on the expected number $\mathfrak m$ of observed points $\Tnm$. In the ideal scenario where $\Xp{j}$ is observed everywhere without error, $R_p \equiv 0$.
We also consider the following assumptions. Below, $B(\Tt; r)$ denote the ball of radius $r$ centered at $\Tt$.

\begin{assumptionH}
	\item\label{ass_D}  Let $X$ belong to the class $ \mathcal {H}^{H_1,H_2}$,
	introduced by Definition \ref{def}, and let $X^{(j)}$, $1\leq j \leq N$,  be independent realizations of $X$.
	
	\item\label{ass_H1}  Three positive constants $\mathfrak{a}$, $\mathfrak{A}$ and $r$ exist such that, for any $\Tt\in\mathcal T$,  
	$$
	\EE\left| 	X^{(j)}\left(\Tt\right)-
	 X^{(j)}\left(\Ss
	 \right)\right|^{2p}\leq  \frac{p!}{2}\mathfrak{a} \mathfrak{A}^{p-2} \|\Tt-\Ss\|^{2p\underline H (\Tt)}
	\qquad \forall \Ss\in B(\Tt; r) ,\; \forall p\geq 1.
	$$
	
	\item\label{ass_H2} Two positive constants $\mathfrak{c}$ and $\mathfrak{D}$, and a function $\rho(\mathfrak m)\leq 1$, exist such that 
	$$
	R_{2p} (\mathfrak m) \leq  \frac{p!}{2}\mathfrak{c} \mathfrak{D}^{p-2}\rho(\mathfrak m)^{2p}, \qquad \forall p\geq 1,\; \forall \mathfrak m>1.
	$$
	\item\label{ass_H3} Two positive constants  $\mathfrak L $ and $\nu$ exist such that 
	$$
	R_2(\mathfrak m) \leq \mathfrak L \mathfrak m ^{-\nu},\qquad \forall \mathfrak m >1.
	$$
\end{assumptionH}

\smallskip

The condition in (H\ref{ass_H1}) imposes  sub-Gaussian  local increments for $X$.  It is satisfied by the processes in the wide class of
multi-fractional Brownian  sheets (MfBs) with a  domain-deformation, as considered  in Section \ref{BfMs}. In the case of noisy, discretely observed realizations $X^{(j)}$, the observable approximation can be obtained from existing bivariate nonparametric smoothing approaches. Under mild conditions, the standard nonparametric smoothers satisfy Assumption \ref{ass_H2}, with $\rho(\mathfrak m) =1$, and Assumption \ref{ass_H3}. See \cite{fan2016multivariate} for the case of local polynomials, and \cite{BELLONI2015} for general series estimators. In the second scenario, where the $X^{(j)}$ are observed  without noise  at discrete points $\Tnm$ in the domain $\mathcal T$, we can simply define $ \widetilde{X}^{(j)}(\Tt)$ as the value of $X^{(j)}$ at the nearest observed point $\Tnm$  to $\Tt$. To provide a simple justification that this simple choice is valid, let us consider that 
a constant $C>0$ exists such that 
$$
C^{-1} \leq  M_j/\mathfrak m\leq C,\quad \forall 1\leq j\leq N.
$$	
Then, with probability exponentially close to 1, there exists at least one point $\Tnm$ in the ball  $B(\Tt; \widetilde r)$, provided $\widetilde r = \mathfrak m^{-\delta}$, for some $\delta \in(1/2,1)$. 
Assumption \ref{ass_H2} is then implied by \ref{ass_H1} with    $\rho(\mathfrak m)=\mathfrak m^{-\delta\underline \beta}$, and $\underline \beta $ from \eqref{low_thres}.   In particular, this also guarantees \ref{ass_H3} with   $\nu = 2\delta\underline \beta$.

\subsection{Concentration bounds for the regularity estimates}

We first derive  the exponential bound for the concentration of the local regularity exponents.  On the one hand, the concentration will depend on the non-stochastic approximation error between
the true parameter and their respective dominating terms. From Proposition  \ref{proprox} these approximation errors are
\begin{equation*}
	R(\underline H )(\Tt) = 	\underline{H}(\Tt) - \frac{\log(\gamma_{\Tt}(2\Delta))-\log(\gamma_{\Tt}(\Delta))}{2\log(2)} ,
\end{equation*}
and
\begin{equation*}
	R(\overline H - \underline H )(\Tt) = \{\overline H -\underline H \}(\Tt) - \frac{\log(\alpha_{\Tt}(2\Delta))-\log(\alpha_{\Tt}(\Delta))}{2\log(2)},
\end{equation*}
respectively. We have
\begin{equation}\label{rates_eRR}
	R(\underline H )(\Tt)  = O(\Delta^{\widetilde \beta -2\underline H(\Tt)}) \quad \text{and} \quad 
	R(\overline H - \underline H )(\Tt) = 	 O(\Delta^{ \beta }).
\end{equation}
On the other hand, the concentration of the  local regularity exponents estimators will also depend on the error between the realizations of $X$ and their observable approximations $\widetilde X^{(j)}$. Finally, since we use Bernstein's inequality, the concentration will also depend on the bound of the moments  in Assumption \ref{ass_H1}. To account for these, let 
$$
\varrho(\Delta,\mathfrak m) = \max\{ \Delta^{2\underline H (\Tt)},\rho^{2}(\mathfrak m)\}^{-1}.
$$ 
Note that $\varrho(\Delta,\mathfrak m) = \Delta^{-2\underline H (\Tt)}$ in the ideal case where $\widetilde X^{(j)} = X^{(j)}$ and thus $\rho (\mathfrak m)=0$. 

\begin{proposition}\label{propCH}
Assumptions (H\ref{ass_D}) to  (H\ref{ass_H3}) hold true. Let $\widehat{\underline{H}}(\Tt) $  and $ \widehat{\overline{H}}(\Tt)$ be the estimators defined in \eqref{est_under} and \eqref{est_over}, respectively. If $\Delta$ is sufficiently small and $\mathfrak m$ sufficiently large,  constants $C_1,\dots,C_5$ exist such that,  
\begin{equation}\label{eq:cdt_eps}
\forall \varepsilon,\tau  \in (0,1) \quad \text{satisfying} \quad \max \{ |\log(\Delta)| |R(\underline H )(\Tt)|, \;|R(\overline H - \underline H )(\Tt)  |\}\leq \varepsilon \leq 2\tau,
\end{equation}
then
\begin{equation}\label{eq:conc-Hhat-around-H_main}
	\mathbb{P}\left[
	|\underline{\widehat{H}}(\Tt)-\underline{H}(\Tt)|\geq 
	\varepsilon 
	\right]
	\leq  p_1,
\end{equation}
and 
\begin{equation}\label{eq:concentration-overlineH_main}
	\PP\left[\left|\widehat{\overline{H}}(\Tt)-\overline{H}(\Tt)\right|\geq \varepsilon\right] 
	\leq C_3\{p_1+p_2+p_3\},
\end{equation}
with
\begin{align}
	p_1&= C_1\exp \left(-C_2N \times \varepsilon^2 \times \Delta ^{4\underline{H}(\Tt)}\varrho(\Delta,\mathfrak m)\right),\\
	p_2 &=	
	\exp\left[ - C_4N\times \varepsilon^2\times \frac{\Delta^{4\overline{H}(\Tt)}\varrho(\Delta,\mathfrak m)}{\log^2(\Delta)}\Delta^{4D(\Tt)}
	\right]\mathbf1_{\{\underline H(\Tt)<\overline H (\Tt)\}},
	\\ p_3&=  \exp\left[
	- C_5N\times \tau^2 \times \frac{\Delta^{4\overline{H}(\Tt)}\varrho(\Delta,\mathfrak m)}{\log^2(\Delta)}\Delta^{4D(\Tt)}
	\right],
\end{align}
where 
$$
D(\Tt)= \overline{H}(\Tt)-\underline{H}(\Tt) \quad \text{ and }  \quad 
\varrho(\Delta,\mathfrak m) = \max\{ \Delta^{2\underline H (\Tt)},\rho(\mathfrak m)^{2}\}^{-1}.
$$ 
\end{proposition}

\medskip 

The term $p_2$ is specific to the anisotropic case, it disappears when $\underline H(\Tt)=\overline H (\Tt)$. We next derive  the bounds for the concentration of the local Hölder constants' estimators. 
In the case where $\underline{H}(\Tt) \neq \overline{H}(\Tt)$, without loss of generality, we set
$$
\underline{H}(\Tt)=H_1(\Tt) < H_2(\Tt)= \overline{H}(\Tt),
$$
such that $L_1^{(1)}(\Tt)$ and $L_1^{(2)}(\Tt)$ are the Hölder constants corresponding to $\underline{H}(\Tt)$.
et 
\begin{equation*}
	R(L_1^{(i)})(\Tt) = 	L_1^{(i)}(\Tt) -\frac{\theta_{\Tt}^{(i)}(\Delta)}{\Delta^{2 H_1 (\Tt)}}=O(\Delta^{2{D}(\Tt)}),
\end{equation*}
and
\begin{equation*}
	R(L_2^{(i)})(\Tt) = 	L_2^{(i)}(\Tt) - \frac{1}{(2^{2D(\Tt)}-1)\Delta^{2D(\Tt)}}\left|\frac{\theta_{\Tt}^{(i)}(2\Delta)}{(2\Delta)^{ 2H_1(\Tt)}}-\frac{\theta_{\Tt}^{(i)}(\Delta)}{\Delta^{ 2H_1(\Tt)}}\right| = O(\Delta^\beta), \quad i=1,2.
\end{equation*}

\medskip

\begin{proposition}\label{conc_Lest}
Assume that the conditions of Proposition \ref{propCH} hold true. Then, for the estimators in \eqref{est_Lcomp}, positive constants $\mathfrak C_1,...,\mathfrak C_4$ exists such that,  for $i=1,2$, and \color{black}	
\begin{equation}\label{eq:cdt_epsL}
	\forall \varepsilon   \in (0,1)  \text{ satisfying }  \max \left\{|R(L_1^{(i)})(\Tt)|,\; |\log(\Delta)| |R(\underline H )(\Tt)|, \; |R(\overline H \!- \!\underline H )(\Tt)|  \right\}\leq \varepsilon ,
\end{equation}
 we have 

\begin{equation}\label{eq:conc_L1_main}
\PP\left(\left|\widehat{L_1^{(i)}}(\Tt)-L_1^{(i)}(\Tt)\right| \geq \varepsilon \right) \leq 
\mathfrak C_1 \exp\left(
	- \mathfrak C_2 N \times \varepsilon^2\times \frac{\Delta^{4\underline H(\Tt)}\varrho(\Delta,\mathfrak m)}{\log^2(\Delta)}
	\right).
\end{equation} 
 Moreover,  if in addition $ |R(L_2^{(i)})(\Tt)|\leq \varepsilon$, $i=1,2$ then
\begin{multline}\label{eq:conc_L2_main}
	\PP\left(\left|\widehat{L_2^{(i)}}(\Tt)-L_2^{(i)}(\Tt)\right|
	\geq \varepsilon \right) \\ \leq\mathfrak  C_3\exp\left(
	-\mathfrak C_4N\times \varepsilon\Delta^{4D(\Tt)}\min\{\varepsilon,\Delta^{4D(\Tt)}\}
	\times \frac{\Delta^{4\overline H (\Tt)}\varrho(\Delta,\mathfrak m)}{\log^4(\Delta)}\times 
(2^{2D(\Tt)}-1)^2	\right).
\end{multline} 
\end{proposition}

\medskip

The second exponential bound in Proposition \ref{conc_Lest} becomes trivial when 
$D(\Tt)=0$, and this reveals that the case $H_1(\Tt)=H_2(\Tt)$  requires special attention. 
In  this case,  the estimator proposed for  $L_1^{(i)}(\Tt)$  becomes an estimator of $L_1^{(i)}(\Tt)+L_2^{(i)}(\Tt)$, $i=1,2$. 
 The indicator of the set defined in \eqref{def_A_N} provides a tool for detecting whether or not $H_1(\Tt)=H_2(\Tt)$, given a data set. In the following, we investigate  the risk associated with this diagnosis tool.

\subsection{A risk bound for the anisotropy detection}\label{info_tau}
Assume without loss of generality that $H_1(\Tt)\leq H_2(\Tt)$. 
Equation \eqref{eq:K1K2} then becomes 
\begin{multline}
\gamma_{\Tt}(\Delta)=\theta_{\Tt}^{(1)}(\Delta)+\theta_{\Tt}^{(2)}(\Delta)\\
=\left(L_1^{(1)}(\Tt)+L_1^{(2)}(\Tt)\right)\Delta^{2{H}_1(\Tt)} +\left(L_2^{(1)}(\Tt)+L_2^{(2)}(\Tt)\right)\Delta^{2{H}_2(\Tt)}+O(\Delta^{2{H}_2(\Tt)+\beta})\\ 
=K_1(\Tt)\Delta^{2H_1(\Tt)}+K_2(\Tt)\Delta^{2H_2(\Tt)} +O(\Delta^{2H_2(\Tt) +\beta}).
\end{multline}
We can now write
$$
\frac{\log(\alpha_{\Tt}(2\Delta))-\log(\alpha_{\Tt}(\Delta))}{2\log2}= D(\Tt)+O\left(\Delta^\beta\right).
$$
Therefore, if $D(\Tt)=H_2(\Tt)- H_1(\Tt)=0$, we get 
$$
\frac{\log(\alpha_{\Tt}(2\Delta))-\log(\alpha_{\Tt}(\Delta))}{2\log2}= O\left(\Delta^\beta\right).
$$ 
We deduce that,  for the event $A_N(\tau)$
introduced in \eqref{def_A_N},  we have to  choose $\tau$ such that $\Delta = o(\tau^{1/\beta})$. The following result proposes an exponential bound for the risk associated to the rule defined by the indicator $\mathbf{1}_{A_N(\tau)}$ in the definition \eqref{est_over}.

\medskip

\begin{proposition}\label{prop5_simple} Assume that the conditions of Proposition \ref{propCH} hold true. Let 
	$$  \max \{ |\log(\Delta)| |R(\underline H )(\Tt)|, \;|R(\overline H - \underline H )(\Tt)  |\}\leq  2\tau  \leq \left\{ \overline H(\Tt)-\underline H(\Tt)\right\} + \mathbf{1}_{\{\underline H(\Tt)=\overline H(\Tt)\}} .
	$$
If $\Delta$ is sufficiently small and $\mathfrak m$ sifficiently large, for $	A_N(\tau)$ defined in \eqref{def_A_N}, we have 
$$
\PP\left( \mathbf{1}_{A_N(\tau)}\neq \mathbf{1}_{\{\underline H(\Tt)<\overline H(\Tt)\}}\right)
\leq  C_3\exp\left[
	- C_5N\times \tau^2 \times \frac{\Delta^{4\overline{H}(\Tt)}\varrho(\Delta,\mathfrak m)}{\log^2(\Delta)}  \Delta^{4D(\Tt)}
	\right],
$$
where $ C_3$ and $ C_5$ are the positive constants defined as in Proposition \ref{propCH}.
\end{proposition}

For a choice of $\Delta$,  Proposition \ref{prop5_simple} allows us to determine the rate of decrease for $\tau$ such that  the indicator of $A_N(\tau)$  detects with high accuracy whether or not $H_1(\Tt)=H_2(\Tt)$. The fastest rate  depends on the approximation errors \eqref{rates_eRR}, which are characteristics of the process $X$. 


\section{Examples}\label{sec:example}

We propose two applications where our estimation approach for the local regularity for multivariate functional data opens the door to new procedures and sharp results.

\subsection{Estimating the characteristics of general Gaussian processes}\label{BfMs}

The multifractional Brownian motion (MfBm) is a generalization of the standard fractional Brownian motion, where the  Hurst parameter is allowed to vary along the path. There are several possible definitions of such a process. They lead to indistinguishable processes, up to a multiplication by a deterministic function. Here, the multi-parameter,  anisotropic multifractional Brownian sheet, which is a multivariate extension, is defined following~\citet{herbin_06}. This definition relies on the so-called harmonizable representation of the MfBm, see \cite{peltier:inria}, \cite{benassi97}, \cite{ayache2011}, \cite{lebo2018} among others. 
\begin{definition}
	Set $d\in \mathbb{N}^\star$ and let $\Eeta=(\eta_1,\dotsc,\eta_d) : [0,\infty)^d\rightarrow (0,1)^d$ be a deterministic map. The multifractional Brownian sheet $W = (W(\Uu) : \Uu \in (0,\infty)^d)$ with Hurst functional parameter $\Eeta$ is defined as follows~:
	$$ 
	W(\Uu)= \left(\prod_{k=1}^d\frac{1}{C(\eta_k(\Uu))}\right)\int_{\mathbb{R}^d}\displaystyle\prod_{k=1}^d\frac{e^{i t_k\zeta_k}-1}{ |\zeta_k|^{\eta_k(\Uu)+\frac{1}{2}}}\widehat{\boldsymbol B}(\D \boldsymbol\zeta), \qquad  \Uu\in (0, \infty)^d,
	$$
	where $\boldsymbol  \zeta =(\zeta_1,\dots,\zeta_d)$ and $\widehat{\boldsymbol B} $ is the Fourier Transform of the white noise in $\mathbb{R}^d$. Here, for any positive $x$,
	$$
	C(x) = \left[ \frac{2\pi}{\Gamma(2x+1)\sin(\pi x)}\right]^{1/2}.
	$$
\end{definition}

Notice that, when $d=1$, the measure $\widehat{\boldsymbol B}(\D \boldsymbol\zeta)$ is the unique, complex-valued Gaussian measure which can be associated with a standard Gaussian measure over $\RR$ by a `stochastic Parseval identity', see~\citet{stoev_taquu_2006}, equation~(2.4). In particular, the construction of $\widehat{\boldsymbol B}(\D \boldsymbol\zeta)$ ensures that $W$ is real-valued. 

We focus on the case $d=2$, and redefine $W=(W(\Uu) : \Uu\in\cU)$ as the restriction to an open subset $\cU\subset (0, \infty)^2$ of the multifractional Brownian sheet with Hurst functional parameter $\Eeta=(\eta_1, \eta_2)$. 
Note that $W$ is a centered Gaussian process with the covariance function
\begin{equation*}
	\EE[W(\Uu)W(\Vv)]
	\!=\! \prod_{ k =1,2}\!\!
	D(\eta_{k }(\Uu),\eta_{k }(\Vv))
	\left[u_{k  }^{\eta_{k  }(\Uu)+\eta_{ k}(\Vv)}\!+\!v_{k }^{\eta_{ k  }(\Uu)\!+\!\eta_{ k  }(\Vv)}\!-|u_{ k }-v_{k }|^{\eta_{ k }(\Uu)+\eta_i(\Vv)}\right],
\end{equation*}
$\Uu =(u_1,u_2),\Vv=(v_1,v_2)\in\cU$, where 
\begin{equation}\label{def_D_func}
	D(x,y) = C^{\,2}((x+y)/2)\cdot(2C(x)C(y))^{-1} \;\;\; \text{ and } \;\;\;  D(x,x)\equiv 1/2.
\end{equation}
In particular, the variance of $W$ is given by $\EE[W^2(\Uu)]=u_1^{2\eta_1(\Uu)}u_2^{2\eta_2(\Uu)}$.

Moreover, we consider a domain deformation $A$, that is a positive  and invertible application $A:\cT\to \cU$. Let  
$$
X= W\circ A \quad \text{ and } \quad \theta(\Tt,\Ss)= \EE\left[\{X(\Tt)-X(\Ss)\}^2\right], \quad \forall \Tt,\Ss\in \cT.
$$

\begin{proposition}\label{mprop}
	 If  $\Eeta:\cU\rightarrow (0,1)^2 $  and $A:\cT\to \cU$   are continuously differentiable,
	\begin{multline*}
		\theta(\Tt, \Ss)
		= |A_1(\Tt)|^{2H_1(\Tt)}|\partial_1A_2(\Tt)(t_1-s_1)+\partial_2A_2(\Tt)(t_2-s_2)|^{2H_2(\Tt)}\\
		+|A_2(\Tt)|^{2H_2(\Tt)}|\partial_1A_1(\Tt)(t_1-s_1)+\partial_2A_1(\Tt)(t_2-s_2)|^{2H_1(\Tt)}+O(\|\Tt-\Ss\|^2)\\
		+O(\|\Tt-\Ss\|^{2\underline{H}(\Tt)+1})+O\left(\|\Tt-\Ss\|^{2H_1(\Tt) +2H_2(\Tt)}\right),
		\qquad
		\Tt, \Ss \in\cT,
	\end{multline*}
	where $\partial_1, \partial_2$ denote the partial derivatives and 
	$$H_1=\eta_1\circ A\quad \text{and}\quad H_2=\eta_2\circ A. $$
\end{proposition}

The proof of the following corollary is immediate, and will thus be omitted.

\begin{corollary}\label{mycoro}
	Assume the conditions of Proposition\ref{mprop}, and that there exist $\rho\in(0,1)$ such that 
	$$
	0\leq \overline{H}(\Tt)-\underline{H}(\Tt)\leq \frac{1-\rho}{2}.$$
	Then  $X=W\circ A \in \mathcal{H}^{H_1,H_2}$
	with  $\mathbf L$ given  by~: 
	$$
	L_1^{(1)}(\Tt)= |A_2(\Tt)|^{2H_2(\Tt)}|\partial_1 A_1(\Tt)|^{2H_1(\Tt)},\quad L_2^{(1)}(\Tt)=|A_1(\Tt)|^{2H_1(\Tt)}|\partial_1 A_2(\Tt)|^{2H_2(\Tt)},$$
	$$L_1^{(2)}(\Tt)=|A_2(\Tt)|^{2H_2(\Tt)}|\partial_2 A_1(\Tt)|^{2H_1(\Tt)},\quad L_2^{(2)}(\Tt)=|A_1(\Tt)|^{2H_1(\Tt)}|\partial_2 A_2(\Tt)|^{2H_2(\Tt)} .
	$$
\end{corollary}

 Let us note that without domain deformation, \emph{i.e.}, when $A$ is the identity, $\mathbf L=(1,0,0,1)$.  The estimation approach introduced in Section \ref{sec_*} allows to estimate $H_1$, $H_2$ and $\mathbf L$  in general. The  estimation of the domain deformation $A$ is  investigated in the following.

\subsubsection {Estimating equations for the domain  deformation}
When one realization of the process is observed on a dense, regular grid, the estimation of the Hurst function of a multifractional Brownian motion was considered by \cite{hsing2020}. See also \cite{hsing2016}.The use of deformation to model non-stationary processes was first introduced  into the spatial statistics literature by \cite{sampson92}. One dimensional deformations behave locally  as a change of scale. In two dimension, deformations can rotate,  as well as scale local coordinates. See \cite{anderes2008}, \cite{anderes2009consistent} and \cite{Clerc2003} for more details. The fact that the deformation can rotate is mainly  related to the identification problem discussed in Section \ref{sec_*}.

As a consequence of our  new approach,  we can build a nonparametric estimator of the deformation $A$ under mild technical conditions. We consider that 
\begin{equation}\label{ini_cd}
	 \text{some   $(t_0,s_0)\in\cT$  is  given for which  $A_1(t_0,s_0)$ and $A_2(t_0,s_0)$ are known. }
\end{equation}
This  initial condition avoids identification issues arising in a fully non parametric setup. We also  assume that the time-deformation $A$ 
is such that 
\begin{equation}\label{simpl_A}
\inf_{\Tt\in\cT} A_k(\Tt) >0, \quad  \inf_{\Tt\in\cT} \partial_{i} A_k(\Tt) \geq 0 \quad \text{ and } \quad \inf_{\Tt\in\cT} \{\partial_{1} A_k(\Tt) +\partial_{2} A_k(\Tt) \} >0, \qquad    i,k=1,2. 
\end{equation}
Finally, we set $H_1(\Tt)<H_2(\Tt)$  and focus on the first coordinate $A_1$ of the deformation $A$.  By Corollary \ref{mycoro}, we have   
$$
L_1^{(1)}(\Tt)= A_2(\Tt)^{2H_2(\Tt)}\partial_1 A_1(\Tt)^{2H_1(\Tt)}.
$$
Since the variance of $X$ is given by
\begin{equation} \label{rel_v_A}
v(\Tt)=\EE[X(\Tt)^2]=A_2(\Tt)^{2H_2(\Tt)}A_1(\Tt)^{2H_1(\Tt)},
\end{equation}
it follows that 
$$ \left(\frac{L_1^{(1)}(\Tt)}{v(\Tt)}\right)^{\frac{1}{2H_1(\Tt)}}=\frac{\partial_1A_1(\Tt)}{A_1(\Tt)}.$$
Integrating both sides we obtain
$$
\log A_1(\Tt)=\int_{t_0}^{t_1} f_1(s,t_2){\D} s+h(t_2),\quad \text{ for} \quad \Tt=(t_1,t_2)\in \cT,
$$
where $h$ is a real-valued function of $t_2$ and 
$$
f_1(\Tt)= \left(\frac{L_1^{(1)}(\Tt)}{v(\Tt)}\right)^{\frac{1}{2H_1(\Tt)}}.
$$
The function   $h$ is determined by  
$$
\frac{h^\prime(t_2)}{h(t_2)}=g_1(t_0,t_2):=\left(\frac{L_1^{(2)}(t_0,t_2)}{v(t_0,t_2)}\right)^{\frac{1}{2H_1(t_0,t_2)}}.
$$
This leads us to the following estimating equation~: 
\begin{equation}\label{est_A1}
A_1(\Tt)=\lambda_1\exp\left(\int_{t_0}^{t_1} f_1(s,t_2){\D} s+\int_{s_0}^{t_2}g_1(t_0,s)\D s\right),
\quad \text{ where } \quad \lambda_1=A_1(t_0,s_0).
\end{equation}
An estimator $\widehat{A}_1$ of    the first component of the domain deformation is easily obtained  
by replacing $f_1$ and $g_1$ by their estimates in \eqref{est_A1}. Estimators of $f_1$ and $g_1$ are  naturally obtained by plugging into their expressions,  the estimators of    $L_1^{(1)}$, $L_1^{(2)}$,  $H_1$ and  an estimator $\widehat v(\Tt)$ of the variance $v(\Tt)$.

To provide a theoretical result for $\widehat{A}_1$, for simplicity, in addition to \eqref{low_thres}, we assume that constants $\underline\beta$, $\overline\beta$ are known such that 
\begin{equation}\label{simpl_bet}
0< \underline\beta \leq \min_{k=1,2} \inf_{\Tt\in\mathcal T} H_k(\Tt) \qquad \text{ and } \qquad 
 \max_{i=1,2}\sup_{\Tt \in \cT} L_1^{(i)}(\Tt) \leq \overline \beta . 
\end{equation}
We then truncate  the estimators, \emph{i.e.}, we replace  $\widehat H_k(\Tt)$ and $\widehat{L^{(i)}_1}(\Tt)$  correspondingly by 
$$
\max \{\widehat H_k(\Tt), \underline \beta \}  \qquad \text{ and  }\qquad  \min\left\{\widehat{L^{(i)}_1}(\Tt), \overline\beta\right\},\qquad \forall \Tt\in \cT,\; k,i=1,2, 
$$  
respectively. Given the relationship \eqref{rel_v_A} and  condition \eqref{simpl_A}, 
the variance $v(\Tt)$ is necessarily  bounded away from zero.
Finally,   for  the estimator of $v(\Tt)$, we assume that, a constant $C_v$ exists such that 
\begin{equation} \label{simpl_v}
\sup_{\Tt\in \cT} \mathbb E\left[ \{v(\Tt)/ \widehat v (\Tt)\}^p\right]	 < C_v^p, \qquad \forall p\geq 1 .
\end{equation}
This  condition   can be satisfied  if, for instance, a positive lower bound $\underline a$ for $A_1$ and $A_2$ is known in \eqref{simpl_A}. By \eqref{rel_v_A}, we then have 
$$
v(\Tt)> \underline v := \min (\underline a ^4,1).
$$
In this case,  $\widehat v (\Tt)$ can simply be defined as maximum between $\underline v$ and the empirical second order moment of the observable approximations $\widetilde X^{(j)}$.

Let
		$$
		F_1 :=\sup_{\Tt\in \cT}\EE\left[|\widehat f_1(\Tt)-f_1(\Tt)|\right]  ,\qquad  G_1 := \sup_{\Tt\in \cT}\EE\left[|\widehat g_1(\Tt)-g_1(\Tt)|\right]  , 
		$$ 
and $\operatorname{diam}(\mathcal T)= \sup_{\Ss^\prime,\Ss\in \cT}\|\Ss^\prime -\Ss\|$.  For the next result,  let  $\Delta = \mathfrak m^{-a}$ and $\rho(\mathfrak m) = \mathfrak m^{-b}$, with $a>0$, $b\geq 0$ and $\rho(\mathfrak m)$ introduced in Assumption \ref{ass_H2}. Moreover, let  
$$
  \chi (\Tt) = a\{2D(\Tt)+1\} - \min\{a\underline H (\Tt), b/2 \} >0. 
$$

\begin{proposition}\label{prop_def_A}
The assumptions of Propositions  \ref{propCH} and \ref{mprop}, and conditions \eqref{ini_cd}, \eqref{simpl_A}, \eqref{simpl_bet} and \eqref{simpl_v} hold true. Moreover, we assume that  constants $\mathfrak a_1$ and $\mathfrak A_1$ exist such that 
$$
\EE[X(\Tt)^{2p}]\leq \frac{p!}{2}\mathfrak a_1\mathfrak A_1^{p-2},\quad \forall p\in\{1,2\ldots\}.
$$ 
Then $F_1+G_1 <\infty$. 
Moreover, let 
\begin{equation}\label{cdt_ell}
\left[ 1/2-   \chi (\Tt)   \liminf_{\mathfrak m, N} \{\log(\mathfrak m)/\log(N)\}\right]_+<\ell < 1/2.
\end{equation}
Then, if $\mathfrak m$ and $N$ are sufficiently large, positive constants $\mathfrak C_v$, $\tilde q_1$ and $\tilde q_2$ exist such that, 
\begin{multline}		\EE \left[\left|\widehat A_1(\Tt)-A_1(\Tt) \right| \right]\\  \leq \mathfrak C_v A_1(\Tt) \operatorname{diam}(\mathcal T)  \max\left\{1,F_1,G_1\right\}  \left\{  
 \frac{a\log(\mathfrak m) }{	\mathfrak m^{  2aD  (\Tt) - \chi (\Tt)   } } 
	N^{\ell-1/2}+\tilde q_1
	\exp(-\tilde q_2N^\ell) \right\}.
\end{multline}
\end{proposition}

A similar representation can be derived for $A_2$, that is 
\begin{equation}\label{est_A2}
	A_2(\Tt)=\lambda_2\exp\left(\int_{t_0}^{t_1} f_2(s,t_2){\D} s+\int_{s_0}^{t_2}g_2(t_0,s)\D s\right),
\end{equation}
where  
$$
f_2(\Tt)= \left(\frac{L_2^{(1)}(\Tt)}{v(\Tt)}\right)^{\frac{1}{2H_2(\Tt)}},\quad g_2(t_0,t_2)= \left(\frac{L_2^{(2)}(t_0,t_2)}{v(t_0,t_2)}\right)^{\frac{1}{2H_2(t_0,t_2)}}\quad \text{and}\quad \lambda_2 =A_2(t_0,s_0).
$$
Estimators of $f_2$ and $g_2$,  are  obtained by plugging into their expressions,  the estimators of $L_2^{(i)}$,  $i=1,2$,  $H_2$, and an estimator of  $v(\Tt)$. Let $\widehat A_2$ be the estimator of $A_2$ obtained by plug-in using \eqref{est_A2}. Under the conditions of Proposition \ref{prop_def_A}, we can show that $F_2+G_2 $ is finite and derive a similar bound for the $\mathbb L^1-$risk of $\widehat A_2$. The arguments are similar and thus omitted.


\subsection{Adaptive optimal bivariate smoothing }\label{sec6}

Consider the problem of nonparame\-tric pointwise estimation of a 2-dimensional  aniso\-tropic regression function from a class of functions which are $\gamma_i-$Hölder continuous in the direction $e_i$, with $\gamma_i\in (0,1]$,  $i=1,2$. It is well-known that, under some conditions on the noise and given an iid sample of size $M_0$, 
the minimax rate of convergence for the estimation of a regression function $f$  over the class is
\begin{equation}\label{minmax_M0}
M_0^{-\frac{\boldsymbol \gamma}{2\boldsymbol \gamma+1}},
\end{equation}
where the effective smoothness $\boldsymbol {\gamma}$ is defined by the formula 
$$
\frac{1}{\boldsymbol \gamma}=\frac{1}{\gamma_1}+\frac{1}{\gamma_2}.
$$
See \cite{Ibra81}, \cite{lepSpo95}, \cite{lep_hoff}, \cite{GK},  \cite{dun_aniso}.

In the context of multivariate functional data, a natural issue is the reconstruction of the realizations of $X$ using the data.
To match the standard nonparametric regression setup, we hereafter consider the case where the set $\mathcal{H}^{H_1,H_2}$ in Definition \ref{def}, is built with the restriction $\boldsymbol L = (L_1,0,0,L_2)$. Fortunately,  the local regularity of a process $X\in \mathcal H^{H_1,H_2}$  is intrinsically linked to the regularity of the sample paths of the process.
Let  $\Tt \in\cT$   
and assume that 
\begin{equation}\label{Yor_us}
	\max_{i=1,2}\sup_{0<\Delta \leq \Delta_0}  \frac{\EE\left[\left\{X\left(\Tt-\Delta e_i/2\right)-X\left(\Tt+\Delta e_i/2\right)\right\}^{2p}\right]}{\EE\left[\left\{X\left(\Tt-\Delta e_i/2\right)-X\left(\Tt+ \Delta e_i/2\right)\right\}^2\right]^p}<\infty, \qquad \forall  p\in \mathbb N.
\end{equation}
By \citet[Theorem 2.1, page 26]{Yor}, 
almost  any realization of $X$  is locally  $\gamma_i-$Hölder  continuous in the direction $e_i$, for any order  $0\leq \gamma_i < H_i(\Tt)$.  See also Lemma SM.\ref{reg_RY_SM} in the Supplementary Material.

Let us notice that, with the simplified structure of $\boldsymbol L$ in the definition of $ \mathcal{H}^{H_1,H_2}$,  the  identification problem mentioned in Section \ref{sec_*} no longer occurs, and we have
\begin{equation*}\label{main_eq}
	\theta_{\Tt}^{(i)}(\Delta)= L_i(\Tt)\Delta^{2H_i(\Tt)}+O(\Delta^{2\overline{H}(\Tt)+\beta}),\quad i=1,2.
\end{equation*}
Following the methodology introduced in Section \ref{sec_*},  we  consider  the  estimating equations for the local regularity exponents~:
\begin{equation*}
	H_i(\Tt) = \frac{\log(\theta^{(i)}_{\Tt}(2\Delta))-\log(\theta^{(i)}_{\Tt}(\Delta))}{2\log(2)} + O(\Delta^{ \beta }),\quad i=1,2.
\end{equation*}
Applying these equations with a learning set of realizations of $X$, we get the  estimates  
$\widehat H_i(\Tt)$. 

Consider now a new  realization
$X^{new}$ of $X$, for which    $(Y^{new}_m,\Tt^{new}_m), 1\leq m\leq M_0$ with 
\begin{equation}
	Y^{new}_m=X^{new}(\Tt^{new}_m)+\varepsilon^{new}_m,\quad \quad 1\leq m\leq M_0,
\end{equation}
 are observed. 
Here, $M_0$ is a realization of the variable $M$, while the $\Tt^{new}_m$ are independent realizations of the bi-dimensional vector $\boldsymbol T$, with $M$ and $\boldsymbol T$  introduced in Section \ref{sec:data}. We use the Nadaraya-Watson estimator to estimate $X^{new}(\Tt)$, and we consider the simpler version with two bandwidths. Formally, let $K:\mathbb R ^2\to \Rplus$ be a density with the support in $[-1,1]\times [-1,1]$,  and $\mathbf{B}=\operatorname{diag}(1/h_1,1/h_2)$ a positive, $2\times 2$ bandwidth matrix. Considering the  vectors $\Tt$ and $\Tt^{new}_m$ as column matrices,  and using the rule $0/0=0$, the Nadaraya-Watson estimator is then
$$
\widehat X^{new}(\Tt;\mathbf{B} )=\sum_{m=1}^{M_0}Y^{new}_m\frac{K\left(\mathbf{B}(\Tt^{new}_m-\Tt)\right)}{\sum_{m=1}^{M_0}K\left(\mathbf{B}(\Tt^{new}_m-\Tt)\right)}
 =:\sum_{m=1}^{M_0}Y^{new}_mW^{new}_m(\Tt). 
$$
To achieve the optimal rate of convergence, the bandwidths have to be selected  according to the regularity  of $X$.  For deriving the properties of  $\widehat X^{new}(\Tt;\mathbf{B})$,  we impose the following mild assumptions.
 $B(\mathbf 0,r)$ denote the ball  centered at the origin of $\mathbb R^2$, with  radius $r$.
 
\begin{assumptionLP}
	\item\label{LP1} Two positive constants $\kappa$ and $r$ exist such that  
	$$\kappa^{-1}\mathbf{1}_{B(0,r)}(\Tt)\leq K(\Tt)\leq\kappa ,\quad \forall \Tt\in \cT. $$
	$h_1,h_2\in\mathcal H$ with $\mathcal H$, a bandwidth range satisfying $\sqrt{\mathfrak m} \inf  \mathcal H \rightarrow \infty$ and $\sup  \mathcal H \rightarrow 0$. 
	
	\item\label{LP2} A constant $c$ exists such that $f_{\mathbf{T}}(\Tt)\geq c>0$,  $\forall \Tt\in \cT$, where $f_{\mathbf{T}}$ is the density function of the random vector $\boldsymbol{T}$ that generated the independent copies  $\Tt^{new}_m$, $1\leq m\leq M_0 $.
	
	\item\label{LP3e} 
	The error terms $\varepsilon^{new}_m$ are iid,  zero mean random variables with  constant variance $\sigma^2$. The variables  $M_0$, $X^{new}$, $\Tt^{new}_m$, and $\varepsilon^{new}_m$,  $1\leq m\leq M_0 $, are mutually independent. 	  A constant $c>0$ exists such that $c^{-1}\leq M_0/\mathfrak m \leq c$.

	\item \label{LP4} The estimators $\widehat H_i(\Tt)$ and $\widehat L_i(\Tt)$ are independent of the variables $M_0$, $X^{new}$, $\Tt^{new}_m$, and $\varepsilon^{new}_m$. Moreover,  
	$$
	\PP\left(\min\left\{|\widehat H_i(\Tt)-H_i(\Tt)|, |\widehat L_i(\Tt)-L_i(\Tt)|\right\}>\log^{-a}(\mathfrak m)\right)\leq \mathfrak k_1 \exp \left(-\mathfrak m\right), \qquad i=1,2,
	$$
	where $\mathfrak k_1$ is some positive constant and $a>1$.
\end{assumptionLP}

In view of our result from Section \ref{sec4},  condition LP\ref{LP4} holds true under mild conditions. 
Let us consider the pointwise, conditional mean square risk of $\widehat X^{new}$,  given the integer $M_0$, that is
\begin{equation}\label{ave_risk}
\mathcal R \left( \Tt;\textbf B, M_0\right)=\EE\left[\left\{\widehat X ^{new}(\Tt;\textbf B)-X^{new}(\Tt)\right\}^2\Big{|} M_0\right].
\end{equation}
We first derive a bound of this risk when $H_1$, $H_2$ and $L_1,L_2$ are given.

\medskip

\begin{proposition}\label{prop_risk1}
	Assume that  (LP\ref{LP1}), (LP\ref{LP2}) and (LP\ref{LP3e})  hold true. Then 
\begin{equation}\label{risk_1}
	\mathcal{R}(\Tt; \textbf B,M_0)\leq \frac{\kappa^2}{c\pi}\frac{\sigma^2}{M_0 h_1h_2}+ 2L_1(\Tt)h_1^{2
		H_1(\Tt)}+ 2L_2(\Tt)h_2^{2H_2(\Tt)}+\text{ negligible terms.}
\end{equation}
\end{proposition}

\medskip

Minimizing the dominating terms in the risk bound yields the optimal  bandwidths. These bandwidths, and the resulting risk rate,  will depend on the regularity of the process and the Hölder constants. These facts are gathered in the following result. Before stating it, let
\begin{equation}\label{eff_smoot}
  \omega(\Tt) = \frac{H_1(\Tt) H_2(\Tt)}{H_1(\Tt)+ H_2(\Tt)}, 
\end{equation}
denote the  effective smoothness of $X$ at $\Tt$, and let 
	$$
\alpha_i(\Tt)=
\frac{\omega(\Tt)}{2\omega(\Tt)+1}\times  \frac{1}{H_i(\Tt)},\qquad i=1,2.
$$
 Moreover, let  $
\mathcal H (\Tt) = 2H_1(\Tt)H_2(\Tt)+H_1(\Tt)+H_2(\Tt).
$

\begin{corollary}\label{cor_ad}
	The minimum of the dominant terms in the risk bound in Proposition \ref{prop_risk1} is attained at $(h_1^*, h_2^*)$, with 
	$$
	h^*_1=	M_0^{ -\alpha_1(\Tt)} \left[\frac{\Lambda_1(\Tt)^{2H_2(\Tt)+1}}{\Lambda_2(\Tt)}\right]^{\frac{1}{2\mathcal H (\Tt)}}
	\quad\text{ and} \quad  
	h^*_2= M_0 ^{ -\alpha_2(\Tt)}
	\left[\frac{\Lambda_2(\Tt)^{2H_1(\Tt)+1}}{\Lambda_1(\Tt)}\right]^{\! \frac{1}{2\mathcal H (\Tt)}},
	$$
	where $\Lambda_i(\Tt)=\kappa^2\sigma^2/\{4c\pi H_i(\Tt)L_i(\Tt)\},$ $i=1,2.$
	Then, up to negligible terms,
$$
		\mathcal{R}(\Tt;\textbf B^*, M_0)\leq M_0^{-\frac{2\omega(\Tt)}{2\omega(\Tt)+1}}\Gamma_1 (\Tt),
$$
	where $\mathbf{B}^*=\operatorname{diag}(1/h^*_1,1/h^*_2)$, and
	$$
\Gamma_1(\Tt)=  \frac{\kappa^2}{\pi} \frac{\sigma^2 }{c}  \Lambda_1(\Tt)^{\frac{H_2(\Tt)}{\mathcal H (\Tt)} }\Lambda_2(\Tt)^{\frac{H_1(\Tt)}{\mathcal H (\Tt)}}
	\left\{1+2H_1(\Tt)+2H_2(\Tt)\right\}.$$
\end{corollary}

The rate of $\mathcal{R}(\Tt;\textbf B^*, M_0)$  is the minimax rate \eqref{minmax_M0} for the effective smoothness $\omega(\Tt) $. In general,  $\omega(\Tt) $ is larger than the local effective smoothness of the realization $X^{new}$, but  arbitrarily close provided  \eqref{Yor_us} is satisfied. 
Finally, we derive the bound of the average, pointwise  risk \eqref{ave_risk}  when the regularity parameters are estimated, following our methodology. Let $\widehat h_1^*$ and $\widehat h_2^*$ be the bandwidths obtained by replacing $H_i(\Tt)$ and $L_i(\Tt)$ by their estimates $\widehat H_i(\Tt)$ and $\widehat L_i (\Tt)$ in the expressions of $h_1^*$ and $h_2^*$, respectively. Let $\widehat {\textbf B}^*= \operatorname{diag}(1/\widehat h_1^*,1/\widehat h_2^*)$.

\begin{proposition}\label{risk_2}
	Assume the conditions of 
	Proposition \ref{prop_risk1} and  \ref{LP4} hold true. Then  
	$$
	\mathcal R(\Tt;\widehat {\textbf B}^*,M_0)\leq \Gamma_2(\Tt) M_0^{-\frac{2\omega(\Tt)}{2\omega(\Tt)+1}+2\log^{-a}(\mathfrak m)}\times \{1+o(\log^{-a}(\mathfrak m))\},
	$$
	where 
	$$
	\Gamma_2(\Tt)= \frac{\kappa^2\sigma^2}{c\pi\Lambda_1^{\alpha_1(\Tt)}\!(\Tt)\Lambda_2^{\alpha_2(\Tt)}\!(\Tt)}+L_1(\Tt)\left(\!\frac{\Lambda_1(\Tt)^{2H_1(\Tt)+1}}{\Lambda_2(\Tt)}\!\right)^{\alpha_1(\Tt)}+L_2(\Tt)\left(\!\frac{\Lambda_2(\Tt)^{2H_2(\Tt)+1}}{\Lambda_1(\Tt)}\!\right)^{\alpha_2(\Tt)}.
	$$
\end{proposition}

Proposition \ref{risk_2} shows that, modulo some constant terms,  the price for the estimation of the local regularity is the factor
$M_0^{2 \log^{-a} (\mathfrak m)}$, for some $a>1$. Since 
$\mathfrak m^{\log^{-1} (\mathfrak m)}=e$ for any $\mathfrak m >0$, 
the factor is essentially equal to 1 under very mild condition. 



\appendix

\section{Proofs}
Below $\sim$ means left hand side bounded above and below by constants times the right hand side.

\begin{proof}[Proof of Proposition \ref{proprox}]
	By definition, we have that 
	$$
	\gamma_{\Tt}(\Delta)=\left(K^{(1)}_1(\Tt)+K_1^{(2)}(\Tt)\right)\Delta^{2\underline H (\Tt)}
	+O(\Delta^{\widetilde \beta})=: \underline K(\Tt)\Delta^{2\underline H (\Tt)}+O(\Delta^{\widetilde \beta}).
	$$
	Moreover,  $\underline K(\Tt) = K_1(\Tt)+K_2(\Tt)$ if $\underline H (\Tt)=\overline H (\Tt)$, and $\underline K(\Tt) = K_1(\Tt)$ otherwise, with 
	$K_1(\Tt)$ and $K_2(\Tt)$ defined \eqref{eq:K1K2}. 
	We deduce
	\begin{multline}\label{eq_gamma1}
		\!\!\!\!\frac{\log(\gamma_{\Tt}(2\Delta))\! -\log(\gamma_{\Tt}(\Delta))}{2\log(2)}= \frac{\log\!\left(\!\underline K(\Tt)(2\Delta)^{2\underline H (\Tt)}\!+O(\Delta^{\widetilde \beta})\!\right)\!-\log\!\left(\!\underline K(\Tt)\Delta^{2\underline H (\Tt)}\!+O(\Delta^{\widetilde \beta})\!\right)}{2\log(2)}    \\
		= \underline H (\Tt)+\frac{\log\left(1+O(\Delta^{\widetilde \beta-2\underline H (\Tt)})\right)-\log(1+O(\Delta^{\widetilde \beta -2\underline H (\Tt)}))}{2\log(2)} 
		= \underline H (\Tt)+O(\Delta^{\widetilde \beta-2\underline H (\Tt)}),
	\end{multline}
	which gives the first part of the statement. For the second part, by the expansion \eqref{eq:K1K2},
	$$
	\gamma_{\Tt}(\Delta)= K_1(\Tt)\Delta^{2\underline H (\Tt)}+K_2(\Tt)\Delta^{2\overline H (\Tt)}+ O(\Delta^{2\overline H (\Tt)+\beta}).
	$$
	Therefore, $\alpha_{\Tt}(\Delta)$ can be written as 
	\begin{align*}
		\alpha_{\Tt}(\Delta)&=  \left|\frac{\gamma_{\Tt}(2\Delta)}{(2\Delta)^{2\underline{H}(\Tt)}}-\frac{\gamma_{\Tt}(\Delta)}{\Delta^{2\underline{H}(\Tt)}}\right|\\
		&= \left| K_2(\Tt)\left(2^{2\overline H(\Tt)-2\underline H (\Tt)}-1\right)\Delta^{2\overline H(\Tt)-2\underline H (\Tt)} +O(\Delta^{2\overline H(\Tt)-2\underline H (\Tt)+\beta})\right|\\
		& =: \left| \overline K(\Tt)\right|\Delta^{2\overline H(\Tt)-2\underline H (\Tt)} +O(\Delta^{2\overline H(\Tt)-2\underline H (\Tt)+\beta}).
	\end{align*}
	Finally, replace $\gamma_{\Tt}$ by $\alpha_{\Tt}$ in
  \eqref{eq_gamma1}, and derive the representation for $\overline H(\Tt)-\underline H (\Tt)$. \end{proof}

\begin{proof}[Proof of Proposition \ref{prop_Lest}]
	Similar to that of Proposition \ref{proprox}. 
\end{proof}

\begin{proof}[Proof of Proposition \ref{propCH}]
 We next simply write $\varrho(\Delta)$ instead of $\varrho(\Delta,\mathfrak m)$. The proof is organized in  several steps. 
First, using Assumptions \ref{ass_H1}, \ref{ass_H2} and \ref{ass_H3}, 	combined with Bernstein's inequality, a constant $\mathfrak u >0$ exists such that, for any $i=1,2$, $\varepsilon \in(0,1)$, and  $0<\Delta\leq \Delta_0$~:  
    \begin{equation}\label{eq:assumption-theta-hat_main}
	\max\left\{
	\PP\left(\widehat{\theta}_{\Tt}^{(i)}(\Delta)-\theta_{\Tt}^{(i)}(\Delta)\geq \varepsilon\right),
	\PP\left(\widehat{\theta}_{\Tt}^{(i)}(\Delta)-\theta_{\Tt}^{(i)}(\Delta)\leq -\varepsilon\right)
	\right\}
	\leq \exp\left( -\mathfrak{u}N\varepsilon^2\varrho(\Delta)\right),
\end{equation}
with $\widehat{\theta}_{\Tt}^{(i)}(\Delta)$ defined in \eqref{eq_theta_hat}, and provided that 
$\mathfrak m$ is sufficiently large.  The proof of \eqref{eq:assumption-theta-hat_main} is provided in the Supplementary Material. 

\noindent
	\textbf{\textit{Step 1 : proof of equation \eqref{eq:conc-Hhat-around-H_main}.}}
	For  $ \max \{R(\underline H )(\Tt),R(\overline H - \underline H )(\Tt)  \}\leq \varepsilon \leq 2\tau$, we have 
	$$
	\mathbb{P}\left[|\underline{\widehat{H}}(\Tt)-\underline{H}(\Tt)|\geq 2\varepsilon \right] \leq 
		\mathbb{P}\left[|\underline{\widehat{H}}(\Tt)-\underline{H}(\Tt)+R(\underline H )(\Tt) |\geq \varepsilon \right] = :A_\varepsilon
	$$
	Using the definitions and elementary inequalities,
	 we have
	\begin{align*}
		A_\varepsilon
		&=  \mathbb{P}\left[\left|\log\left(\frac{\widehat{\gamma}_{\Tt}(2\Delta)\gamma_{\Tt}(\Delta)}{\gamma_{\Tt}(2\Delta)\widehat{\gamma}_{\Tt}(\Delta)}\right)\right|\geq 2\varepsilon\log2 \right]\\
		&\leq \mathbb{P}\left[\frac{\widehat{\gamma}_{\Tt}(2\Delta)\gamma_{\Tt}(\Delta)}{\gamma_{\Tt}(2\Delta)\widehat{\gamma}_{\Tt}(\Delta)}\geq 2^{2\varepsilon}\right]+\mathbb{P}\left[\frac{\widehat{\gamma}_{\Tt}(2\Delta)\gamma_{\Tt}(\Delta)}{\gamma_{\Tt}(2\Delta)\widehat{\gamma}_{\Tt}(\Delta)}\leq 2^{-2\varepsilon}\right]\\
		&\leq  \mathbb{P}\left[\frac{\widehat{\gamma}_{\Tt}(2\Delta)}{\gamma_{\Tt}(2\Delta)}\geq 2^\varepsilon\right]+\mathbb{P}\left[\frac{\widehat{\gamma}_{\Tt}(2\Delta)}{\gamma_{\Tt}(2\Delta)}\leq 2^{-\varepsilon}\right] +\mathbb{P}\left[\frac{\widehat{\gamma}_{\Tt}(\Delta)}{\gamma_{\Tt}(\Delta)}\geq 2^\varepsilon\right]+\mathbb{P}\left[\frac{\widehat{\gamma}_{\Tt}(\Delta)}{{\gamma_{\Tt}}(\Delta)}\leq 2^{-\varepsilon}\right]\\
		&\leq 4\exp\left(-\mathfrak{u}N(2^\varepsilon-1)^2  \gamma_*(\Delta) \varrho(\Delta)\right) +4 \exp\left( - \mathfrak{u}N(1-2^{-\varepsilon})^2\gamma_*(\Delta)\varrho(\Delta) \right),
	\end{align*}
where $\gamma_*(\Delta) = \min\{\gamma^2_{\Tt}(2\Delta),\gamma^2_{\Tt}(\Delta)\}$, and the last line is a  direct consequence of~\eqref{eq:assumption-theta-hat_main}.   
	
By elementary algebra and the fact that, for small $\Delta$,  we have $\gamma_{\Tt}(\Delta) = K_1(\Tt)\Delta^{2\underline H (\Tt)}\{1+o(1)\}$,  we deduce that   positive constants  $C_1$ and $C_2$  exist such that 
	\begin{equation}\label{eq:concentration-Hhat-around-H}
		\mathbb{P}\left[
		|\underline{\widehat{H}}(\Tt)-\underline{H}(\Tt)|
		\geq \varepsilon 
		\right]
		\leq C_1\exp \left(-C_2N\varepsilon^2\Delta ^{4\underline{H}(\Tt)}\varrho(\Delta)\right).
	\end{equation}
	
	

	\noindent
	\textbf{\textit{Step 2.}} This step consists in proving that 
	constants $\tilde L_5$, $\tilde L_6$ exist such that
	\begin{equation}\label{eq:bound-alpha_main}
		\PP(|\widehat{\alpha}_{\Tt}(\Delta)-\alpha_{\Tt}(\Delta)|\geq \varepsilon)
		\leq \tilde L_5
		\exp\left(
		-\tilde L_6N\varepsilon^2\frac{\Delta^{4\overline{H}(\Tt)}\varrho(\Delta)}{\log^2(\Delta)}
		\right),
	\end{equation} 
	provided $\Delta$ is sufficiently small.
	The proof is relegated to the Supplementary Material.

	
	\noindent
	\textbf{\textit{Step 3.}} To prove equation~\eqref{eq:concentration-overlineH_main}, we recall that
	\begin{equation*}
		\overline H(\Tt) = \underline H(\Tt) + D(\Tt)
		\qquad\text{and}\qquad
		\widehat{\overline{H}}(\Tt)=\widehat{\underline{H}}(\Tt)+\widehat{D}(\Tt)\mathbf{1}_{A_N(\tau)},
	\end{equation*}
	where the event $A_N(\tau)$ is defined in~\eqref{def_A_N}, and $\widehat{D}(\Tt)$ is the estimator of the difference ${\overline{H}}(\Tt)-{\underline{H}}(\Tt)$. We simply write $A_N$ instead of $A_N(\tau)$ in the sequel. Two  cases can be distinguished~: the isotropic case, where $\underline{H}(\Tt)=\overline{H}(\Tt)$, and the anisotropic case, where $\underline{H}(\Tt)<\overline{H}(\Tt)$. In the anisotropic  situation, we use  \eqref{eq:concentration-Hhat-around-H}, with $\varepsilon $ as in \eqref{eq:cdt_eps}, to get 
	\begin{multline*}
		\PP\left[\left|\widehat{\overline{H}}(\Tt)-\overline{H}(\Tt)\right|\geq 2\varepsilon\right]
		\leq
		\PP\left[\left|\widehat{\underline{H}}(\Tt)-\underline H(\Tt)\right|\geq \varepsilon\right]
		+ \PP\left[\left|\widehat D(\Tt) - D(\Tt)\right|\geq \varepsilon\right]
		+ \PP[\bar A_N]\\
		\leq
		C_1\exp \left(-C_2N \varepsilon^2\Delta ^{4\underline{H}(\Tt)}\varrho(\Delta)\right) 
		+ \PP\left[\left|\widehat D(\Tt) - D(\Tt)\right|\geq \varepsilon\right] 
		+ \PP[\bar A_N].
	\end{multline*}
Here, for a set $A$, $\overline{A}$ denotes its complement. We now remark that,  for $\varepsilon $ as in \eqref{eq:cdt_eps},
	\begin{multline*}
\hspace{-.2cm} \PP\left[\left|\widehat D(\Tt)-D(\Tt)\right|\geq \varepsilon\right]
\leq \PP\left[\left|\log\left(\frac{\widehat{\alpha}_{\Tt}(2\Delta)\alpha_{\Tt}(\Delta)}{\widehat{\alpha}_{\Tt}(\Delta)\alpha_{\Tt}(2\Delta)}\right)\right|\geq \varepsilon\log(2)\right]\\
\leq \PP\!\left[\frac{\widehat{\alpha}_{\Tt}(2\Delta)}{\alpha_{\Tt}(2\Delta)}\geq2^{\varepsilon/2}\right]+\PP\!\left[\frac{\widehat{\alpha}_{\Tt}(\Delta)}{\alpha_{\Tt}(\Delta)}\leq2^{-\varepsilon/2}\right]+\PP\!\left[\frac{\widehat{\alpha}_{\Tt}(2\Delta)}{\alpha_{\Tt}(2\Delta)}\leq2^{-\varepsilon/2}\right]+\PP\!\left[\frac{\widehat{\alpha}_{\Tt}(\Delta)}{\alpha_{\Tt}(\Delta)}\geq 2^{\varepsilon/2}\right].
	\end{multline*}
 We focus on the first term, the other three can be bounded similarly.  By~\eqref{eq:bound-alpha_main}, for small $\Delta$,
	\begin{align*}
		\PP\left[\frac{\widehat{\alpha}_{\Tt}(2\Delta)}{\alpha_{\Tt}(2\Delta)}\geq2^{\varepsilon/2}\right]
		&=\PP\left[\widehat{\alpha}_{\Tt}(2\Delta)-\alpha_{\Tt}(2\Delta)\geq (2^{\varepsilon/2}-1)\alpha_{\Tt}(2\Delta)\right]\\
		&\leq \tilde L_5
		\exp\left[
		-\tilde L_6N\left((2^{\varepsilon/2}-1)\alpha_{\Tt}(2\Delta)\right)^2\frac{\Delta^{4\overline{H}(\Tt)}\varrho(\Delta)}{\log^2(\Delta)}
		\right].
	\end{align*}
	Since $(2^{\varepsilon/2}-1)^2 \geq \varepsilon^2\log^2(2)/4$, we obtain that
	\begin{equation*}
		\PP\left[\frac{\widehat{\alpha}_{\Tt}(2\Delta)}{\alpha_{\Tt}(2\Delta)}\geq2^{\varepsilon/2}\right]
		\leq \tilde L_5
		\exp\left[
		-\tilde L_7N\varepsilon^2\frac{\Delta^{4\overline{H}(\Tt)}\varrho(\Delta)}{\log^2(\Delta)}\Delta^{4D(\Tt)}
		\right],
	\end{equation*}
	for some positive constant $\tilde L_7$. The same inequality, with possibly different constants, remains valid for the three other terms. Therefore, a  constant $\tilde L_8$ exists such that,  for $\varepsilon $ as in \eqref{eq:cdt_eps},
	\begin{equation}\label{eq:concentration-Difference}
		\PP\left[\left|\widehat{D}(\Tt)-D(\Tt)\right|\geq \varepsilon\right]
		\leq \tilde L_8
		\exp\left[
		-\tilde L_7N\varepsilon^2\frac{\Delta^{4\overline{H}(\Tt)}\varrho(\Delta)}{\log^2(\Delta)}\Delta^{4D(\Tt)}
		\right].
	\end{equation}
	
	Finally, it remains to bound $\PP[\bar A_N]$. Since $\tau \leq D(\Tt)/2$, we obtain 
	\begin{equation*}
		\PP\left(\bar A_N\right)=\PP\left(\widehat D (\Tt)- D(\Tt) \leq \tau -D(\Tt)\right)\leq \PP\left(\widehat D (\Tt)- D(\Tt) \leq -\tau \right).
	\end{equation*}
	Using \eqref{eq:concentration-Difference} with $2\tau$ in place of $\varepsilon$ (which is allowed by the condition $\varepsilon \leq 2\tau$)   leads to 
	\begin{equation}\label{A_N_aniso}
		\PP[\bar A_N] \leq \tilde L_8
		\exp\left[
		-4\tilde L_7N\tau^2\frac{\Delta^{4\overline{H}(\Tt)}\varrho(\Delta)}{\log^2(\Delta)}\Delta^{4D(\Tt)}
		\right].
	\end{equation}
	This implies, for $\varepsilon $ as in \eqref{eq:cdt_eps},
	\begin{multline}\label{eq:anisotropic}
		\PP\left[\left|\widehat{\overline{H}}(\Tt)-\overline{H}(\Tt)\right|\geq \varepsilon\right]
		\leq
		  C_1\exp \left[-C_2N (\varepsilon^2/4)\Delta ^{4\underline{H}(\Tt)}\varrho(\Delta)\right] 
		\\ + \tilde L_8
		\exp\left[
		-\tilde L_7N(\varepsilon^2/4)\frac{\Delta^{4\overline{H}(\Tt)}\varrho(\Delta)}{\log^2(\Delta)}\Delta^{4D(\Tt)}
		\right]
		+ \tilde L_8
		\exp\left[
		-4\tilde L_7N\tau^2\frac{\Delta^{4\overline{H}(\Tt)}\varrho(\Delta)}{\log^2(\Delta)}\Delta^{4D(\Tt)}
		\right],
	\end{multline}
and this concludes the proof  of the anisotropic case. 

For the isotropic case, where $\underline H(\Tt)=\overline H(\Tt)$, we use \eqref{eq:concentration-Hhat-around-H} and decompose as follows~: 
	\begin{multline*}
		\PP\left[\left|\widehat{\overline{H}}(\Tt)-\overline{H}(\Tt)\right|\geq \varepsilon\right] \leq  \PP\left[\left|\widehat{\underline{H}}(\Tt)-\underline{H}(\Tt)\right|\geq \varepsilon\right]+ \PP[A_N]\\
		\leq C_1\exp \left[-C_2N\varepsilon^2\Delta ^{4\underline{H}(\Tt)}\varrho(\Delta)\right] +\PP[A_N].
	\end{multline*}
We now have to bound $\PP[A_N]$, instead of  $\PP[\Bar A_N]$. For this, we can simply write 
	\begin{equation*}
		\PP[A_N]= \PP[\widehat D (\Tt)\geq \tau]=\PP[\widehat D (\Tt)- D(\Tt)\geq \tau].
	\end{equation*}
	Using \eqref{eq:concentration-Difference}  with $2\tau$ in place of $\varepsilon$ we then obtain 
	\begin{equation}\label{A_N_iso}
		\PP[ A_N] \leq \tilde L_8
		\exp\left[
		-4\tilde L_7N\tau^2\frac{\Delta^{4\overline{H}(\Tt)}\varrho(\Delta)}{\log^2(\Delta)}\Delta^{4D(\Tt)}
		\right],
	\end{equation}
	which leads to
	\begin{multline}\label{eq:isotropic}
		\PP \left[\left|\widehat{\overline{H}}(\Tt) - \overline{H}(\Tt)\right| \geq \varepsilon\right] 
	 	\leq   C_1\exp  \left[
		-C_2N\varepsilon^2\Delta ^{ 4\underline{H}(\Tt)}\varrho(\Delta)
		\right] 
		\\ +  \tilde L_8   \exp \left[ 
		-4 \tilde L_7N\tau^2\frac{\Delta^{ 4\overline{H}(\Tt)}\varrho(\Delta)}{\log^2(\Delta)}\Delta^{ 4D(\Tt)}
		\right] .
	\end{multline}
	Combining  \eqref{eq:anisotropic} and \eqref{eq:isotropic}, three positive constants $L_3,L_4$ and $L_5$ exist such that  
	\begin{multline*}
		\PP\left[\left|\widehat{\overline{H}}(\Tt)-\overline{H}(\Tt)\right|\geq \varepsilon\right] 
		\leq L_3\Big\{\exp \left[
		-L_2N\varepsilon^2\Delta ^{4\underline{H}(\Tt)}\varrho(\Delta)
		\right] 
		\\ \left. +\exp\left[
		- L_4N\tau^2\Delta^{4\overline{H}(\Tt)}\varrho(\Delta)\log^{-2}(\Delta)\Delta^{4D(\Tt)}
		\right]+P\right\},
	\end{multline*}
	 for any $\varepsilon $ as in \eqref{eq:cdt_eps}, where 
	$$
	P= \exp\left[
	- L_5N\varepsilon^2\Delta^{4\overline{H}(\Tt)}\varrho(\Delta)\log^{-2}(\Delta)\Delta^{4D(\Tt)}
	\right]\mathbf1_{\{\underline H(\Tt)<\overline H (\Tt)\}}.
	$$
	Let us note that $P$ is a term which only occurs in the anisotropic case. 
	\end{proof}

\begin{proof}[Proof of Proposition \ref{conc_Lest}] 

	\textbf{\textit{Proof of  \eqref{eq:conc_L1_main}}.} Here is  the anisotropic case,  and we consider $H_1(\Tt)=\underline H (\Tt)$. Set   $\varepsilon$ as in \eqref{eq:cdt_epsL}  and, for $i=1,2$,   define  
\begin{equation*}
\PP\left(\left|\widehat{L_1^{(i)}}(\Tt)-L_1^{(i)}(\Tt)\right|\geq 2\varepsilon  \right) \leq 
\PP\left(\left|\widehat{L_1^{(i)}}(\Tt)-L_1^{(i)}(\Tt)+ R(L_1^{(i)})(\Tt)\right|\geq \varepsilon  \right)=:\mathfrak A_{\varepsilon}^{(i)}.
\end{equation*}
Using the definition of $\widehat{L_1^{(i)}}(\Tt)$ we can decompose~:
\begin{equation}\label{eq:decom-holder-const}
	\mathfrak A_{\varepsilon}^{(i)} \leq 
			\PP\left(\frac{\widehat{\theta}_{\Tt}^{(i)}(\Delta)}{\theta_{\Tt}^{(i)}(\Delta)}\geq\left(1+ \varepsilon\frac{\Delta^{2\underline{H}(\Tt)}}{\theta_{\Tt}^{(i)}(\Delta)}\right)^{\!\frac{1}{2}}\right)
			+ \PP\left(\frac{\Delta^{2\widehat{\underline{H}}(\Tt)}}{\Delta^{2\underline{H}(\Tt)}}\leq\left(1+ \varepsilon\frac{\Delta^{2\underline{H}(\Tt)}}{\theta_{\Tt}^{(i)}(\Delta)}\right)^{\!-\frac{1}{2}}\right).
\end{equation}
To bound these two terms, we first notice that a constant $K$ exists such that,  $\forall \varepsilon$ as in \eqref{eq:cdt_epsL},
\begin{equation*}
\left(1+ \varepsilon\frac{\Delta^{2\underline{H}(\Tt)}}{\theta_{\Tt}^{(i)}(\Delta)}\right)^{\frac{1}{2}}\geq
1+\varepsilon \frac{1}{2\sqrt{ 1+ L_1^{(i)} (\Tt)}}+\varepsilon O(\Delta^{2\overline{H}(\Tt)-\underline H(\Tt)})
\geq 1+\varepsilon K,
\end{equation*}
  provided $\Delta$ is sufficiently small. Then, by  \eqref{eq:assumption-theta-hat_main},  for sufficiently large $\mathfrak m$,  
\begin{equation*}
    \PP\left(\frac{\widehat{\theta}_{\Tt}^{(i)}(\Delta)}{\theta_{\Tt}^{(i)}(\Delta)}\geq\left(1+ \varepsilon\frac{\Delta^{2\underline{H}(\Tt)}}{\theta_{\Tt}^{(i)}(\Delta)}\right)^{\frac{1}{2}}\right) 
    \leq \exp \left(-\mathfrak{u}KN\varepsilon^2\{ \theta_{\Tt}^{(i)}(\Delta)\}^2\varrho(\Delta)\right).
\end{equation*}
Since $\theta_{\Tt}^{(i)}(\Delta)\sim L_1^{(i)}(\Tt)\Delta^{2\underline H (\Tt)}$ we obtain~:
\begin{equation}\label{eq:first-term-L}
  \PP\left(
  	\frac{\widehat{\theta}_{\Tt}^{(i)}(\Delta)}{\theta_{\Tt}^{(i)}(\Delta)}\geq\left(1+ \varepsilon\frac{\Delta^{2\underline{H}(\Tt)}}{\theta_{\Tt}^{(i)}(\Delta)}\right)^{\frac{1}{2}}
	\right)
 \leq\exp \left(
 	-\mathcal L N\varepsilon^2\Delta^{4\underline H (\Tt)}\varrho(\Delta)
	\right),
\end{equation}
for some  constant $\mathcal L$.  For the second term, since $\Delta$ is  small, and $\log(x)\leq x-1$, $x>0$,  
\begin{equation*}
    \PP\!\left[\!\frac{\Delta^{2\widehat{\underline{H}}(\Tt)}}{\Delta^{2\underline{H}(\Tt)}}\leq\! \left\{\!1\!+ \varepsilon\frac{\Delta^{2\underline{H}(\Tt)}}{\theta_{\Tt}^{(i)}(\Delta)}\right\}^{\!\!-\frac{1}{2}}\right] \!\leq \PP\!\left[\!\frac{\Delta^{2\widehat{\underline{H}}(\Tt)}}{\Delta^{2\underline{H}(\Tt)}}\leq 1\! - \! K \varepsilon\right]
    \!\leq \PP\left[\!\widehat{\underline{H}}(\Tt)\!-\underline{H}(\Tt)\!\geq - \frac{K/2}{\log(\Delta)} \varepsilon\right]\!,
\end{equation*}
provided $\varepsilon$ satisfies \eqref{eq:cdt_epsL}.
Using \eqref{eq:concentration-Hhat-around-H}, we get
\begin{equation}\label{eq:second-term-L} 
\PP\left(
	 \frac{\Delta^{2\widehat{\underline{H}}(\Tt)}}{\Delta^{2\underline{H}(\Tt)}}\leq\left(1+ \varepsilon\frac{\Delta^{2\underline{H}(\Tt)}}{\theta_{\Tt}^{(i)}(\Delta)}\right)^{-\frac{1}{2}}
	  \right) 
\leq   C_1\exp\left(
	-\frac{C_2K^2}{4}N\varepsilon^2\frac{\Delta^{4\underline H(\Tt)}\varrho(\Delta)}{\log^2( \Delta)}
	\right). 
\end{equation}
Finally, combining \eqref{eq:decom-holder-const}, \eqref{eq:first-term-L}, \eqref{eq:second-term-L}, and considering, without loss of generality, $\Delta_0\leq e^{-1}$ in Definition \ref{def}, positive constants $ \mathfrak C _1$ and  $ \mathfrak C_2$ exist such that,    $\forall \varepsilon$ as in \eqref{eq:cdt_epsL}, we have~:
\begin{equation*}
\mathfrak A_{\varepsilon}^{(i)} \leq  \mathfrak C_1 
\exp\left(
	- \mathfrak C_2 N \varepsilon^2\frac{\Delta^{4\underline H(\Tt)}\varrho(\Delta)}{\log^2(\Delta)}
	\right),\qquad i=1,2.
\end{equation*} 

\noindent \textbf{\textit{Proof of \eqref{eq:conc_L2_main}:}} 
For $\Tt\in \cT$ and $i=1,2$, let
\begin{equation*}
\alpha_{\Tt}^{(i)}(\Delta)= \left|\frac{\theta_{\Tt}^{(i)}(2\Delta)}{(2\Delta)^{2\underline H(\Tt)}}-\frac{\theta_{\Tt}^{(i)}(\Delta)}{\Delta^{2\underline H(\Tt)}}\right|
\quad \text{ and } \quad
\widehat{ \alpha}_{\Tt}^{(i)}(\Delta)= \left|\frac{\widehat{\theta}_{\Tt}^{(i)}(2\Delta)}{(2\Delta)^{2\widehat{\underline H}(\Tt)}}-\frac{\widehat{\theta}_{\Tt}^{(i)}(\Delta)}{\Delta^{2\widehat{\underline H}(\Tt)}}\right|.
\end{equation*}
For any $\varepsilon$ such that  $ |R(L_2^{(i)})(\Tt)|\leq \varepsilon$ (recall $|R(L_2^{(i)})(\Tt)|=O(\Delta^\beta)$)  and $i=1,2$, we decompose as follows~:
\begin{equation*}
 \PP\left(\widehat{L_2^{(i)}}(\Tt)-L_2^{(i)}(\Tt)\geq 2 \varepsilon,  \widehat D(\Tt) \neq 0 \right)\leq B^{(i)}_1+B^{(i)}_2+B^{(i)}_3+B^{(i)}_4,
\end{equation*}
 where 
 \begin{equation*}
 B^{(i)}_1=\PP\!\left(\!\widehat{\alpha}_{\Tt}^{(i)}(\Delta)-\alpha_{\Tt}^{(i)}(\Delta)\geq \sqrt{\varepsilon/3}\right),\quad 
 B^{(i)}_2=\PP\!\left(\frac{\widehat{\alpha}_{\Tt}^{(i)}(\Delta)-\alpha_{\Tt}^{(i)}(\Delta)}{\{4^{D(\Tt)}-1\}\Delta^{2D(\Tt)}}\geq \varepsilon/3\!\right),
 \end{equation*}
 \begin{align}\label{B3B4_main}
 B^{(i)}_3&=\PP\left(\alpha_{\Tt}^{(i)}(\Delta)\left(\frac{1}{\{4^{\widehat{D}(\Tt)}-1\}\Delta^{2\widehat{D}(\Tt)}} -\frac{1}{\{4^{D(\Tt)}-1\}\Delta^{2D(\Tt)}}\right)\geq \varepsilon/3 ,\widehat D(\Tt) \neq 0 \right),\\
 B^{(i)}_4&=\PP\left(\frac{1}{\{4^{\widehat{D}(\Tt)}-1\}\Delta^{2\widehat{D}(\Tt)}} -\frac{1}{\{4^{D(\Tt)}-1\}\Delta^{2D(\Tt)}}\geq \sqrt{ \varepsilon/3}, \widehat D(\Tt) \neq 0 \right).\notag
 \end{align}
 By the arguments used for \eqref{eq:bound-alpha_main}, $\mathfrak m$ sufficiently large,  constants $\tilde C_1$ and $\tilde C_2$ exists  such that, 
$\forall \varepsilon$ such that $|\log(\Delta)| |R(\underline H )(\Tt)| \leq \varepsilon$ (recall $|R(\underline H )(\Tt)|=O(\Delta^{2{D}(\Tt)})$ in the anisotropic case),   
\begin{multline}\label{eq:bound-B1}
B^{(i)}_1 \leq   \tilde C_1 \exp  \left[
	  -\tilde C_2   N\varepsilon\frac{\Delta^{4\overline H (\Tt)}\varrho(\Delta)}{\log^2(\Delta)}
	 \right]  \\
\text{ and }\qquad 
B_2^{(i)}  \leq   \tilde C_1 \exp \left[ 
	- \tilde C_2  N\varepsilon^2(4^{D(\Tt)}  -1)^2 \Delta^{4D(\Tt)}\frac{\Delta^{4\overline H (\Tt)}\varrho(\Delta)}{\log^2(\Delta)}
	  \right] .
\end{multline}
To bound $B^{(i)}_3$ and $B^{(i)}_4$ in \eqref{B3B4_main}, we  use the fact that $\alpha_{\Tt}^{(i)}(\Delta)\sim L_2^{(i)}(\Tt)(4^{D(\Tt)}-1)\Delta^{2D(\Tt)}$,  Lemma SM.\ref{lemma_L2_use} in the Supplement, and the fact that  $L_2^{(1)}$, $L_2^{(2)}$ are bounded functions. Moreover, we show that the probability of the event $\{\widehat D(\Tt)=0\}$ is negligible, see \eqref{bound_Dhat_negli}.  The details are given in the Supplementary. From that and  \eqref{eq:bound-B1}, the proof follows. \end{proof}


\begin{proof}[Proof of Proposition \ref{prop5_simple}]
	A direct consequence of \eqref{A_N_aniso}  and \eqref{A_N_iso}.
\end{proof}



\begin{proof}[Proof of Proposition \ref{mprop}]
	First, let us denote  $$B(\Tt,\Ss) =2D(H_1(\Tt),H_1(\Ss))D(H_2(\Tt),H_2(\Ss))\qquad \forall \Tt,\Ss\in \cT,$$ with $D(x,y)$ defined in \eqref{def_D_func}. By construction the function $B(\cdot,\cdot)$ is symmetric. Moreover, we show in the Supplementary Material that 
\begin{equation}\label{eq:B-diag}
	B(\Tt,\Ss)= \frac{1}{2}+O(\|\Tt-\Ss\|^2).
\end{equation}
	
Next,  using the covariance function structure of the process $X$, we can write~: 
	\begin{multline}
			\theta(\Tt,\Ss)= \EE[X^2(\Tt)]+\EE[X^2(\Ss)]- 2\EE[X(\Tt)X(\Ss)]\\
			=  \mathfrak B_1(\Tt,\Ss) +\mathfrak B_1(\Ss,\Tt) +\mathfrak B_2(\Tt,\Ss)+\mathfrak B_3(\Tt,\Ss)-\mathfrak B_4(\Tt,\Ss), \quad \forall \Tt,\Ss\in \cT,
		\end{multline}
		where 
		\begin{align*}
			\mathfrak B_1 (\Tt,\Ss)&= |A_1(\Tt)|^{2H_1(\Tt)}|A_2(\Tt)|^{2H_2(\Tt)}-B(\Tt,\Ss)|A_1(\Tt)|^{H_1(\Tt)+H_1(\Ss)}\\
			&\times\left(|A_2(\Tt)|^{H_2(\Tt)+H_2(\Ss)}\!+|A_2(\Ss)|^{H_2(\Tt)+H_2(\Ss)}\right),\\
			\mathfrak B_2(\Tt,\Ss) &= B(\Tt,\Ss)\left(|A_1(\Tt)|^{H_1(\Tt)+H_1(\Ss)}\!+|A_1(\Ss)|^{H_1(\Tt)+H_1(\Ss)}\right)|A_2(\Tt)\!-\!A_2(\Ss)|^{H_2(\Tt)+H_2(\Ss)},\\
			\mathfrak B_3(\Tt,\Ss)&=B(\Tt,\Ss)\left(|A_2(\Tt)|^{H_2(\Tt)+H_2(\Ss)}\!+|A_2(\Ss)|^{H_2(\Tt)+H_2(\Ss)}\right)|A_1(\Tt)\!-\!A_1(\Ss)|^{H_1(\Tt)+H_1(\Ss)},\\
			\mathfrak B_4(\Tt,\Ss)&= B(\Tt,\Ss)|A_1(\Tt)\!-\!A_1(\Ss)|^{H_1(\Tt)+H_1(\Ss)}|A_2(\Tt)-A_2(\Ss)|^{H_2(\Tt)+H_2(\Ss)}.
		\end{align*}
		Let 
		\begin{equation}
			a(\Tt,\Ss)=\frac{|A_1(\Tt)|^{H_1(\Tt)-H_1(\Ss)}|A_2(\Tt)|^{H_2(\Tt)-H_2(\Ss)}- B(\Tt,\Ss)}{B(\Tt,\Ss)}.
		\end{equation}
		We then have~:
		\begin{multline}\label{B1+B1}
			\mathfrak B_1(\Tt,\Ss) +\mathfrak B_1(\Ss,\Tt)=B(\Tt,\Ss) \left(a(\Tt,\Ss)|A_2(\Tt)|^{H_2(\Tt)+H_2(\Ss)}-|A_2(\Ss)|^{H_2(\Tt)+H_2(\Ss)}\right)\\
			\times\left(|A_1(\Tt)|^{H_1(\Tt)+H_1(\Ss)}-a(\Ss,\Tt)|A_1(\Ss)|^{H_1(\Tt)+H_1(\Ss)}\right)\\
			+|A_1(\Ss)|^{H_1(\Tt)+H_1(\Ss)}|A_2(\Tt)|^{H_2(\Tt)+H_2(\Ss)}B(\Tt,\Ss) \left\{ a(\Tt,\Ss)a(\Ss,\Tt) -1\right\}.
		\end{multline}
Using \eqref{eq:B-diag}, we show in the Supplementary Material that
		\begin{equation}\label{eq:double-a}
	a(\Tt,\Ss)a(\Ss,\Tt)-1 =O(\|\Tt-\Ss\|^2).
\end{equation}	
We can next deduce that 
\begin{align*}
	B(\Tt,\Ss) &\left(a(\Tt,\Ss)|A_2(\Tt)|^{H_2(\Tt)+H_2(\Ss)}-|A_2(\Ss)|^{H_2(\Tt)+H_2(\Ss)}\right)\times\\
	&\left(|A_1(\Tt)|^{H_1(\Tt)+H_1(\Ss)}-a(\Ss,\Tt)|A_1(\Ss)|^{H_1(\Tt)+H_1(\Ss)}\right)=O(\|\Tt-\Ss\|^2).
\end{align*}
From this and \eqref{eq:double-a}, equation \eqref{B1+B1} becomes~:
\begin{equation*}
	\mathfrak B_1(\Tt,\Ss) +\mathfrak B_1(\Ss,\Tt)= O(\|\Tt-\Ss\|^2).
\end{equation*}
For the terms $B_3(\Tt,\Ss)$ and $B_4(\Tt,\Ss)$, we apply \eqref{eq:B-diag}. Finally,  we get 
\begin{multline*}
	\EE[(X(\Tt)-X(\Ss))^2] \\ =\mathfrak B_4(\Tt,\Ss)+ \frac{1}{2}\left(|A_1(\Tt)|^{H_1(\Tt)+H_1(\Ss)}+|A_1(\Ss)|^{H_1(\Tt)+H_1(\Ss)}\right)|A_2(\Tt)-A_2(\Ss)|^{H_2(\Tt)+H_2(\Ss)}\\
	+\frac{1}{2}\! \left(\!|A_2(\Tt)|^{H_1(\Tt)+H_1(\Ss)}+|A_2(\Ss)|^{H_1(\Tt)+H_1(\Ss)}\!\right)\!|A_1(\Tt)\!-A_1(\Ss)|^{H_1(\Tt)+H_1(\Ss)}
	+O(\|\Tt-\Ss\|^2).
\end{multline*}
The last expression and the Taylor expansion then imply
\begin{align*}
	\EE[(X(\Tt)-X(\Ss))^2]=&|A_1(\Tt)|^{2H_1(\Tt)}|\partial_1 A_2(\Tt)(t_1-s_1)+\partial_2 A_2(\Tt)(t_2-s_2)|^{2H_2(\Tt)}\\
	+&|A_2(\Tt)|^{2H_2(\Tt)}|\partial _1A_1(\Tt)(t_1-s_1)+\partial_2 A_1(\Tt)(t_2-s_2)|^{2H_1(\Tt)}\\
	+&O(\|\Tt-\Ss\|^2)+O(\|\Tt-\Ss\|^{2\underline{H}(\Tt)+1})+ O\left(\|\Tt-\Ss\|^{2\underline H(\Tt)+2\overline H (\Tt)}\right).
\end{align*}
\end{proof}


\begin{proof}[Proof of  Proposition \ref{prop_def_A}]
We  start by showing that 
\begin{equation}\label{eq_f1_g1}
\underset{\Tt\in \cT}{\sup}\EE\left[ |\widehat f_1(\Tt)-f_1(\Tt)|^2\right]<\infty \qquad \text{ and }\qquad \underset{\Tt\in \cT}{\sup}\EE\left[ |\widehat g_1(\Tt)-g_1(\Tt)|^2\right]<\infty.
\end{equation}
Since 
$
\EE|\widehat f_1(\Tt)\! -f_1(\Tt)|^2\leq 2f_1^2(\Tt)\!+2\EE|\widehat f_1(\Tt)|^2
$, it suffices to bound $f_1^2(\Tt)$ and $\EE|\widehat f_1(\Tt)|^2$.  
By \eqref{simpl_A} and \eqref{simpl_bet},
\begin{equation*}
f_1^2(\Tt)= \left(L_1^{(1)}(\Tt) \big/v(\Tt)\right)^{\frac{1}{H_1(\Tt)}}\leq \left(\overline \beta / \underline v\right)^{\frac{1}{\underline \beta}}.
\end{equation*}
Moreover,
\begin{equation*}
\widehat f_1^2(\Tt)=\left( \widehat L_1^{(1)}(\Tt)\big / \widehat v(\Tt) \right) ^{\frac{1}{\widehat H_1(\Tt)}} \leq \left( \overline \beta/\underline v\right)^{\frac{1}{\underline \beta}}\times \left\{v(\Tt)/\widehat v(\Tt)\right\}^{\frac{1}{\underline \beta}}.
\end{equation*}
Therefore, by  \eqref{simpl_v} we obtain 
\begin{equation*}
\EE |\widehat f_1(\Tt)|^2 \leq  \ \left( \overline \beta/\underline v\right)^{\frac{1}{\underline \beta}} \times\EE\left[ \left\{v(\Tt)/\widehat v(\Tt)\right\}^{\frac{1}{\underline \beta}}\right]<\infty.
\end{equation*}
The first part of \eqref{eq_f1_g1} follows, the second part can be obtained with similar arguments.  

Next, by the condition \eqref{simpl_v} and Fubini's Theorem, a constant $\mathfrak C_v$ exists such that 
	\begin{multline}\label{deux_int}
		\EE\left[|A_1(\Tt)-\widehat A_1(\Tt)|\right]\leq \mathfrak C_v A_1(\Tt)\left[  \int_{t_0}^{t_1} \EE|\widehat f_1(s,t_2)-f_1(s,t_2)|{\D} s  \right. \\ \left. + \int_{s_0}^{t_2}\EE|\widehat g_1(t_0,s)-g_1(t_0,s)|\D s \right].
	\end{multline}
A detailed justification of the last inequality is provided in the Supplementary Material. We next bound the two integrals in \eqref{deux_int}. For $\lambda \in (0,1)$, we define the set 
	$$
	\mathcal O(\lambda)=\left\{|\widehat f_1(s,t_2)-f_1(s,t_2)|\leq \lambda \right\}.
	$$
By Cauchy-Schwarz inequality,
	\begin{multline}\label{exp}
		\int_{t_0}^{t_1}\EE\left[ |\widehat f_1(s,t_2)-f_1(s,t_2)|\right]{\D} s\leq \int_{t_0}^{t_1}\left(\lambda+\EE\left[ |\widehat f_1(s,t_2)-f_1(s,t_2)|\mathbf{1}_{\overline{\mathcal O}(\lambda)}\right]\right){\D} s\\
		\leq \int_{t_0}^{t_1}\left(\lambda+\EE\left[ |\widehat f_1(s,t_2)-f_1(s,t_2)|^2\right]^{1/2}\PP^{1/2}(\overline{\mathcal O}(\lambda))\right){\D} s.
	\end{multline}
We now want to apply  Lemma SM.\ref{lem} in the Supplementary Material with
$$
a_N=\widehat{ L_1^{(1)}}(\Tt), \quad b_N=\widehat v(\Tt), \quad c_N= \frac{1}{\widehat H_1(\Tt)} \quad \text{and } \quad 
a=L_1^{(1)}(\Tt), \quad b= v(\Tt), \quad c= \frac{1}{ H_1(\Tt)}.
$$
and $\mathfrak C= \overline \beta/\underline v$.  We first note that, by Proposition \ref{conc_Lest} with $i=1$, for any $\varepsilon$ satisfying \eqref{eq:cdt_epsL},   
$$
\PP\left(\left|\widehat{L_1^{(1)}}(\Tt)-L_1^{(1)}(\Tt)\right|\geq \varepsilon  \right)\leq \mathfrak{C}_1\exp\left(-\mathfrak{C}_2N\varepsilon^2\frac{\Delta ^{4H_1(\Tt)}\varrho(\Delta,\mathfrak m)}{\log^2(\Delta)}\right).
$$
Moreover, by Lemma SM.\ref{concentration:variance}, and for sufficiently large $\mathfrak m$, there exists a constant $\mathfrak e$ such that 
$$
  \forall \eta \in (0,1),\qquad \PP(|\widehat v(\Tt)-v(\Tt)|\geq \eta)\leq 2\exp(- \mathfrak  e N\eta ^2). 
$$
It remains to derive a bound for the concentration of $c_N$, which follows from that on the concentration of $\widehat{H}_1(\Tt)$, and the fact that $H_1\leq 1$~: 
 for any $\varepsilon$ satisfying \eqref{eq:cdt_epsL},  and thus \eqref{eq:cdt_eps},  
\begin{multline*}
	\PP\left(\pm \{c_N-c \} \geq\varepsilon\right)
	=\PP\left( \pm \{H_1(\Tt)-\widehat{H}_1(\Tt)\}\{1\pm \varepsilon H_1(\Tt)\}\geq \varepsilon H_1^2(\Tt) \right)\\
	\leq \PP\left(\pm \{H_1(\Tt)-\widehat{H}_1(\Tt)\}\geq  H_1^2(\Tt)\varepsilon/\{1+ \varepsilon H_1(\Tt)\}\right),
\end{multline*}
and thus, by Proposition \ref{propCH}, 
$$
\PP\left(|c_N-c | \geq\varepsilon\right) \leq \widetilde C_1\exp\left(-\widetilde C_2N\left\{\underline \beta^4/4\right\}\varepsilon^2\Delta^{4H_1(\Tt)}\varrho(\Delta,\mathfrak m)\right) .
$$
Finally, by Lemma SM.\ref{lem} we obtain 
	$$
	\PP(\overline{\mathcal O}(\lambda))= \PP\!\left (|\widehat f_1(s,t_2)\!-f_1(s,t_2)|\geq \lambda\!\right)\!\leq q_1 \exp\left\{\! -q_2N\lambda^2\Delta^{2H_1(s,t_2)}\!\varrho(\Delta,\mathfrak m)\log^{-2}(\Delta) \right\}\!,
	$$
	for some  constants $q_1$, $q_2$, provided  \eqref{eq:cdt_epsL} is satisfied with $\varepsilon = \lambda$. By \eqref{exp}, we now  write 
	\begin{multline*}
		\int_{t_0}^{t_1}\EE\left[ |\widehat f_1(s,t_2)-f_1(s,t_2)|\right]{\D} s
		\\
	\hspace{-1cm} 	\leq  \int_{t_0}^{t_1}\left(\lambda+\underset{\Tt\in \cT}{\sup}\EE\left[ |\widehat f_1(\Tt)-f_1(\Tt)|^2\right]^{1/2} q_1\exp\left\{ -q_2N\lambda^2\frac{\Delta^{2H_1(s,t_2)}\varrho(\Delta,\mathfrak m)}{\log^2(\Delta)} \right\}\right){\D} s\\
		\leq\operatorname{ diam }(\cT)\max\left\{1,\underset{\Tt\in \cT}{\sup}\EE\left[ |\widehat f_1(\Tt)-f_1(\Tt)|^2\right]^{1/2}\right \}\!\left [\lambda+q_1\exp\left\{ -q_2N\lambda^2\frac{\Delta^{2}\varrho(\Delta,\mathfrak m)}{\log^2(\Delta)} \right\}\right].
	\end{multline*}
A simple choice of $\lambda$ can be defined as follows~: for some suitable   $\ell\in (0,1/2)$, let
	$$
	\lambda=\frac{|\log(\Delta)|}{\Delta \sqrt{\varrho(\Delta,\mathfrak m)}}N^{\ell-1/2},
	$$
	 which satisfies \eqref{eq:cdt_epsL}, provided that $\ell $ satisfies \eqref{cdt_ell}.
 We then obtain 
	$$
	\int_{t_0}^{t_1}\EE\left[ |\widehat f_1(s,t_2)-f_1(s,t_2)|\right]{\D} s\leq C \left( \frac{|\log(\Delta)|}{\Delta \sqrt{\varrho(\Delta,\mathfrak m)}}N^{\ell-1/2}+q_1e^{-q_2N^\ell}\right),
	$$
	for some constant $C$. Up to a change of the constants $C$, $q_1$, $q_2$, a similar bound holds true for the second integral on the RHS on \eqref{deux_int}.  It remains to replace   $\Delta $ and $\rho(\mathfrak m) $ by $\mathfrak m^{-a}$ and $\mathfrak m ^{-b}$, respectively.  
\end{proof}


\begin{proof}[Proof of Proposition \ref{prop_risk1}.]
	The risk $\mathcal R(\Tt; \boldsymbol B ,M_0)$ is the sum of the squared bias and the variance, for which we derive the 
 following bounds in Lemma SM.\ref{lemma_BV}  in the Supplement~: 
\begin{multline*}
 \EE\!\left[\!\left(\sum_{m=1}^{M_0}\{X^{new}(\Tt^{new}_m;\!\boldsymbol B)-X^{new}(\Tt)\} W^{new}_m(\Tt) \right)^2\! \Big{|}M_0\! \right]\\
 \leq 2\left\{\!L_1(\Tt) h_1^{2H_1(\Tt)}\!+L_2(\Tt) h_2^{2H_2(\Tt)}\!\right\}\!\{1+o(1)\},
\end{multline*}
and, for some $a_1>1$,
 $$
 \EE\left[\left(\sum_{m=1}^{M_0}\varepsilon^{new}_m  W^{new}_m(\Tt)  \right)^2 \Big{|}M_0\right]\leq \frac{\kappa^2\sigma^2}{c\pi}\frac{1}{M_0 h_1h_2}\left\{1+a_1 M_0^{-1/4}\right\}.
 $$
\end{proof}

\begin{proof}[Proof of Proposition \ref{risk_2}.]
Let $\widehat \omega (\Tt)$,  $\widehat \alpha_i(\Tt)$ and  $\widehat \Lambda_i (\Tt)$ be the estimators obtained by replacing $H_i(\Tt)$ and $L_i(\Tt)$ in the definitions of     $ \omega (\Tt)$,  $ \alpha_i(\Tt)$ and $ \Lambda_i(\Tt)$, respectively. We define the sets 
		$$
		\mathcal{F}=\bigcap_{i=1,2} \left\{ |\widehat \alpha _i(\Tt)-\alpha_i(\Tt)|\leq \log^{-a}(\mathfrak m)\right\}
\quad\text{		and } \quad 
		\mathcal E= \bigcap_{i=1,2}\left\{ \left|  \widehat \Lambda_i(\Tt)/\Lambda_i(\Tt)-1\right| \leq \log^{-a}(\mathfrak m)\right\},
		$$
with $a$ from assumption \ref{LP4}.

		Since  $\widehat h^*_1$ and $\widehat h^*_2$ are independent of the new realization $X^{new}$,   by  \eqref{risk_1} we obtain 
		\begin{multline*}
\EE\left[\{\widehat X^{new}(\Tt;\widehat{\mathbf{B}}^* )-X^{new}(\Tt)\}^2 \big{|}M_0,\widehat h_1^*,\widehat h_2^*  \right]\leq  \frac{\kappa^2}{c\pi} \frac{\sigma^2}{M_0 \widehat h^*_1 \widehat h^*_2}\\+ 2L_1(\Tt)\{\widehat h_1^*\}^{2
				H_1(\Tt)}+ 2L_2(\Tt)\{\widehat h_2^*\}^{2H_2(\Tt)}.
		\end{multline*}
		Replacing the expressions of $\widehat h^*_1$ and $\widehat h^*_2$, we have 
		$$
					\frac{\kappa^2}{c\pi} \frac{\sigma^2}{M_0 \widehat h_1^* \widehat h_2^*}=  \frac{\kappa^2\sigma^2}{c\pi} M_0^{ \alpha_1(\Tt)+\alpha_2(\Tt)-1}\Lambda_1^{-\alpha_1(\Tt)}(\Tt)\Lambda_2^{-\alpha_2(\Tt)}(\Tt)\times \widehat {\Upsilon}(\Tt),
		$$
	where 	
	\begin{equation*}
\widehat {\Upsilon}(\Tt)=
	M_0^{ \widehat \alpha_1(\Tt)+\widehat\alpha_2(\Tt)-\alpha_1(\Tt)-\alpha_2(\Tt)}\frac{\widehat \Lambda_1^{-\widehat\alpha_1(\Tt)}(\Tt)\widehat \Lambda_2^{-\widehat\alpha_2(\Tt)}(\Tt)}{ \Lambda_1^{-\alpha_1(\Tt)}(\Tt)\Lambda_2^{-\alpha_2(\Tt)}(\Tt)}.
		\end{equation*}
	Let $\EE_{M_0}[\cdot] = \EE[\cdot \mid M_0]$.	Then, on the event $\mathcal F\cap \mathcal E$,  using Cauchy-Schwarz inequality, 
\begin{multline*}
\EE_{M_0\!\!}\left[ \frac{\kappa^2}{c\pi} \frac{\sigma^2}{M_0 \widehat h^*_1 \widehat h^*_2}  \mathbf{1}_{\mathcal{F}\cap\mathcal{E} }  \right]\leq \frac{\kappa^2\sigma^2}{c\pi} M_0^{ \alpha_1(\Tt)+\alpha_2(\Tt)-1}\Lambda_1^{-\alpha_1(\Tt)}(\Tt)\Lambda_2^{-\alpha_2(\Tt)}(\Tt) \EE_{M_0}\left[  \widehat {\Upsilon}(\Tt)\mathbf{1}_{\mathcal{F}\cap\mathcal{E} } \right]
\\ \leq\frac{\kappa^2\sigma^2}{c\pi} M_0^{ \alpha_1(\Tt)+\alpha_2(\Tt)-1+2\log^{-a}(\mathfrak m)}\Lambda_1^{\!-\alpha_1(\Tt)}(\Tt)\Lambda_2^{\!-\alpha_2(\Tt)}(\Tt)\EE_{M_0\!\!} \left[\frac{\widehat \Lambda_1^{-2\widehat\alpha_1(\Tt)}(\Tt)\widehat \Lambda_2^{-2\widehat\alpha_2(\Tt)}(\Tt)}{ \Lambda_1^{-2\alpha_1(\Tt)}(\Tt)\Lambda_2^{-2\alpha_2(\Tt)}(\Tt)}\mathbf{1}_{\mathcal{F}\cap\mathcal{E} }\! \right]\!.
\end{multline*}
	Next, on the event $\mathcal E$, for $i=1,2$, we have
		$$
		\frac{\widehat \Lambda_i^{-2\widehat\alpha_i(\Tt)}(\Tt)}{\Lambda_i^{-2\alpha_i(\Tt)}(\Tt)}= \frac{\widehat \Lambda_i^{-2\widehat\alpha_i(\Tt)}(\Tt)}{\Lambda_i^{-2\widehat\alpha_i(\Tt)}(\Tt)}\frac{\Lambda_i^{-2\widehat\alpha_i(\Tt)}(\Tt)}{\Lambda_i^{-2\alpha_i(\Tt)}(\Tt)}\leq \left(1+\log^{-a}(\mathfrak m)\right)^{-2\widehat \alpha_i(\Tt)}\Lambda_i(\Tt)^{2\log^{-a}(\mathfrak m)}.
		$$
	Note that, by definition, $2\alpha_i(\Tt)<1$, $i=1,2$. 	It follows that on the event $\mathcal F \cap \mathcal E$, we have  
		$$
		\frac{\widehat \Lambda_1^{-2\widehat\alpha_1(\Tt)}(\Tt)\widehat \Lambda_2^{-2\widehat\alpha_2(\Tt)}(\Tt)}{ \Lambda_1^{-2\alpha_1(\Tt)}(\Tt)\Lambda_2^{-2\alpha_2(\Tt)}(\Tt)}= 1+O_\PP(\log^{-a}(\mathfrak m)).
		$$
		Consequently, under $\mathcal F \cap \mathcal E$ we obtain 
		$$
		\EE_{M_0\!\!} \left[ \frac{\kappa^2}{c\pi} \frac{\sigma^2}{M_0 \widehat h^*_1 \widehat h^*_2} \mathbf{1}_{\mathcal{F}\cap\mathcal{E} } \right]\leq \frac{\kappa^2 \sigma^2}{c\pi\Lambda_1^{\alpha_1(\Tt)}(\Tt)\Lambda_2^{\alpha_2(\Tt)}(\Tt)} M_0^{ -\frac{\omega(\Tt)}{2\omega(\Tt)+1}+2\log^{-a}(\mathfrak m)}  \{1+O(\log^{-a}(\mathfrak m))\}.
		$$
		By similar argument, we can also show that on the event $\mathcal F \cap \mathcal E$, we have the bound 
\begin{multline*}
			\EE_{M_0\!\!}\left[ L_1(\Tt) \{\widehat h^*_1\}^{2H_1(\Tt)}  \mathbf{1}_{\mathcal{F}\cap\mathcal{E} }\right] \\ \leq L_1(\Tt)\left[\Lambda_1(\Tt)^{2H_1(\Tt)+1}\big/\Lambda_2(\Tt)\right]^{\!\alpha_1(\Tt)}M_0^{-\frac{\omega(\Tt)}{2\omega(\Tt)+1}+2\log^{-a}(\mathfrak m)}\!\times \{1+O(\log^{-a}(\mathfrak m))\},	 
\end{multline*}
		and symmetrically the bound for $ \EE_{M_0\!\!}\left[ L_2(\Tt) \{\widehat h^*_2\}^{2H_2(\Tt)} \mathbf{1}_{\mathcal{F}\cap\mathcal{E} }\big{|}M_0\right]$.
	Since
		\begin{multline*}
			\EE_{M_0\!\!} \left[ \{\widehat X^{new}(\Tt;\widehat{\mathbf{B}}^* )-X^{new}(\Tt)\}^2 \right]\leq \EE_{M_0\!\!}\left[ \{\widehat X^{new}(\Tt;\widehat{\mathbf{B}}^* )-X^{new}(\Tt)\}^2\mathbf{1}_{\mathcal F}\mathbf{1}_{\mathcal E} \right]\\
			+\EE_{M_0\!\!}\left[ \{\widehat X^{new}(\Tt;\widehat{\mathbf{B}}^* )-X^{new}(\Tt)\}^2\mathbf{1}_{\overline{\mathcal F}} \right]+\EE_{M_0\!\!}\left[ \{\widehat X^{new}(\Tt;\widehat{\mathbf{B}}^* )-X^{new}(\Tt)\}^2\mathbf{1}_{\overline{\mathcal E}} \right],
		\end{multline*}
	and given the facts above, it remains to investigate the last two expectations in the last display. 
		Using \ref{LP4}, \ref{ass_D} and Cauchy-Schwarz inequality,  we get 
		$$
		\EE_{M_0\!\!}\left[ \{\widehat X^{new}(\Tt;\widehat{\mathbf{B}}^* )\!-\!X^{new}(\Tt)\}^2\mathbf{1}_{\overline{\mathcal F}} \right]\!+\EE_{M_0\!\!}\left[ \{\widehat X^{new}(\Tt;\widehat{\mathbf{B}}^* )\! -\! X^{new}(\Tt)\}^2\mathbf{1}_{\overline{\mathcal E}} \right]\! \!= o(\log^{-a}(\mathfrak m)).
		$$ 
		We finally deduce 
		$$
		\EE_{M_0\!\!}\left[ \{\widehat X^{new}(\Tt;\widehat{\mathbf{B}}^* )-X^{new}(\Tt)\}^2 \right]\leq \Gamma_2(\Tt) M_0^{-\frac{\omega(\Tt)}{2\omega(\Tt)+1}+2\log^{-a}(\mathfrak m)}\times \{1+O(\log^{-a}(\mathfrak m))\},
		$$
		with $\Gamma_2(\Tt)$ defined in Proposition \ref{risk_2}.		
		\end{proof}

\medskip

\noindent \textbf{\large Acknowledgements:}
V. Patilea acknowledges support from the grant of the Ministry of Research, Innovation and Digitization, CNCS/CCCDI-UEFISCDI,  number PN-III-P4-ID-PCE-2020-1112, within PNCDI III.

\medskip

\noindent\textbf{\large Supplementary Material:}
In the Supplement we provide complements for the proofs  of Propositions \ref{propCH}, \ref{conc_Lest}, \ref{mprop}, \ref{prop_risk1}, and we prove some technical lemmas.  
Moreover, the justification for the local  Hölder continuity of the  realizations of $X$,  stated in Section \ref{sec6} above, is provided.

\bibliographystyle{apalike}

\bibliography{biblio_final.bib}


\medskip

\noindent \textbf{\large Acknowledgements:}
V. Patilea acknowledges support from the grant of the Ministry of Research, Innovation and Digitization, CNCS/CCCDI-UEFISCDI,  number PN-III-P4-ID-PCE-2020-1112, within PNCDI III.

\medskip

\noindent\textbf{\large Supplementary Material:}
In the Supplement we provide complements for the proofs  of Propositions \ref{propCH}, \ref{conc_Lest}, \ref{mprop}, \ref{prop_risk1}, and we prove some technical lemmas.  
Moreover, the justification for the local  Hölder continuity of the  realizations of $X$,  stated in Section \ref{sec6} above, is provided.


\bibliographystyle{apalike}

\bibliography{biblio_final.bib}

\end{document}


\maketitle
	
	\begin{abstract}
		In this Supplementary Material we provide complementary arguments for the proofs of the propositions in the Appendix of the main manuscript. More precisely, in Sections \ref{secSM_propCH}, \ref{secSM_conc_Lest}, \ref{secSM_mprop}, \ref{secSM_prop_risk1} below we complete the proofs of Propositions \ref{propCH}, \ref{conc_Lest}, \ref{mprop}, \ref{prop_risk1}, respectively.  
		Moreover, in Section \ref{tech_lem} below we provide the justification for the local  Hölder continuity of the  realizations of $X$, as stated in Section \ref{sec6} of the main manuscript, and we prove some technical lemmas. 
	\end{abstract}


\section{Complements for the proof of Proposition \ref{propCH}}\label{secSM_propCH}
\begin{proof}[Proof of \eqref{eq:assumption-theta-hat_main}]
We show that  a constant $\mathfrak u >0$ exists such that, for any $i=1,2$  $\varepsilon \in(0,1)$, and  $0<\Delta\leq \Delta_0$~: 
\begin{equation}\label{eq:assumption-theta-hat_SM}
	\max\left\{
	\PP\left(\widehat{\theta}_{\Tt}^{(i)}(\Delta)-\theta_{\Tt}^{(i)}(\Delta)\geq \varepsilon\right),
	\PP\left(\widehat{\theta}_{\Tt}^{(i)}(\Delta)-\theta_{\Tt}^{(i)}(\Delta)\leq -\varepsilon\right)
	\right\}
	\leq \exp\left( -\mathfrak{u}N\varepsilon^2\varrho(\Delta,\mathfrak m)\right),
\end{equation}
with $\widehat{\theta}_{\Tt}^{(i)}(\Delta)$ defined in \eqref{eq_theta_hat}, and provided that 
$\mathfrak m$ is sufficiently large.  
 We  decompose
 $$
 \widehat{\theta}^{(i)}_{\Tt}(\Delta)-\theta^{(i)}_{\Tt}(\Delta)= \underbrace{\frac{1}{N}\sum_{j=1}^N \Bar{Z}_j }_{\text{stochastic term}}
 +\underbrace{\EE\left[\widehat{\theta}^{(i)}_{\Tt}(\Delta)\right]-\theta^{(i)}_{\Tt}(\Delta)}_{\text{bias term}},\qquad i=1,2,
 $$
where $\Bar{Z}_j= Z_j- \EE[Z_j]$, and 
$$
Z_j= \left(\widetilde X^{(j)}\left(\Tt-\frac{\Delta}{2}e_i\right)-\widetilde X^{(j)}\left(\Tt+\frac{\Delta}{2}e_i\right)\right)^2,  \qquad j\in\{1,\cdots, N\}.
$$ 

\medskip

\noindent \textbf{Bounding the bias term :} Recall that $\xi^{(j)}(\Tt)=\widetilde X^{(j)}(\Tt)- X^{(j)}(\Tt)$.
Using the identities $a^2 - b^2 = (a-b)^2 + 2b(a-b)$ and  $(a-b)^2\leq 2a^2+2b^2$,  and  Cauchy-Schwarz inequality,  we obtain
\begin{multline*}
    \EE\left[\widehat {\theta}^{(i)}_{\Tt}(\Delta)\right]-\theta^{(i)}_{\Tt}(\Delta)\\
    =   \EE\!\left[\! \frac{1}{N}\sum_{j=1}^N \! \left(\!\widetilde X^{(j)}\!\left(\! \Tt-\frac{\Delta}{2}e_i\right) \!- \!\widetilde X^{(j)}\!\left(\Tt+\frac{\Delta}{2}e_i\right)\!\right)^{\!2}\right]\!
     -\EE\!\left[ \left(\! X\!\left(\Tt-\frac{\Delta}{2}e_i\right)\!-X\!\left(\!\Tt+\frac{\Delta}{2}e_i\right)\!\right)^{\!2}\right]
   \\
    \hspace{-5.5cm}  =  \frac{1}{N}\sum_{j=1}^N  \EE\left[\left\{\xi^{(j)}\left(\Tt-\frac{\Delta}{2}e_i\right)-\xi^{(j)}\left(\Tt+\frac{\Delta}{2}e_i\right)\right\}^2\right]\\
    +\frac{2}{N}\sum_{j=1}^N \EE\left[\left\{\xi^{(j)}\left(\Tt-\frac{\Delta}{2}e_i\right)-\xi^{(j)}\left(\Tt+\frac{\Delta}{2}e_i\right)\right\}\right.\\ \hspace{4cm} \times \left.\left\{ X^{(j)}\left(\Tt-\frac{\Delta}{2}e_i\right)- X^{(j)}\left(\Tt+\frac{\Delta}{2}e_i\right)\right\}\right]\\
    \leq  4R_2(\mathfrak m) +4\left( R_2(\mathfrak m)\times \theta^{(i)}_{\Tt}(\Delta)\right)^{1/2}  := \frac{\eta^*_{\Tt}(\Delta)}{2}.
\end{multline*}

\medskip

\noindent \textbf{Bounding the stochastic term~:}
Using the convexity of the function $x\mapsto x^p$, \ref{ass_H1} and \ref{ass_H2} we obtain :
\begin{align*}
    \EE\left[|\Bar{Z}_j|^p\right]
    \leq& 2^p \EE[|Z_j|^p]\\
    =& 2^p \EE\left[\left(\widetilde X^{(j)}\left(\Tt-\frac{\Delta}{2}e_i\right)-\widetilde X^{(j)}\left(\Tt+\frac{\Delta}{2}e_i\right)\right)^{2p}\right]\\
    =& 2^p \EE\left[\left\{X^{(j)}\!\left(\!\Tt-\frac{\Delta}{2}e_i\right)\!- X^{(j)}\!\left(\Tt+\frac{\Delta}{2}e_i\right)\! +\xi^{(j)}\!\left(\!\Tt-\frac{\Delta}{2}e_i\right)-\xi^{(j)}\!\left(\!\Tt+\frac{\Delta}{2}e_i\right)\right\}^{2p}\right] \\
    \leq & \frac{18^p}{3}\EE\left[\left(X^{(j)}\left(\Tt-\frac{\Delta}{2}e_i\right)- X^{(j)}\left(\Tt+\frac{\Delta}{2}e_i\right)\right)^{2p}\right]\\
    &+\frac{18^p}{3}\EE\left[\left|\xi^{(j)}\left(\Tt-\frac{\Delta}{2}e_i\right)\right|^{2p}+\left|\xi^{(j)}\left(\Tt+\frac{\Delta}{2}e_i\right)\right|^{2p}\right]\\
    \leq& \frac{18^p}{3}\left(\EE\left[ \left(X^{(j)}\left(\Tt-\frac{\Delta}{2}e_i\right)- X^{(j)}\left(\Tt+\frac{\Delta}{2}e_i\right)\right)^{2p} \right]+ 2R_{2p}(\mathfrak m) \right)\\
    \leq & \frac{18^p}{3}\left( \frac{p!}{2}\mathfrak{a} \mathfrak{A}^{p-2}\Delta^{2p\underline H (\Tt)}+p!\mathfrak{c} \mathfrak{D}^{p-2} \rho(\mathfrak m)^{2p}\right)\\
    \leq& \mathfrak{g}\frac{p!}{2}\mathbf{\mathfrak{G}}^{p-2}\times \max\{ \Delta^{2p\underline H (\Tt)},\rho(\mathfrak m)^{2p}\}=\mathfrak{g}\frac{p!}{2}\mathbf{\mathfrak{G}}^{p-2}\times\varrho(\Delta,\mathfrak m)^{-p} ,
\end{align*}
for some positive constants $\mathfrak G$ and $\mathfrak g$. The bounds on the moments of $\Bar{Z}_j$ allows to derive exponential bounds for the concentration $\widehat{\theta}^{(i)}_{\Tt}(\Delta)-\theta^{(i)}_{\Tt}(\Delta)$. 
Let $\mathfrak{m}$ sufficiently large 
such that  $0<\eta^*_{\Tt}(\Delta) \leq  1 $, and let $\eta^*_{\Tt}(\Delta)\leq \eta \leq 1$.   By Bernstein's inequality,  we get  
\begin{align*}
    \PP \left[\widehat { \theta} ^{(i)}_{\Tt}(\Delta) -\ \theta^{(i)}_{\Tt}(\Delta)\geq \eta \right]=& \PP \left [\frac{1}{N} \sum_{j=1}^N\Bar{Z}_j+\EE\left[\widehat { \theta}^{(i)}_{\Tt}(\Delta)\right]- \theta^{(i)}_{\Tt}(\Delta)\geq \eta\right]\\
    \leq & \PP\left [ \frac{1}{N} \sum_{j=1}^N\Bar{Z}_j+\frac{\eta^*_{\Tt}(\Delta)}{2}\geq \eta \right]\\
    \leq & \PP\left [ \frac{1}{N} \sum_{j=1}^N\Bar{Z}_j\geq \frac{\eta}{2}\right]\\
    \leq & \exp\left(-\frac{N\eta^2}{8\mathfrak{g}\varrho(\Delta,\mathfrak m)^{-2}+4\mathfrak{G}\varrho(\Delta,\mathfrak m)^{-1}\eta} \right).
\end{align*}
Then   
\begin{equation*}
\PP \left[\hat  \theta^{(i)}_{\Tt}(\Delta) -\ \theta^{(i)}_{\Tt}(\Delta)\geq \eta \right] \leq \exp\left(-\mathfrak{u}N\eta^2\varrho(\Delta,\mathfrak m) \right), 
\end{equation*}
where $\mathfrak{u}=\{8\mathfrak{g}+4\mathfrak{G}\}^{-1}$, because  $\eta,\varrho^{-1}(\Delta,\mathfrak m)\leq 1$. A similar bound can be derived for $\mathbb{P}\left[\widehat{\theta}^{(i)}_{\Tt}(\Delta)-\theta^{(i)}_{\Tt}(\Delta)\leq - \eta  \right] $. 
\end{proof}

\medskip

	Let us notice that, for proving \eqref{eq:concentration-Hhat-around-H} we have to apply \eqref{eq:assumption-theta-hat_main} with 
	$ \max \{|R(\underline H )(\Tt)|,|R(\overline H - \underline H )(\Tt) | \}\leq \varepsilon $. The stronger condition $  \max \{|\log(\Delta)||R(\underline H )(\Tt)|,|R(\overline H - \underline H )(\Tt) | \}\leq \varepsilon $ will be needed later, for deriving \eqref{eq:second-term}. 
	
	\medskip
	
\begin{proof}[Proof of  \eqref{eq:bound-alpha_main}]
 Let $\varepsilon$ such that $ \max \{|\log(\Delta)| |R(\underline H )(\Tt)|, \; |R(\overline H - \underline H )(\Tt) | \}\leq \varepsilon $. Under the conditions of Proposition \ref{propCH},  we have to prove that 
constants $\tilde L_5$, $\tilde L_6$ exist such that, for sufficiently small  $\Delta$, 
\begin{equation}\label{eq:bound-alpha_SM}
	\PP(|\widehat{\alpha}_{\Tt}(\Delta)-\alpha_{\Tt}(\Delta)|\geq \varepsilon)
	\leq \tilde L_5
	\exp\left(
	-\tilde L_6N\varepsilon^2\frac{\Delta^{4\overline{H}(\Tt)}\varrho(\Delta)}{\log^2(\Delta)}
	\right).
\end{equation} 
 Here, $\varrho(\Delta)$ is a short notation for $\varrho(\Delta)=\varrho(\Delta,\mathfrak m)$.  Let us define the event
\begin{equation*}
    \Omega_1 = \left\{ \frac{\widehat{\gamma_{\Tt}}(2\Delta)}{(2\Delta)^{2\widehat{\underline{H}}(\Tt)}} \neq \frac{\widehat{\gamma_{\Tt}}(\Delta)}{(\Delta)^{2\widehat{\underline{H}}(\Tt)}} \right\}.
\end{equation*}
Using the definition of $\widehat\alpha_{\Tt}(\Delta)$, we have
\begin{align}
    \PP(|\widehat{\alpha}_{\Tt}(\Delta)-\alpha_{\Tt}(\Delta)|\geq \varepsilon, \Omega_1)
    &\leq
    \PP\left(\left|
        \left|\frac{\widehat\gamma_{\Tt}(2\Delta)}{(2\Delta)^{2\widehat{\underline{H}}(\Tt)}}-\frac{\widehat\gamma_{\Tt}(\Delta)}{\Delta^{2\widehat{\underline{H}}(\Tt)}}\right|
        -
        \left|\frac{\gamma_{\Tt}(2\Delta)}{(2\Delta)^{2\underline{H}(\Tt)}}-\frac{\gamma_{\Tt}(\Delta)}{\Delta^{2\underline{H}(\Tt)}}\right|
    \right| \geq \varepsilon\right)\nonumber \\
    &\leq
    \PP\left(\left|
        \frac{\widehat\gamma_{\Tt}(2\Delta)}{(2\Delta)^{2\widehat{\underline{H}}(\Tt)}}- \frac{\gamma_{\Tt}(2\Delta)}{(2\Delta)^{2\underline{H}(\Tt)}}\right|
        +
        \left| \frac{\widehat\gamma_{\Tt}(\Delta)}{\Delta^{2\widehat{\underline{H}}(\Tt)}}-\frac{\gamma_{\Tt}(\Delta)}{\Delta^{2\underline{H}(\Tt)}}\right|
    \geq \varepsilon\right)\nonumber\\
    &\leq B_\varepsilon + B_\varepsilon' + B_\varepsilon'' + B_\varepsilon''',
    \label{eq:decomposition-concentration-alpha}
\end{align}
where
\begin{equation*}
    B_\varepsilon 
    = \PP\left(\frac{\widehat{\gamma}_{\Tt}(\Delta)}{\Delta^{2\widehat{\underline{H}}(\Tt)}} -\frac{{\gamma_{\Tt}}(\Delta)}{\Delta^{2\underline{H}(\Tt)}}\geq \frac{\varepsilon}{2}\right),
    \qquad
    B_\varepsilon' 
    = \PP\left(\frac{\widehat{\gamma}_{\Tt}(\Delta)}{\Delta^{2\widehat{\underline{H}}(\Tt)}} -\frac{{\gamma_{\Tt}}(\Delta)}{\Delta^{2\underline{H}(\Tt)}}\leq -\frac{\varepsilon}{2}\right).
\end{equation*}
\begin{equation*}
    B_\varepsilon''
    = \PP\left(\frac{\widehat{\gamma}_{\Tt}(2\Delta)}{(2\Delta)^{2\widehat{\underline{H}}(\Tt)}} -\frac{{\gamma_{\Tt}}(2\Delta)}{(2\Delta)^{2\underline{H}(\Tt)}}\geq \frac{\varepsilon}{2}\right),
    \qquad
    B_\varepsilon''' 
    = \PP\left(\frac{\widehat{\gamma}_{\Tt}(2\Delta)}{(2\Delta)^{2\widehat{\underline{H}}(\Tt)}} -\frac{{\gamma_{\Tt}}(2\Delta)}{(2\Delta)^{2\underline{H}(\Tt)}}\leq -\frac{\varepsilon}{2}\right).
\end{equation*}
Since all these terms can be bounded in a similar way, we  focus on $B_\varepsilon$. By rearranging the different terms and using simple algebra, we get
\begin{align*}
    B_\varepsilon
    &=\PP\left(\frac{{\gamma_{\Tt}}(\Delta)}{\Delta^{2\underline{H}(\Tt)}}\left(\frac{\widehat{\gamma}_{\Tt}(\Delta) \Delta ^{2\underline{H}(\Tt)}}{\gamma_{\Tt}(\Delta)\Delta^{2\widehat{\underline{H}}(\Tt)}} -1\right)\geq \frac{\varepsilon}{2}\right)\\
    &\leq \PP\left(\frac{\widehat{\gamma}_{\Tt}(\Delta)}{\gamma_{\Tt}(\Delta)}\geq\left(1+ \frac{\varepsilon}{2}\frac{\Delta^{2\underline{H}(\Tt)}}{\gamma_{\Tt}(\Delta)}\right)^{\frac{1}{2}}\right)+ \PP\left(\frac{\Delta^{2\widehat{\underline{H}}(\Tt)}}{\Delta^{2\underline{H}(\Tt)}}\leq\left(1+ \frac{\varepsilon}{2}\frac{\Delta^{2\underline{H}(\Tt)}}{\gamma_{\Tt}(\Delta)}\right)^{-\frac{1}{2}}\right)\\
    &\leq \PP\left(\frac{\widehat{\gamma}_{\Tt}(\Delta)}{\gamma_{\Tt}(\Delta)}\geq 1+\delta(\Delta)\frac{\varepsilon}{2}\right)
    + \PP\left(\frac{\Delta^{2\widehat{\underline{H}}(\Tt)}}{\Delta^{2\underline{H}(\Tt)}}\leq1-\delta(\Delta) \frac{\varepsilon}{2}\right),
\end{align*}
where
\begin{equation*}
    \delta^2(\Delta)=
    \frac{\Delta^{4\underline{H}(\Tt)}/\gamma^2_{\Tt}(\Delta)}%
    {4\{1+\Delta^{2\underline{H}(\Tt)}/\gamma_{\Tt}(\Delta)\}} 
    = K_1^*(\Tt)\{1+o(1)\}
    \quad \text{with} \quad
    K_1^*(\Tt)=\frac{1}{4\{1+K_1^2(\Tt)\}}.
\end{equation*}
 Using~\eqref{eq:assumption-theta-hat_SM},   the first term can easily be bounded by
\begin{equation}\label{eq:first-term}
    \PP\left(\frac{\widehat{\gamma}_{\Tt}(\Delta)}{\gamma_{\Tt}(\Delta)}\geq 1+\delta(\Delta)\frac{\varepsilon}{2}\right) 
    \leq  2\exp \left(-\frac{\mathfrak{u}}{4}N\varepsilon^2\delta^2(\Delta)\Delta^{4\underline{H}(\Tt)}\varrho(\Delta)\right).
\end{equation}
The study of the second term boils down to the concentration of $\widehat{\underline{H}}(\Tt)$ around $\underline{H}(\Tt)$. Indeed, using~\eqref{eq:concentration-Hhat-around-H}, and since $\log(1+x) \leq x$, $\forall x>-1$, we have
\begin{align}
    \PP\left(\frac{\Delta^{2\widehat{\underline{H}}(\Tt)}}{\Delta^{2\underline{H}(\Tt)}}\leq1-\delta(\Delta) \frac{\varepsilon}{2}\right)
    &\leq \PP\left(2\widehat{\underline{H}}(\Tt)-2\underline{H}(\Tt)\geq - \frac{\delta(\Delta)}{\log(\Delta)} \frac{\varepsilon}{2}\right)\nonumber\\
    &\leq  C_1\exp \left(-C_2N\varepsilon^2\frac{\delta^2(\Delta)}{16\log^2(\Delta)}\Delta ^{4\underline{H}(\Tt)}\varrho(\Delta)\right)\label{eq:second-term}.
\end{align} 
Therefore, taking together~\eqref{eq:first-term} and~\eqref{eq:second-term},  positive constants $\tilde L_1$ and $\tilde L_2$ exist such that, for $\Delta$ sufficiently small,
\begin{equation}\label{eq:control-Bvarepsilon}
    B_\varepsilon
    \leq 
    \widetilde{L}_1\exp\left(-\widetilde{L}_2N\varepsilon^2\frac{K_1^*(\Tt)}{\log^2(\Delta)}\Delta ^{4\underline{H}(\Tt)}\varrho(\Delta)\right).
\end{equation}
By~\eqref{eq:decomposition-concentration-alpha}, this implies that two positive constants $L_3$ and $L_4$ exist such that, for $\Delta$ sufficiently small,
\begin{equation*}
    \PP(|\widehat{\alpha}_{\Tt}(\Delta)-\alpha_{\Tt}(\Delta)|\geq \varepsilon)
    \leq 
    L_3\exp\left(-L_4N\varepsilon^2\frac{\Delta ^{4\underline{H}(\Tt)}\varrho(\Delta)}{\log^2(\Delta)}\right)
    + \PP(\overline{\Omega}_1).
\end{equation*}
It remains to bound $\PP(\overline{\Omega}_1)$. For this, assume without loss of generality, that 
$$
 {{\gamma_{\Tt}}(2\Delta)}/{(2\Delta)^{2{\underline{H}}(\Tt)}}  \geq {{\gamma_{\Tt}}(\Delta)}/{\Delta^{{2\underline{H}}(\Tt)}}.
$$ 
By definition of $\overline{\Omega}_1$, we then have~: 
\begin{align*}
    \PP(\overline{\Omega}_1) 
    &= \PP\left(
        \frac{\widehat{\gamma}_{\Tt}(2\Delta)}{(2\Delta)^{2\widehat{\underline{H}}(\Tt)}}  - \frac{\widehat{\gamma}_{\Tt}(\Delta)}{(\Delta)^{2\widehat{\underline{H}}(\Tt)}}
        = 0
    \right) \\
    &\leq \PP\left(
        \frac{\widehat{\gamma}_{\Tt}(2\Delta)}{(2\Delta)^{2\widehat{\underline{H}}(\Tt)}}  - \frac{\widehat{\gamma}_{\Tt}(\Delta)}{(\Delta)^{2\widehat{\underline{H}}(\Tt)}}
        \leq (1-\varepsilon) 
        \left(
        \frac{{\gamma_{\Tt}}(2\Delta)}{(2\Delta)^{2{\underline{H}}(\Tt)}}  - \frac{{\gamma_{\Tt}}(\Delta)}{(\Delta)^{2{\underline{H}}(\Tt)}}
        \right)
    \right)\\
    &\leq \PP\left(\frac{\widehat{\gamma}_{\Tt}(2\Delta)}{(2\Delta)^{2\widehat{\underline{H}}(\Tt)}}-\frac{{\gamma_{\Tt}}(2\Delta)}{(2\Delta)^{2\underline{H}(\Tt)}}\leq \frac{\varepsilon}{2}\alpha_{\Tt}(\Delta)\right)+\PP\left(\frac{\widehat{\gamma}_{\Tt}(\Delta)}{\Delta^{2\widehat{\underline{H}}(\Tt)}}-\frac{{\gamma_{\Tt}}(\Delta)}{\Delta^{2\underline{H}(\Tt)}}\geq- \frac{\varepsilon}{2}\alpha_{\Tt}(\Delta)\right)\\ 
    &\leq \widetilde{L}_3\exp \left(-\widetilde{L}_4N\varepsilon^2\alpha_{\Tt}^2(\Delta)\frac{\Delta ^{4\underline{H}(\Tt)}\varrho(\Delta)}{\log^2(\Delta)}\right).
\end{align*}
 For the last inequality,  we use inequalities similar to~\eqref{eq:control-Bvarepsilon} applied with $\Delta$ and $2\delta$, and where $\varepsilon$ was replaced by $\varepsilon\alpha_{\Tt}(\Delta)$.  Now, let us notice that
\begin{equation*}
    \alpha_{\Tt}(\Delta) = \frac{K_2(\Tt)}{2^{2D(\Tt)}-1} \Delta^{2D(\Tt)} (1+o(1)),
\end{equation*}
where $K_2(\Tt)$ is defined by \eqref{eq:K1K2} and, recall that $D(\Tt)= \overline H(\Tt) - \underline H(\Tt)$.
This implies that constants $\tilde L_5$ and $\tilde L_6$ exist such that, for sufficiently small $\Delta$,
\begin{equation}\label{eq:bound-alpha}
    \PP(|\widehat{\alpha}_{\Tt}(\Delta)-\alpha_{\Tt}(\Delta)|\geq \varepsilon)
    \leq \tilde L_5
    \exp\left(
        -\tilde L_6N\varepsilon^2\frac{\Delta^{4\overline{H}(\Tt)}\varrho(\Delta)}{\log^2(\Delta)}
    \right).
\end{equation}
\end{proof}


\section{Complements for the proof of Proposition \ref{conc_Lest}}\label{secSM_conc_Lest}
 
We first bound the probability of the event $\{\widehat D(\Tt)=0\}$ when $D(\Tt) >0$, where
$$
\widehat D(\Tt) =\widehat{(\overline{H}-\underline{H})}(\Tt)=\frac{\log(\widehat{\alpha}_{\Tt}(2\Delta))-\log(\widehat{\alpha}_{\Tt}(\Delta))}{2\log(2)}.
$$
By \eqref{eq:concentration-Difference}, for $  D(\Tt) \geq  |\log(\Delta)| \max \{|R(\underline H )(\Tt)|, \;|R(\overline H - \underline H )(\Tt)  |\}  $, which is trivially satisfied for small $\Delta$,
\begin{multline}\label{bound_Dhat_negli}
\PP\left(\widehat D(\Tt)=0\right) \leq \PP\left(\left|\widehat D(\Tt) -  D(\Tt)\right| \geq  D(\Tt)\right) 
\\ \leq 
\tilde L_8
\exp\left(
-\tilde L_7N\times D^2(\Tt)
\times 
\frac{\Delta^{4\overline{H}(\Tt)}\varrho(\Delta)}{\log^2(\Delta)}\Delta^{4D(\Tt)}
\right).
\end{multline}
The bound in \eqref{bound_Dhat_negli} will be shown to be negligible compared to the bound of 
$$
\PP\left(\left|\widehat{L_2^{(i)}}(\Tt)-L_2^{(i)}(\Tt)\right|\geq \varepsilon, \widehat D(\Tt) \neq 0 \right),
$$
given in \eqref{proba-interm} below.

Next, the following result implies the bounds for the terms $B^{(i)}_3$, $B^{(i)}_4$, $i=1,2$, in \eqref{B3B4_main} in the proof of Proposition \ref{conc_Lest} in the Appendix. 

\begin{lemSM}\label{lemma_L2_use}
Let $\varepsilon\in (0,1)$ and set 
\begin{equation*}
\mathfrak D_{\varepsilon}=\PP\left( \frac{1}{(2^{2\widehat{D}(\Tt)}-1)\Delta^{2\widehat{D}(\Tt)}} -\frac{1}{(2^{2D(\Tt)}-1)\Delta^{2D(\Tt)}}\geq \varepsilon, \widehat D(\Tt) \neq 0\right).
\end{equation*}
Then two positive constants $\tilde C_3$ and $\tilde C_4$ exists such that~:
\begin{equation}\label{proba-interm}
\mathfrak D_{\varepsilon}\leq \tilde C_3 \exp\left(
	-\tilde C_4 N\varepsilon^2\left(2^{2D(\Tt)}-1\right)^2\frac{\Delta^{4\overline H (\Tt)}\varrho(\Delta)}{\log^4(\Delta)}\Delta^{8D(\Tt)}
	\right).
\end{equation}
\end{lemSM}

\begin{proof}[Proof of Lemma SM.\ref{lemma_L2_use}]
	Let introduce the following notation~: for any event $\mathcal A$, 
	$$
	\PP_{\neq} (\mathcal A) = \PP (\mathcal A , \widehat D(\Tt) \neq 0).
	$$
 We first  decompose as follow~:
\begin{multline}
\mathfrak D_{\varepsilon}= \PP_{\neq}\left(\Delta^{2D(\Tt)-2\widehat{D}(\Tt)}\frac{2^{2\widehat{D}(\Tt)}-1}{2^{2D(\Tt)}-1}-1\geq \varepsilon\Delta^{2D(\Tt)}(2^{2D(\Tt)}-1)\right)\\
\leq \PP_{\neq}\left(\Delta^{2D(\Tt)-2\widehat{D}(\Tt)}\geq 1+\varepsilon\Delta^{2D(\Tt)}\frac{2^{2D(\Tt)}-1}{2}\right)\\
+ \PP_{\neq}\left(\frac{2^{2\widehat{D}(\Tt)}-1}{2^{2D(\Tt)}-1}\geq 1+\varepsilon\Delta^{2D(\Tt)}\frac{2^{2D(\Tt)}-1}{2}\right).
\end{multline}
Since $\Delta $ is assumed to be sufficiently small, we first have~:
\begin{equation*}
    \PP_{\neq} \left(\Delta^{2D(\Tt)-2\widehat{D}(\Tt)}\geq 1+\varepsilon\Delta^{2D(\Tt)}\frac{2^{2D(\Tt)}-1}{2}\right)\leq \PP_{\neq}\left(D(\Tt)-\widehat{D}(\Tt)\leq \varepsilon\Delta^{2D(\Tt)}\frac{2^{2D(\Tt)}-1}{4\log(\Delta)}\right).
\end{equation*}
Now using \eqref{eq:concentration-Difference} with small $\Delta$, we obtain for any  $|R(\underline H )(\Tt)| \leq \varepsilon\Delta^{2D(\Tt) }/|\log(\Delta)|^2 $~: 
\begin{multline}\label{T1}
 \PP_{\neq}\left(\Delta^{2D(\Tt)-2\widehat{D}(\Tt)}\geq 1+\varepsilon\Delta^{2D(\Tt)}\frac{2^{2D(\Tt)}-1}{2}\right) \\
 \leq  \PP_{\neq}\left( 2 |D(\Tt)-\widehat{D}(\Tt) |
  \geq \varepsilon \frac{\Delta^{2D(\Tt)}}{|\log(\Delta)|}\frac{2^{2D(\Tt)}-1}{2}\right)
 \\
\leq\tilde L_8 \exp\left(
	-\frac{\tilde L_7}{4} N\varepsilon^2\left(2^{2D(\Tt)}-1\right)^2\frac{\Delta^{4\overline{H}(\Tt)}\varrho(\Delta)}{\log^4(\Delta)}\Delta^{8D(\Tt)}
	\right).
\end{multline}

For the second term, and by using simple algebra we have~:
\begin{equation}
 \PP_{\neq}\left(\frac{2^{2\widehat{D}(\Tt)}-1}{2^{2D(\Tt)}-1}\geq 1+\varepsilon\Delta^{2D(\Tt)}\frac{2^{2D(\Tt)}-1}{2}\right)= \PP_{\neq}\left(2^{2\widehat{D}(\Tt)-2D(\Tt)}\geq \varepsilon\Delta^{2D(\Tt)}\frac{\left(2^{2D(\Tt)}-1\right)^2}{2^{2D(\Tt)+1}}+1\right)
\end{equation}
Taking the logarithm on both sides and using the  inequality $\log(x)\leq x-1$,  $\forall x>0$, we obtain~:
\begin{equation}
 \PP_{\neq}\left(\frac{2^{2\widehat{D}(\Tt)}-1}{2^{2D(\Tt)}-1}\geq 1+\varepsilon\Delta^{2D(\Tt)}\frac{2^{2D(\Tt)}-1}{2}\right)\leq  \PP_{\neq}\left(\widehat{D}(\Tt)-D(\Tt)\geq\varepsilon\Delta^{2D(\Tt)}\frac{\left(2^{2D(\Tt)}-1\right)^2}{2^{2D(\Tt)+2}\log2} \right).
\end{equation}
Finally, using \eqref{eq:concentration-Difference}
\begin{multline}\label{T2}
 \PP_{\neq}\left(\frac{2^{2\widehat{D}(\Tt)}-1}{2^{2D(\Tt)}-1}\geq 1+\varepsilon\Delta^{2D(\Tt)}\frac{2^{2D(\Tt)}-1}{2}\right)\\\leq \tilde L_8\exp\left(-\frac{\tilde L_7}{2^{4D(\Tt)+2}\log^22} N\varepsilon^2\left(2^{2D(\Tt)}-1\right)^4\frac{\Delta^{4\overline H(\Tt)}\varrho(\Delta)}{\log^2(\Delta)}\Delta^{8D(\Tt)}\right).
\end{multline}
Taking together \eqref{T1} and \eqref{T2}, we deduce that  positive constants $\tilde C_3$ and $\tilde C_4$ exist such that, $\forall \varepsilon\in (0,1)$, 
\begin{equation}\label{proba-interm}
\mathfrak D_{\varepsilon}\leq \tilde C_3 \exp\left(
	-\tilde C_4 N\varepsilon^2\left(2^{2D(\Tt)}-1\right)^2\frac{\Delta^{4\overline H (\Tt)}\varrho(\Delta)}{\log^4(\Delta)}\Delta^{8D(\Tt)}
	\right).
\end{equation}
\end{proof}

\medskip

We can now derive the bounds for the terms $B_3$, $B_4$ in
\eqref{B3B4_main}
 from the proof of Proposition \ref{conc_Lest}.
We first notice that  
$$
B_3^{(i)}= \mathfrak D_{\frac{\varepsilon}{3\alpha_{\Tt}^{(i)}(\Delta)}}\quad \text{ and }\quad B_4^{(i)}= \mathfrak D_{\sqrt{\varepsilon/3}}.
$$
Since  $\alpha_{\Tt}^{(i)}(\Delta)\sim L_2^{(i)}(\Tt)(2^{2D(\Tt)}-1)\Delta^{2D(\Tt)}$, we obtain~: with small $\Delta$ and  $|R(\underline H )(\Tt)| \leq \varepsilon\Delta^{2D(\Tt) }/|\log(\Delta)|^2 $,  
\begin{equation}\label{eq:bound-B3}
B_3^{(i)} \leq \tilde C_3 \exp\left(
	-\tilde C_4 N\left(\frac{\varepsilon \Delta^{-2D(\Tt)}}{3L_2^{(i)}(\Tt)(2^{2D(\Tt)}-1)} \right)^2\left(2^{2D(\Tt)}-1\right)^2\frac{\Delta^{4\overline H (\Tt)}\varrho(\Delta)}{\log^4(\Delta)}\Delta^{8D(\Tt)}
	\right), 
\end{equation}
\begin{equation}\label{eq:bound-B4}
B_4^{(i)}\leq  \tilde C_3 \exp\left(
	-\tilde C_4 N\varepsilon\left(2^{2D(\Tt)}-1\right)^2\frac{\Delta^{4\overline H (\Tt)}\varrho(\Delta)}{3\log^4(\Delta)}\Delta^{8D(\Tt)}
	\right).   
\end{equation}
Finally, by combining \eqref{eq:bound-B1}, 
 \eqref{eq:bound-B3} and \eqref{eq:bound-B4},  
 two positive constants $C_3$ and $C_4$ exists such that for any $\varepsilon\in (0,1)$ we have 
\begin{equation}
\PP\left(\left|\widehat{L_2^{(i)}}(\Tt)-L_2^{(i)}(\Tt)\right|\geq \varepsilon\right)
 \leq C_3\exp\left(\!
	-C_4N\varepsilon\min\{\varepsilon,\Delta^{4D(\Tt)}\}(4^{D(\Tt)}\!-1)^2\frac{\Delta^{4\overline H (\Tt)}\varrho(\Delta)}{\log^4(\Delta)}\Delta^{4D(\Tt)}\!
	\right)\!,
\end{equation} 
$i=1,2$. Here, we also use the fact that $L_2^{(1)}, L_2^{(2)}$ are a bounded function.


\section{Complements for the proof of Proposition \ref{mprop}}\label{secSM_mprop}

\subsection{Proof of equation (\ref{eq:B-diag})}

We have to prove  
\begin{equation}\label{eq:B-dia_SM}
	B(\Tt,\Ss)= \frac{1}{2}+O(\|\Tt-\Ss\|^2).
\end{equation}
Here  $B(\Tt,\Ss) =2D(H_1(\Tt),H_1(\Ss))D(H_2(\Tt),H_2(\Ss))$.
Let us note that the map $(x,y)\mapsto D(x,y)$ admits partial derivatives of any order on $(0,1)\times (0,1)$. Next, let
$$
g(x) = \log(\Gamma (2x+1)) + \log( \sin(\pi x)) =:g_1(x) - g_2(x).
$$
We notice that $g^{\prime\prime}(x) < 0$, for any $x\in (0,1)$. Indeed, using the expression of the derivative of the digamma function, cf. \citet[page 260]{AS1964}, we have
$$
g^{\prime\prime}(x) = 4\sum_{k\geq 0} \frac1{(2x+1+k)^2} - \frac{\pi^2}{\sin^2 (\pi x)}
= g_1^{\prime\prime}(x) - g_2^{\prime\prime}(x).
$$
We deduce that $g^{\prime\prime}$ is decreasing on $[1/2, 1)$ and, since $g^{\prime\prime}(0+) = -\infty$, the function $g^{\prime\prime}_1$ is decreasing on $(0,1/2]$ with
$$
g^{\prime\prime}_1(0) =\pi^2/6, \;\; g^{\prime\prime}_1(1/4) =4 \{\pi^2/2 -4\}, \;\; g^{\prime\prime}_1(1/2) = 4 \{\pi^2/6 - 1\} ,
$$
and the function $g^{\prime\prime}_2$ is decreasing on $(0,1/2]$ with
$$
g_2^{\prime\prime}(0+) =\infty, \;\; g_2^{\prime\prime}(1/4) =2\pi^2, \;\; g^{\prime\prime}_2(1/2) = \pi^2 ,
$$
we conclude that $g^{\prime\prime}<0$ on $(0,1]$. In other words, $x\mapsto g(x)$ is log-concave, and thus
$$
2D(x,y)=\frac{\sqrt{\exp(g(x))\times \exp(g(y))}} {\exp(g((x+y)/2))} <1,\qquad \forall 0<x\neq y \leq 1.
$$
We deduce that $1/2$ is the maximum value that the function  $B$ could reach. Which means that the gradient of $B$ on the diagonal is zero and therefore, we have \eqref{eq:B-dia_SM}.


\subsection{Proof of equation (\ref{eq:double-a})}
We here  prove that
\begin{equation}\label{eq:double-a_SM}
	a(\Tt,\Ss)a(\Ss,\Tt)-1 =O(\|\Tt-\Ss\|^2).
\end{equation}
For two $2-$dimensional vectors $\Uu$ and $\Vv$, let $\Uu \cdot \Vv$ denote their scalar product. Moreover, let $\nabla H_1$ and $\nabla H_2$ be the $2-$dimensional vectors of partial derivatives of $H_1$ and $H_2$, respectively. Since the functions $H_1$ and $H_2$ are supposed to be continuously differentiable, using \eqref{eq:B-diag}, we get~:
\begin{equation*}
a(\Tt,\Ss)=\frac{\frac{1}{2}+\left\{\log(|A_1(\Tt)|)\nabla H_1(\Tt)+\log(|A_2(\Tt)|)\nabla H_2(\Tt)\right\}\cdot (\Tt-\Ss)+O(\|\Tt-\Ss\|^2)}{B(\Tt,\Ss)}.
\end{equation*}  
Similar arguments lead to~:
\begin{equation*}
a(\Ss,\Tt)=\frac{\frac{1}{2}-\left\{\log(|A_1(\Ss)|)\nabla H_1(\Tt)+\log(|A_2(\Ss)|)\nabla H_2(\Tt)\right\} \cdot (\Tt-\Ss)+O(\|\Tt-\Ss\|^2)}{B(\Tt,\Ss)}.
\end{equation*}
Thus, by multiplying  the two quantities, we obtain when $\Ss$ is close  to $\Tt$~:
\begin{multline}\label{eq:petit-a}
a(\Tt,\Ss)a(\Ss,\Tt)
=\frac{1}{B^2(\Tt,\Ss)}\\
\times \left(\frac{1}{4}+\frac{1}{2} \left\{\log\left(\frac{|A_1(\Tt)|}{|A_1(\Ss)|}\right)\nabla H_1(\Tt)+\log \left(\frac{|A_2(\Tt)|}{|A_2(\Ss)|}\right)\nabla H_2(\Tt)\right\} \cdot (\Tt-\Ss)+O(\|\Tt-\Ss\|^2)\right).
\end{multline}
On the another hand we have~:
$$
\log\left(\frac{|A_1(\Tt)|}{|A_1(\Ss)|}\right)=\frac{\nabla A_1(\Tt)}{|A_1(\Tt)|} \cdot (\Tt-\Ss)+O(\|\Tt-\Ss\|^2),$$
 means that 
\begin{equation}\label{inq:firstone}
\log\left(\frac{|A_1(\Tt)|}{|A_1(\Ss)|}\right)\nabla H_1(\Tt) \cdot (\Tt-\Ss)=O(\|\Tt-\Ss\|^2).
\end{equation}
Similar arguments yield~:
\begin{equation}\label{inq:secondeone}
\log\left(\frac{|A_2(\Tt)|}{|A_2(\Ss)|}\right)\nabla H_2(\Tt)\cdot (\Tt-\Ss)=O(\|\Tt-\Ss\|^2).
\end{equation}
By plugging \eqref{inq:firstone} and \eqref{inq:secondeone} into \eqref{eq:petit-a}, and using \eqref{eq:B-dia_SM}, we obtain 
\eqref{eq:double-a_SM} and thus \eqref{eq:double-a}.


\subsection{Proof of equation (\ref{deux_int})}
For any $a,b\in \RR$ we have 
\begin{equation*}
|e^a-e^b|\leq e^ae^{|a-b|}|a-b|.
\end{equation*}
We therfore obtain
\begin{multline*}
		\left|A_1(\Tt)-\widehat A_1(\Tt)\right| \leq A_1(\Tt) \exp\left[\int_{t_0}^{t_1} \left|\widehat f_1(s,t_2)-f_1(s,t_2)\right|{\D} s+\int_{s_0}^{t_2}\left|\widehat g_1(t_0,s)-g_1(t_0,s)\right|\D s\right]\\
		\times \left\{\int_{t_0}^{t_1} \left |\widehat f_1(s,t_2)-f_1(s,t_2)\right|{\D} s+\int_{s_0}^{t_2}\left|\widehat g_1(t_0,s)-g_1(t_0,s)\right|\D s\right\}.
\end{multline*}
It finally remains to show that 
$$
\EE \exp\left[\int_{t_0}^{t_1} \left|\widehat f_1(s,t_2)-f_1(s,t_2)\right|{\D} s+\int_{s_0}^{t_2}\left|\widehat g_1(t_0,s)-g_1(t_0,s)\right|\D s\right]<\infty.
$$
By the Taylor expansion of the exponential function, and using the convexity of the function $x\mapsto x^p$, $p\in \NN$, we obtain
\begin{multline*}
\EE \exp\left[\int_{t_0}^{t_1} \left|\widehat f_1(s,t_2)-f_1(s,t_2)\right|{\D} s+\int_{s_0}^{t_2}\left|\widehat g_1(t_0,s)-g_1(t_0,s)\right|\D s\right]\\
\leq \EE\sum_{p\in\NN}\frac{2^{p-1}}{p!} \left\{  \left[\int_{t_0}^{t_1} \left|\widehat f_1(s,t_2)-f_1(s,t_2)\right|{\D} s\right]^p + \left[\int_{s_0}^{t_2}\left|\widehat g_1(t_0,s)-g_1(t_0,s)\right|\D s\right]^p \right\}.
\end{multline*}
By Jensen's inequality we now obtain 
 \begin{equation}\label{Lpf}
 \left[\int_{t_0}^{t_1} \left|\widehat f_1(s,t_2)-f_1(s,t_2)\right|{\D} s\right]^p\leq(2 |t_1-t_0|)^{p-1}\int_{t_0}^{t_1} \left\{\left|\widehat f_1(s,t_2)\right|^p+\left|f_1(s,t_2)\right|^p\right\}{\D} s.
 \end{equation}
 We also have by \eqref{simpl_v} and \eqref{simpl_A} that
  \begin{equation}
	 f_1(\Tt)\leq  \left(\frac{\overline \beta}{\underline v}\right)^{\frac{1}{2\underline \beta}}, \quad \text{ and }\quad \EE\left[\left|\widehat f_1(\Tt)\right|^p \right]\leq  \left(\frac{\overline \beta}{\underline v}\right)^{\frac{p}{2\underline \beta}}\EE\left[\left(\frac{v(\Tt)}{\widehat v(\Tt)}\right)^{\frac{p}{2\beta}}\right]\leq \left(\frac{\overline \beta}{\underline v}\right)^{\frac{1}{2\underline \beta}}C_v^{p\lceil 1/(2 \underline\beta)\rceil},
 \end{equation} 
 where $\lceil 1/(2\underline \beta)\rceil$ denote the ceiling of $1/(2\underline\beta)$. Therefore, \eqref{Lpf} becomes
 
 \begin{equation}
 \left[\int_{t_0}^{t_1} \left|\widehat f_1(s,t_2)-f_1(s,t_2)\right|{\D} s\right]^p\leq(2 |t_1-t_0| )^{p}\left(\frac{\overline \beta}{\underline v}\right)^{\frac{p}{2\underline \beta}}\times\left\{ 1+C_v^{\lceil 1/(2\underline\beta)\rceil}\right\}^p.
 \end{equation}
 A similar  inequality remains valid for the integral of $|\widehat g_1(t_0,s)-g_1(t_0,s)|$.  We therefore  get 
  \begin{multline}
 \EE \exp\left[\int_{t_0}^{t_1} \left|\widehat f_1(s,t_2)-f_1(s,t_2)\right|{\D} s+\int_{s_0}^{t_2}\left|\widehat g_1(t_0,s)-g_1(t_0,s)\right|\D s\right]\\
\leq \exp\left( 4\operatorname{diam}(\cT)\left(\frac{\overline \beta}{\underline v}\right)^{\frac{1}{2\underline \beta}}\times\left\{ 1+C_v^{\lceil	1/2\underline\beta \rceil }\right\}\right) := \mathfrak C_v.
 \end{multline}


\section{Proof of Proposition \ref{prop_risk1}}\label{secSM_prop_risk1}

The proof of  Proposition \ref{prop_risk1} is a direct consequence of the following lemmas.

\begin{lemSM}\label{Tsy}
Assumptions  \ref{LP1}, \ref{LP2} and \ref{LP3e} hold true.
Then  a positive constant $a_1>1$ exists such that  
$$
\EE\left[\underset{1\leq m\leq M_0}{\max} W^{new}_m(\Tt)\Big{|}M_0\right]\leq \frac{\kappa^2}{ c\pi M_0h_1h_2}\left\{1+a_1M_0^{-1/4}\right\}.
$$
\end{lemSM}

\begin{proof}[Proof of Lemma SM.\ref{Tsy}]
  Recall that, for $1\leq m\leq M_0$ and $\Tt\in \cT$,     
$$
W_m^{new}(\Tt)= \frac{K\left(\boldsymbol{B}(\Tt^{new}_m-\Tt)\right)}{\sum_{m=1}^{M_0}K
	\left(\boldsymbol{B}(\Tt^{new}_m-\Tt)\right)}, \quad \text{with}\quad \frac{0}{0}=0.
$$
Let
$$
N_{h_1,h_2}(\Tt) =  \operatorname{Card} \left\{ \Tt^{new}_m, 1\leq m \leq M_0, \|\boldsymbol{B}(\Tt^{new}_m-\Tt)\| \leq 1\right\}.
$$
Given $M_0$, $N_{h_1,h_2}(\Tt)$ can also be written  as 
$$
N_{h_1,h_2}(\Tt)=B_1(\Tt)+B_2(\Tt)+\dots +B_{M_0}(\Tt) \leq M_0,
$$
where $\left(B_m(\Tt)\right)_{1\leq m\leq M_0} $ are iid Bernoulli random variables with  parameter $p(\Tt)$ given by $$
p(\Tt)=p(\Tt;\boldsymbol{B}) = \PP\left(\boldsymbol{B}(\Tt^{new}_m-\Tt)\in B(0,1)\big{|}M_0\right).
$$ 
Using \ref{LP1} we have
$$
W_m^{new}(\Tt)\leq \kappa^2\frac{\mathbf{1}_{\{N_{h_1,h_2}(\Tt)\geq 1\}}}{N_{h_1,h_2}(\Tt)},\quad \forall 1\leq m\leq M_0,
$$
which implies 
\begin{equation}
\EE\left[\max_m W^{new}_m(\Tt)\Big{|}M_0\right]
=\EE\left[ \max_m W^{new}_m(\Tt)\mathbf{1}_{\{N_{h_1,h_2}(\Tt)\geq 1\}}\Big{|}M_0\right]
\leq\kappa^2 \EE\left[\frac{\mathbf{1}_{\{N_{h_1,h_2}(\Tt)\geq 1\}}}{N_{h_1,h_2}(\Tt)}\Big{|}M_0\right].
\end{equation}
Let
$$
\mathcal G:=\left\{p(\Tt)M_0\frac{\mathbf{1}_{\{N_{h_1,h_2}(\Tt)\geq 1\}}}{N_{h_1,h_2}(\Tt)}\leq 1+M_0^{-1/4} \right\}.
$$
Since \ref{LP2} holds true, by Lemma SM.\ref{Ber_par} we have that  $p(\Tt)\geq \pi ch_1h_2$, and therefore~:
\begin{align*}
\EE\left[\frac{\kappa^2\mathbf{1}_{\{N_{h_1,h_2}(\Tt)\geq 1\} }}{N_{h_1,h_2}(\Tt)}\Big{|}M_0\right]\leq&\frac{\kappa^2}{ c\pi M_0h_1h_2} \EE\left[\frac{M_0p(\Tt)}{N_{h_1,h_2}(\Tt)}
\mathbf{1}_{\{N_{h_1,h_2}(\Tt)\geq 1\}}\Big{|}M_0\right]\\
\leq & \frac{\kappa^2}{ c\pi M_0h_1h_2}\left(1+M_0^{-1/4}\right)+\EE\left[\frac{M_0p(\Tt)}{N_{h_1,h_2}(\Tt)}
\mathbf{1}_{\{N_{h_1,h_2}(\Tt)\geq 1\}}\mathbf{1}_{\overline{ \mathcal G}}\Big{|}M_0\right]\\
\leq& \frac{\kappa^2}{ c\pi M_0h_1h_2}\left(1+M_0^{- 1/4}\right)+M_0 p(\Tt)\PP\left(\overline{ \mathcal G} \mid M_0\right)
\end{align*}
It remains to bound $\PP\left(\overline{ \mathcal G} \mid M_0\right)$. By Hoeffding’s inequality we obtain
\begin{multline}
\PP\left(\overline{ \mathcal G} \mid M_0\right)
= \PP\left(\frac{p(\Tt)M_0}{N_{h_1,h_2}(\Tt)}\geq 1+M_0^{-1/4} \Big{|}M_0\right)\\
= \PP\left( \left\{-N_{h_1,h_2}(\Tt)+ p(\Tt)M_0\right\}(1+M_0^{-1/4})\geq p(\Tt)M_0^{3/4}\right)\\
\leq  \PP\left( -N_{h_1,h_2}(\Tt)+ p(\Tt)M_0\geq \frac{p(\Tt)M_0^{3/4}}{1+M_0^{-1/4}}\right)\\
\leq\exp\left(-2\left(\frac{p(\Tt)M_0^{3/4}}{1+M_0^{-1/4}}\right)^2M_0 \right)
\leq \exp\left(-\frac{\left\{M^{5/4}_0p(\Tt)\right\}^2}{2}\right).
\end{multline}
Since $\mathfrak m h_1h_2$ tends to infinity, there exists $a_1 >1$ such that   
$$
\EE\left[\underset{1\leq m\leq M_0}{\max} W^{new}_m(\Tt)\Big{|}M_0\right]\leq \frac{\kappa^2}{ c\pi M_0h_1h_2}\left(1+a_1M_0^{-1/4}\right).
$$
\end{proof}

\medskip

\begin{lemSM}\label{lemma_BV}
Under the assumptions of Lemma SM.\ref{Tsy}, a positive constant $a_1>1$ exists  such that 
\begin{multline}
    \EE\left[\left(\sum_{m=1}^{M_0}\left\{X^{new}(\Tt^{new}_m)\!-X^{new}(\Tt) \right\} W^{new}_m(\Tt)\right)^2 \Big{|}M_0 \right] \\ \leq 2\left\{L_1(\Tt) h_1^{2H_1(\Tt)}+L_2(\Tt) h_2^{2H_2(\Tt)}\right\}\{1+o(1)\},
\end{multline}
and
$$
\EE\left[\left(\sum_{m=1}^{M_0}\varepsilon^{new}_m W^{new}_m(\Tt)\right)^2 \Big{|}M_0\right]\leq \frac{\kappa^2\sigma^2}{c\pi}\frac{1}{M_0 h_1h_2}\left\{1+ a_1M_0^{-1/4}\right\}.
$$
\end{lemSM}

\begin{proof}[Proof of Lemma SM.\ref{lemma_BV}]
We start by proving the second inequality. For $1\leq m\neq m^\prime\leq M_0$,  $\varepsilon_m^{new}$ and $\varepsilon_{m^\prime}^{new}$ are independent. Moreover, the $\varepsilon_m^{new} $ and $W_m^{new} (\Tt)$ are independent, and the $\varepsilon_m^{new}$ are centered. From these, the fact that  
$\sum_{m=1}^{M_0}W_m^{new}(\Tt)=1 $, and Lemma SM.\ref{Tsy}, we deduce that 
\begin{multline}
\EE\left[\left(\sum_{m=1}^{M_0}\varepsilon_m^{new} W_m^{new} (\Tt) \right)^2\Big{|}M_0 \right]=\EE\left[\sum_{m=1}^{M_0}\{\varepsilon_m^{new} W_m^{new} \}^2(\Tt) \Big{|}M_0 \right]= \sigma^2 \EE\left[\sum_{m=1}^{M_0}\{W^{new}_m\}^2(\Tt) \Big{|}M_0 \right]\\
\leq  \sigma^2 \EE\left[\max_{1\leq m \leq M_0}W_m^{new} (\Tt)  \times \sum_{m=1}^{M_0}W_m^{new} (\Tt) \Big{|}M_0 \right]\leq \frac{\sigma^2 \kappa^2}{ c\pi M_0h_1h_2}\left\{1+a_1M_0^{-1/4}\right\}.
\end{multline}

For the first inequality, denote $\Tt^{new}_m=(t^{new}_{m,1},t^{new}_{m,2})$ and  $\mathcal T ^{new}_{obs}=\left\{\Tt^{new}_m, 1\leq m\leq M_0 \right\}$.
Let $\Tt=(t_1,t_2)\in \cT$. Recall that $K$ is a density with the support in $[-1,1]^2$, and thus $W_m^{new}(\Tt)=0$ as soon as  $|t_1-t_{m,1}^{new}|\geq h_1$ or $|t_2-t_{m,2}^{new}|\geq h_2$. By Jensen's inequality,
\begin{multline}
 \EE\left[\left\{\sum_{m=1}^{M_0}[X^{new}(\Tt^{new}_m)\!-\!X^{new}(\Tt)] W^{new}_m(\Tt)\right\}^{\!2} \Big{|}M_0 \right]\\ \leq  \EE\left[\sum_{m=1}^{M_0}\left\{X^{new}(\Tt^{new}_m)-X^{new}(\Tt)\right\}^2 W^{new}_m(\Tt) \Big{|}M_0 \right].
\end{multline}
Next, for $1\leq m\leq M_0$, we have
\begin{multline*}
	\EE\left[\left\{X^{new}(\Tt^{new}_m)-X^{new}(\Tt)\right\}^2\big{|}M_0,\mathcal T ^{new}_{obs}\right]\\
	\leq
	2\EE\left[\left\{X^{new}(\Tt^{new}_m)-X^{new}(t_1,t^{new}_{m,2}) \right\}^2 \big{|} M_0,\mathcal T ^{new}_{obs} \right]\\
	+2\EE\left[\left\{X^{new}(t_1,t^{new}_{m,2})-X^{new}(\Tt)\right\}^2\big{|}M_0,\mathcal T ^{new}_{obs}\right] .
\end{multline*}
Let $\rho_i=:|t_i-t^{new}_{m,i}|$,  $i=1,2$.
Then by Lemma SM.\ref{lemme_tec_1}, we have 
\begin{multline*}
	\EE\left[\left\{ X^{new}(\Tt^{new}_m)-X^{new}(t_1,t^{new}_{m,2}) \right\}^2 \big{|} M_0,\mathcal T ^{new}_{obs},\rho_1 \leq h_1\right]\leq  L_1(\Tt^{new}_m) \rho_1^{2H_1(\Tt)} \mathbf 1_{\{\rho_1 \leq h_1\}}\\
	\leq L_1(\Tt)\rho_1^{2H_1(\Tt)}\mathbf 1_{\{\rho_1 \leq h_1\}}+|L_1(\Tt)-L_1(\Tt^{new}_m)|\rho_1^{2H_1(\Tt)}\mathbf 1_{\{\rho_1 \leq h_1\}}\\
	\leq L_1(\Tt) h_1^{2H_1(\Tt)}+ k_1h_1^{2H_1(\Tt)}\|\Tt-\Tt_m^{new}\|.
\end{multline*}
For the last inequality we used the fact that  $L_1$ is Lipschitz continuous with Lipschitz constant $k_1$. Using the same type of arguments, and the fact that $L_2$ is Lipschitz continuous with Lipschitz constant $k_2$, we get 
$$
\EE\left[\left\{X^{new}(t_1,t^{new}_{m,2})-X^{new}(\Tt)\right\}^2\big{|}M_0,\mathcal T ^{new}_{obs},\rho_2 \leq h_2\right] \leq   L_2(\Tt) h_2 ^{2H_2(\Tt)}+ k_2h_2^{2H_2(\Tt)}\|\Tt-\Tt_m^{new}\|.
$$
Gathering facts,
$$
\EE\left[\sum_{m=1}^{M_0}\left\{X^{new}(\Tt^{new}_m)\!-\!X^{new}(\Tt)\right\}^2 W_m^{new}(\Tt) \Big{|}M_0 \right]\!\leq 2\left\{L_1(\Tt)h_1^{2H_1(\Tt)}+L_2(\Tt)h_2^{2H_2(\Tt)}\right\}\{1+o(1)\}.
$$
\end{proof}

\section{Technical lemmas}\label{tech_lem}


\subsection{Proof of the local  Hölder continuity of the  realizations of $X$}

In Section \ref{sec6} of the main manuscript, we linked the  local regularity in quadratic mean to the (analytic) regularity of the  realizations of $X$. Let us provide a detailed justification for our statement.

\begin{lemSM}\label{reg_RY_SM}
	Consider $\mathcal{H}^{H_1,H_2}$ in Definition \ref{def} is built with restricted $\boldsymbol L = (L_1,0,0,L_2)$, and let  $X\in \mathcal H^{H_1,H_2}$ 
	Assume that 
	\begin{equation}\label{Yor_us_SM}
		\max_{i=1,2}\sup_{0<\Delta \leq \Delta_0}  \frac{\EE\left[\left\{X\left(\Tt-\Delta e_i/2\right)-X\left(\Tt+\Delta e_i/2\right)\right\}^{2p}\right]}{\EE\left[\left\{X\left(\Tt-\Delta e_i/2\right)-X\left(\Tt+ \Delta e_i/2\right)\right\}^2\right]^p}<\infty, \qquad \forall  p\in \mathbb N.
	\end{equation}
	The, almost  any realization of $X$  is locally $\alpha-$Hölder continuous in the direction $e_i$, for any order $0\leq \alpha < H_i(\Tt)$, $i=1,2$.
\end{lemSM}

\begin{proof}[Proof of Lemma SM.\ref{reg_RY_SM}]
	Condition \eqref{Yor_us_SM} means that, for any $p\in \mathbb N$, there exists $C_p(\Tt)$, independent of $\Delta$, such that
	$$
	\EE\left[\left\{X\left(\Tt-\frac{\Delta}{2}e_i\right)-X\left(\Tt+\frac{\Delta}{2}e_i\right)\right\}^{2p}\right]\leq C_p(\Tt)\EE\left[\left\{X\left(\Tt-\frac{\Delta}{2}e_i\right)-X\left(\Tt+\frac{\Delta}{2}e_i\right)\right\}^2\right]^p.
	$$
	Since   $\boldsymbol L = (L_1,0,0,L_2),$  we obtain
	$$
	\EE\left[\left\{X\left(\Tt-\frac{\Delta}{2}e_i\right)-X\left(\Tt+\frac{\Delta}{2}e_i\right)\right\}^{2p}\right]\leq C_p(\Tt)L_i^p(\Tt) \Delta^{2pH_i(\Tt)}, \quad \forall \Delta\in [0,\Delta_0].  $$
	Using \citet[Theorem 2.1, page 26]{Yor}, we then get that 
	$$
	\EE\!\left[ \underset{0<\Delta\leq \Delta_0}{\sup}\!\left(\frac{|X\left(\Tt-\Delta e_i/2 \right)-X\left(\Tt+ \Delta e_i/2 \right)|}{\Delta^{\nu}}\right)^{2p}\right]< \infty,\quad \forall p\in \mathbb N,  \; \forall \nu \in \left[0, \frac{2pH_i(\Tt)-2}{2p}\right).
	$$
	Letting $p\rightarrow \infty$, we obtain the stated local Hölder property. \end{proof}


\begin{lemSM}\label{lem}
Let $a_N,b_N, c_N$ be three positives sequences such that 
$$
c_N\geq 1 \quad \text{ and  } \quad \max\{a/b, a_N/b_N\}\leq \mathfrak C,
$$ 
for some constant $\mathfrak C >1$.
Let $a,b>0$ and $c\geq 1$. For $0 \leq u\leq 1$, let 
$$
\varphi_1(u)=\PP\left(|a_N-a|\geq u\right), \quad 
\varphi_2(u)=\PP\left(|b_N-b|\geq u\right)
\quad \text{ and } \quad  \varphi_3(u)=\PP\left(|c_N-c|\geq u\right).
 $$
Define
$$
r_1=\frac{a\mathfrak C}{2  \min\left\{1, 1/c\right\}} \quad ,\quad r_2= \frac{b}{6\mathfrak C} \min\left\{1,\frac{1}{c} \right\}  \quad \text{and} \quad r_3= \min\left\{\frac{1}{\sqrt{2\mathfrak C}}, \frac{1}{\left(a/b\right)^c|\log\left(a/b\right)|}\right\}.
$$
We then have 
$$\PP\left(\left|\left(\frac{a_N}{b_N}\right)^{c_N}-\left(\frac{a}{b}\right)^{c}\right|\geq u\right)\leq 2\varphi_1\left(r_1u\right)+2\varphi_2\left( r_2u\right)+2\varphi_3\left(r_3 u\right).
$$

\end{lemSM}

\begin{proof}[Proof of Lemma SM.\ref{lem}]
	By elementary algebra, and using the conditions $c_N,c\geq 1$, and $0 \leq u\leq 1$, 
\begin{multline}
    \PP\left(\left|\left(\frac{a_N}{b_N}\right)^{c_N}-\left(\frac{a}{b}\right)^{c}\right|\geq 2u\right)
    \leq  \PP\left(\left|\left(\frac{a_N}{b_N}\right)^{c_N}-\left(\frac{a}{b}\right)^{c_N}\right|\geq u\right)+\PP\left(\left|\left(\frac{a}{b}\right)^{c_N}-\left(\frac{a}{b}\right)^{c}\right|\geq u\right)\\
    = \PP\left(\left|\int_{\frac{a}{b}}^{\frac{a_N}{b_N}}c_Nx^{c_N-1}{\rm{d}}x\right|\geq u\right)+\PP\left(\left(\frac{a}{b}\right)^{c}\left|\left(\frac{a}{b}\right)^{c_N-c}-1\right|\geq u\right)\\
    \leq \PP\left(c_N\left|\frac{a}{b}-\frac{a_N}{b_N}\right|\mathfrak C\geq u\right)+\PP\left(\left(\frac{a}{b}\right)^{c}\left|\left(\frac{a}{b}\right)^{c_N-c}-1\right|\geq u\right)\\
    \leq \PP\left((c_N-c)\left|\frac{a}{b}-\frac{a_N}{b_N}\right|\geq u/(2 \mathfrak C)\right)+\PP\left(c\left|\frac{a}{b}-\frac{a_N}{b_N}\right|\geq u/(2 \mathfrak C)\right)\\ \hspace{6cm} + \PP\left(\left|\left(\frac{a}{b}\right)^{c_N-c}-1\right|\geq u\left(\frac{a}{b}\right)^{-c}\right)\\
    \leq\varphi_3(\sqrt{u/(2\mathfrak C)})+\PP\left(\left|\frac{a}{b}-\frac{a_N}{b_N}\right|\geq \sqrt{u/(2 \mathfrak C)}\right)+\PP\left(c\left|\frac{a}{b}-\frac{a_N}{b_N}\right|\geq u/(2\mathfrak C)\right)\\\hspace{6cm} + \PP\left(\left|\left(\frac{a}{b}\right)^{c_N-c}-1\right|\geq u\left(\frac{a}{b}\right)^{-c}\right)\\
    \leq \varphi_3(\sqrt{u/(2 \mathfrak C)})+\PP\left(\left|\frac{a}{b}-\frac{a_N}{b_N}\right|\geq \min\left\{1,\frac{1}{c}\right\}u/(2\mathfrak C)\right)+ \PP\left(\left|\left(\frac{a}{b}\right)^{c_N-c}-1\right|\geq u\left(\frac{a}{b}\right)^{-c}\right).
\end{multline}
Moreover, using the inequality 
$$
|xy-1|\leq |x-1||y-1| + |x-1|+|y-1|,
$$
for $\varepsilon>0$, we have

\begin{multline}
\PP\left(\left|\frac{a}{a_N}\times \frac{b_N}{b}-1\right|\geq  \varepsilon \right)\\
\leq \PP\left(\left| \frac{a}{a_N}-1\right|\times\left| \frac{b_N}{b}-1\right|\geq\varepsilon/3 \right)+\PP\left(\left| \frac{a}{a_N}-1\right|\geq\varepsilon/3 \right)+\PP\left(\left| \frac{b_N}{b}-1\right|\geq\varepsilon/3 \right)\\
 \leq \PP\left(\left| \frac{a}{a_N}-1\right|\geq\sqrt{\varepsilon/3} \right)+ \PP\left(\left| \frac{b_N}{b}-1\right|\geq\sqrt{\varepsilon/3} \right)+\PP\left(\left| \frac{a}{a_N}-1\right|\geq\varepsilon/3 \right)+\varphi_2\left(b\varepsilon/3 \right)\\
 = \PP\left(\left| \frac{a}{a_N}-1\right|\geq\sqrt{\varepsilon/3} \right)+\PP\left(\left| \frac{a}{a_N}-1\right|\geq\varepsilon/3 \right)+ \varphi_2\left(b\sqrt{\varepsilon/3} \right)+\varphi_2\left(b\varepsilon/3 \right).
\end{multline}
We next relate 
$$ 
\PP\left(\left| \frac{a}{a_N}-1\right|\geq\sqrt{\varepsilon/3} \right)+\PP\left(\left| \frac{a}{a_N}-1\right|\geq\varepsilon/3 \right),
$$ 
to the function $\varphi_1$. We have 
\begin{multline}
 \PP\left(\left| \frac{a}{a_N}-1\right|\geq\varepsilon/3 \right)= \PP\left(\frac{a}{a_N}-1\geq\varepsilon/3 \right)+\PP\left( \frac{a}{a_N}-1<-\varepsilon/3 \right)\\
 = \PP\left( a-a_N\geq a_N \varepsilon/3 \right)+ \PP\left(a-a_N<-a_N\varepsilon/3 \right)\\
 =  \PP\left( a-a_N\geq \frac{a\varepsilon}{3+\varepsilon} \right)+ \PP\left(a-a_N<-\frac{a\varepsilon}{3-\varepsilon}\right).
\end{multline}
Since $a>0$ and $-a\varepsilon/\{3-\varepsilon\} \leq -a\varepsilon/\{3+\varepsilon\} $ we obtain 
\begin{equation}
 \PP\left(\left| \frac{a}{a_N}-1\right|\geq\varepsilon/3 \right)\leq  \PP\left( a-a_N\geq \frac{a\varepsilon}{3+\varepsilon} \right)+ \PP\left( a-a_N< -\frac{a\varepsilon}{3+\varepsilon} \right)= \varphi_1\left(\frac{a\varepsilon}{3+\varepsilon} \right).
\end{equation}
By similar arguments we obtain 
\begin{equation}
\PP\left(\left| \frac{a}{a_N}-1\right|\geq\sqrt{\varepsilon/3} \right)\leq \varphi_1\left(\frac{a\sqrt \varepsilon}{\sqrt 3+\sqrt\varepsilon} \right).
\end{equation}
We finally deduce that 
\begin{equation}
 \PP\left(\left| \frac{a}{a_N}-1\right|\geq\varepsilon/3 \right)\leq 2\varphi_1\left( \frac{a\varepsilon}{3+\varepsilon}\right)+2\varphi_2(\varepsilon/3).
\end{equation}
Which means that for $u\in (0,1)$ we have 
\begin{equation}
\PP\left(\left|\frac{a}{b}-\frac{a_N}{b_N}\right|\geq \min\left\{1,\frac{1}{c}\right\}u/(2\mathfrak C)\right)\leq 2\varphi_1(r_1 u)+2\varphi_2(r_2 u),
\end{equation}
where, 
$$
r_1=\frac{a\mathfrak C}{2  \min\left\{1,\frac{1}{c}\right\}} \quad \text{and} \quad r_2= \min\left\{1,\frac{1}{c} \right\}\frac{b}{6\mathfrak C}.
$$
It finally remains to control $$
\PP\left(\left|\left(\frac{a}{b}\right)^{c_N-c}-1\right|\geq u\left(\frac{a}{b}\right)^{-c}\right).
$$ 
We have 
\begin{multline}
\PP\left(\left|\left(\frac{a}{b}\right)^{c_N-c}-1\right|\geq u\left(\frac{a}{b}\right)^{-c}\right)\\
= \PP\left(\left(\frac{a}{b}\right)^{c_N-c}\geq 1+u\left(\frac{a}{b}\right)^{-c}\right)+\PP\left(\left(\frac{a}{b}\right)^{c_N-c}\leq 1-u\left(\frac{a}{b}\right)^{-c}\right).
\end{multline} 
Since for $x>0$ we have $\log(x)\leq x-1$, for $u$ small such that  $1 -\left(a/b\right)^{-c}u>0$, we obtain~:
\begin{equation}\label{eq_C11}
\PP\left((c_N-c)\log\frac{a}{b}\leq \log\left(1-u\left(\frac{a}{b}\right)^{-c}\right)\right)\leq \PP\left((c_N-c)\log\frac{a}{b}\leq -\left(\frac{a}{b}\right)^{-c}u\right),
\end{equation}
and since  for $x\in[0,1]$ we have $\log(1+x)\geq x/2$ we obtain~:
\begin{equation}\label{eq_C22}
\PP\left((c_N-c)\log\frac{a}{b}\geq \log\left(1+u\left(\frac{a}{b}\right)^{-c}\right)\right)\leq \PP\left((c_N-c)\log\frac{a}{b}\geq \left(\frac{a}{b}\right)^{-c}u/2\right).
\end{equation}
Combining \eqref{eq_C11} and \eqref{eq_C22}, whatever the sign of $\log(a/b)$ is, we obtain 
\begin{equation}
\PP\left(\left|\left(\frac{a}{b}\right)^{c_N-c}-1\right|\geq u\left(\frac{a}{b}\right)^{-c}\right)\leq \PP\left(\left|c_N-c\right|\geq\frac{u}{2(a/b)^{c}\left|\log(a/b)\right|} \right)= \varphi_3(r_3u),
\end{equation}
where $$
r_3=\min\left\{1/\sqrt{2} , \left(\left(\frac{a}{b}\right)^{c}\left|\log(\frac{a}{b})\right|\right)^{-1}\right\}
$$

\end{proof}

\medskip

\begin{lemSM}\label{concentration:variance}
 Under the conditions of Proposition \ref{prop_def_A} and for sufficiently large $\mathfrak m$,  a positive constant $\mathfrak e$ exists such that, for any  $\eta \in (0,1)$, 
\begin{equation}
\PP\left(|\widehat v(\Tt) -v(\Tt)|\geq \eta \right)\leq 2\exp\left(-\mathfrak e N\eta^2 \right),
\end{equation}
where $\widehat v(\Tt) =\max\{\widetilde v(\Tt), \underline v\}$ with
$$
\widetilde v(\Tt) = \frac{1}{N}\sum_{j=1}^N\{\widetilde X ^{(j)}(\Tt)\}^2.
$$
\end{lemSM}

\begin{proof}[Proof of Lemma SM.\ref{concentration:variance}]
Recall the notation $\xi^{(j)}(\Tt)= \widetilde X^{(j)}- X^{(j)}$ for $1\leq j\leq N.$ We decompose
\begin{equation}
    \widetilde v(\Tt) -v(\Tt)=\frac{1}{N}\sum_{j=1}^N\left(\{\widetilde X ^{(j)}(\Tt)\}^2-\EE[\{\widetilde X ^{(j)}(\Tt)\}^2]\right)+\EE[\widehat v(\Tt)]- \EE[\{X(\Tt)\}^2].
\end{equation}

\textbf{The bias term :}  We have 
\begin{equation}
    \EE\left[\widetilde v(\Tt)-\{X(\Tt)\}^2\right]
        =\frac{1}{N}\sum _{j=1}^N\EE\left[\{\widetilde X(\Tt)-X(\Tt)\}\times \{\widetilde X(\Tt)+X(\Tt)\}\right],
\end{equation}
and by  Cauchy-Schwartz inequality,
\begin{equation}
     \EE\left[\widetilde v(\Tt)-\{X(\Tt)\}^2\right]\leq 
     R_2(\mathfrak m)^{1/2} \times \frac{1}{N}\sum _{j=1}^N \EE\left[\{\widetilde X^{(j)}(\Tt)+X^{(j)}(\Tt)\}^2\right]^{1/2}.
\end{equation}
In addition, we have 
\begin{multline}
\EE\left[\{\widetilde X^{(j)}(\Tt)+X^{(j)}(\Tt)\}^2\right]= \EE\left[\{\xi^{(j)}(\Tt)+2X^{(j)}(\Tt)\}^2\right]
\\ \leq 2 R_2(\mathfrak m)+8\EE[\{X^{(j)}(\Tt)\}^{2}]\leq 2 R_2(\mathfrak m)+8\mathfrak a_1. 
\end{multline}
Therefore,
\begin{equation}
\EE\left[\widehat v(\Tt)-\{X(\Tt)\}^2\right]\leq R_2(\mathfrak m)^{1/2}\left( 2 R_2(\mathfrak m)+8\mathfrak a_1\right)^{1/2}:= \eta_*.
\end{equation}

\textbf{The variance term :} 
\begin{multline}
    \EE\left[\left|\{\widetilde X(\Tt)\}^2-\EE[\{\widetilde X(\Tt)\}^2]\right|^p\right]\leq 2^p\EE\left[\{\widetilde X (\Tt)\} ^{2p}\right]\\
    \leq 2^{3p-1}\left(\EE[\{\xi^{(j)}(\Tt)\}^{2p}]+\EE[\{X^{(j)}(\Tt)\}^{2p}]\right)\\
    \leq 2^{3p-1}\left( \frac{p!}{2} \mathfrak c\mathfrak D ^{p-2}\rho(\mathfrak m)^{2p}+ \frac{p!}{2}\mathfrak a_1\mathfrak A_1^{p-2} \right).
\end{multline}
 Since,  for sufficiently large $\mathfrak m$,  $\rho(\mathfrak m)\leq 1$,  positive constants $\widetilde {\mathfrak {a}}$ and $\widetilde {\mathfrak A}$ exist such that 
\begin{equation}
  \EE\left[\left|\{\widetilde X(\Tt)\}^2-\EE[\{\widetilde X(\Tt)\}^2]\right|^p\right]\leq \frac{p!}{2}\widetilde {\mathfrak{a}}\widetilde{\mathfrak{A}}^{p-2}.
\end{equation}

\textbf{Exponential bound:}
Let $0<2\eta\leq\eta_* $, We have
\begin{multline}
\PP\left(\widetilde v(\Tt) -v(\Tt)\geq \eta \right)\leq \PP\left(\frac{1}{N}\sum_{j=1}^N\left(\{\widetilde X ^{(j)}(\Tt)\}^2-\EE[\{\widetilde X ^{(j)}(\Tt)\}^2]\right)+\EE[\widehat v(\Tt)]- \EE[\{X(\Tt)\}^2]\geq \eta \right)\\
\leq \PP\left(\frac{1}{N}\sum_{j=1}^N\left(\{\widetilde X ^{(j)}(\Tt)\}^2-\EE[\{\widetilde X ^{(j)}(\Tt)\}^2]\right)+\eta_*\geq \eta \right)\\
\leq \PP\left(\frac{1}{N}\sum_{j=1}^N\left(\{\widetilde X ^{(j)}(\Tt)\}^2-\EE[\{\widetilde X ^{(j)}(\Tt)\}^2]\right)\geq \eta/2 \right).
\end{multline} 
Bernstein's inequality then implies 
\begin{equation}
\PP\left(\widetilde v(\Tt) -v(\Tt)\geq \eta \right)\leq \exp\left(-\frac{N\eta^2}{8\widetilde{\mathfrak a}+4\widetilde{\mathfrak A}\eta}\right).
\end{equation}
 For sufficiently large $\mathfrak m$ and $\eta\in (0,1)$,  we have 
\begin{equation}
\PP\left(\widetilde v(\Tt) -v(\Tt)\geq \eta \right)\leq \exp\left(-\mathfrak e N\eta^2 \right),
\end{equation}
for some positive constant $\mathfrak e.$ The same bound is valid for $\PP\left(\widetilde v(\Tt) -v(\Tt)\leq -\eta \right).$
Finally, it suffices to notice that by definition
$$
\PP\left(|\widehat v(\Tt) -v(\Tt)|\geq \eta \right) \leq \PP\left(|\widetilde v(\Tt) -v(\Tt)|\geq \eta \right).
$$
\end{proof}

\medskip

\begin{lemSM}\label{lemme_tec_1}
   Let  $X\in \mathcal H^{H_1,H_2}$ and $\boldsymbol L =(L_1,0,0,L_2)$. For any $\rho\leq \Delta_0$, 
    \begin{equation*}
             \EE\left[\left\{X(\Tt)-X(\Tt+\rho e_i)\right\}^2\right] \leq L_i(\Tt) \rho^{2H_i(\Tt)}\times\{1 + O(\rho\log(\rho))\},\quad \forall\Tt\in\cT,\quad i=1,2.
    \end{equation*}
\end{lemSM}

\begin{proof}[Proof of Lemma SM.\ref{lemme_tec_1}]
    We fix $i=1,2$ and  denote $\tilde \Tt_i =  \Tt+\rho/2 e_i$. Since $X\in \mathcal H^{H_1,H_2}$ and $\boldsymbol L =(L_1,0,0,L_2)$, we have 
    \begin{equation}
     \EE\left[\left\{X(\Tt)-X(\Tt+\rho e_i)\right\}^2\right]=  \EE\left[\left\{X\left(\tilde \Tt_i -\frac{\rho}{2}e_i\right)-X\left(\tilde\Tt_i+\frac{\rho}{2} e_i\right)\right\}^2\right]
     \leq L_i(\tilde \Tt_i)\rho^{2H_i(\tilde\Tt_i)}.
    \end{equation}
    The function $H_i$ is continuously differentiable, and thus 
    \begin{equation}
    \rho^{2H_i(\tilde\Tt_i)}= \rho^{2H_i(\Tt)}\rho^{2H_i(\tilde\Tt_i)-2H_i(\Tt)}
    = \rho^{2H_i(\Tt)}\rho^{2R(\Tt)},
     \end{equation}   
    with
    $$
    R(\Tt) = \nabla H_i(\Tt_*) \cdot (\tilde \Tt_i-\Tt),
    $$
    where $\Tt_*$ is a point on the segment between $\tilde \Tt_i$ and $\Tt$. (Here, for two $2-$dimensional vectors $\Uu$ and $\Vv$, let $\Uu \cdot \Vv$ denote their scalar product.)    
    Since by the assumptions, a constant $C>0$ exists such that  $|R(\Tt)|\leq C \rho$  and 
    $$
    \left|\rho^{2R(\Tt)} - 1\right| = \left|\exp\{ 2R(\Tt)\log(\rho) \}- 1\right| \leq C 2R(\Tt)\log(\rho),
    $$
    we deduce 
    $$
    \rho^{2H_i(\tilde\Tt_i)}
      = \rho^{2H_i(\Tt)}\times \{1+O\left(\rho\log(\rho)\right)\}.
$$
    On the hand, the function $L_i$ is Lipschitz continuous, and thus 
    \begin{align}
    \EE\left[\left\{X(\Tt)-X(\Tt+\rho e_i)\right\}^2\right] &\leq L_i(\tilde \Tt_i) \rho^{2H_i(\Tt)}\times \{1+O\left(\rho\log(\rho)\right)\}\\
    &\leq \left (L_i(\Tt) \rho^{2H_i(\Tt)}+ |L_i(\tilde \Tt_i)-L_i(\Tt)| \rho^{2H_i(\Tt)}\right)\times \{1+O\left(\rho\log(\rho)\right)\}\\
    &=\left (L_i(\Tt) \rho^{2H_i(\Tt)}+O(\rho^{2H_i(\Tt)+1})\right)\times \{1+O\left(\rho\log(\rho)\right)\}\\
    &= L_i(\Tt) \rho^{2H_i(\Tt)} + O(\rho^{2H_i(\Tt)+1}\log(\rho))
    \end{align}
    We finally obtain 
    \begin{equation}
     \EE\left[\left\{X(\Tt)-X(\Tt+\rho e_i)\right\}^2\right] \leq L_i(\Tt) \rho^{2H_i(\Tt)}\times\{1 + O(\rho\log(\rho))\}.
    \end{equation}
\end{proof}


\begin{lemSM}\label{Ber_par}
 We have that 
 $$
 p(\Tt)\geq \pi c h_1h_2.
 $$
 \end{lemSM}
\begin{proof}[Proof of Lemma SM.\ref{Ber_par}]
We have 
\begin{equation}
p(\Tt)=\PP\left(\mathbf{B}(\Tt^{[0]}_m-\Tt)\in B(0,1)\big{|}M_0\right)
=\int_{\mathbf{B}(\Uu-\Tt)\in B(0,1)}f_{\TT}(\Uu)\D \Uu.
\end{equation}
By substitution $(\Ss=\mathbf{B}(\Uu-\Tt))$,  using \ref{LP2} and the Lebesgue measure of the unit ball, we obtain
$$
p(\Tt)= \int_{\Ss \in B(0,1)}f_{\TT}(\mathbf{B}^{-1}\Ss +\Tt)|\mathbf{B}^{-1}| \D \Ss \geq c \pi h_1h_2.
$$
\end{proof}

\bibliographystyle{apalike}

\bibliography{biblio_final.bib}